\documentclass[11pt]{article}
\usepackage{amsmath,amssymb}
\usepackage{pstricks}
\usepackage{pst-node}
\usepackage{pst-poly}
\usepackage{multido}
\usepackage{tikz}

\usetikzlibrary{snakes}

\newtheorem{result}{Theorem}

\newtheorem{define}{Definition}
\newtheorem{support}{Lemma}
\newtheorem{propo}{Proposition}

\newtheorem{note}{Remark}

\renewcommand*{\thetheorem}{\Alph{theorem}}

\newcommand{\qed}{%
\ifmmode 
\else \leavevmode\unskip\penalty9999 \hbox{}\nobreak\hfill \fi
\quad\hbox{\qedsymbol}}
\newcommand{\openbox}{\leavevmode \hbox to.77778em{%
\hfil\vrule
\vbox to.675em{\hrule width.6em\vfil\hrule}%
\vrule\hfil}}
\newcommand{\qedsymbol}{\openbox}

\newcommand{\showgrid}{}
\newcommand{\gridon}{\renewcommand{\showgrid}{\psset{subgriddiv=1,griddots=10,gridlabels=6pt}\psgrid}}
\gridon

\begin{document}

\begin{center} {\bf\LARGE  On Convex Graphs Having Plane Spanning Subgraph of Certain Type} \end{center}

\vskip5pt

\centerline{\large   Niran Abbas Ali$^a$, Gek L. Chia$^{b}$,  \  Hazim Michman  Trao$^{c}$ \ and \   Adem Kilicman$^d$ }

\begin{center}
\itshape\small  $^{a, c, d}\/$Department of Mathematics, \\ Universiti Putra Malaysia, 43400 Serdang, Malaysia,  \\
 \vspace{1mm}
 $^{b}\/$Department of Mathematical and Actuarial Sciences, \\  Universiti Tunku Abdul Rahman, Sungai Long Campus,    Malaysia \\
 \vspace{1mm}
 $^{a, c}\/$Department of Mathematics, College of Science, \\ Al Mustansiriyah University, 10052 Filastin Street, Iraq
\end{center}

\begin{abstract}
Motivated by  a result of \cite{alta:refer}, we determine necessary and sufficient conditions on $F\/$ with $|E(F)| \leq n-1\/$ for which  $K_n - F\/$ admits a $g$-angulation. For $|E(F)| \geq n\/$,  we investigate the possibility of placing $F\/$ in $K_n\/$ such that $K_n -F\/$ admits a $g$-angulation for certain families of graphs $F\/$.
\end{abstract}

\vspace{1mm}
 \section{Introduction and Preliminary}

By a {\em geometric graph\/} we mean a graph whose edges are straight line segments. By a {\em convex graph\/}, we mean a geometric graph whose vertices are in convex position. 

\vspace{5mm}  Let $S\/$ be a set of $n\/$ points in general position in the plane. A {\em   $g$-angulation\/ } of $S\/$ is a plane graph in which each face interior to the convex hull of $S\/$ is a $g$-cycle.
 A {\em   convex $g$-angulation\/ } is a $g$-angulation on $S\/$ of $n\/$ points in convex position in the plane.
 We say that $G_n$ is a $g$-angulation of a graph $G(V, E)$ if $E(G_n)\subseteq E$.
 In particular, the $g$-angulation is a triangulation if $g=3$.

\vspace{5mm} The triangulation existence problem is the following: On a given graph $G$, decide whether there exists a triangulation of $G$.
This problem is NP-complete  (see \cite{llo:refer}).  \
This article extends the problem of an article \cite{alta:refer} on triangulability of convex geometric graphs to the  $g$-angulation existence problem of a convex geometric graph $G\/$  by considering a spanning subgraph $F_n$  of $K_n$ with $G=K_n-F_n$ is the convex graph obtained from $K_n\/$ by deleting the set of  edges of $F_n\/$.

\vspace{5mm} To decide whether $G\/$ admits a $g$-angulation or not, we first characterize the forbidding configurations for any possible $g$-angulation of $K_n$.
A configuration of $n$ vertices in convex position (vertices of $K_n$) having a common edge with each $g$-angulation of $K_n$, is a forbidding configuration for $g$-angulations of $K_n$.

\vspace{5mm}  Determining the smallest size of the forbidding configurations for $g$-angulations of $K_n$ is a natural Tur$\acute{a}$n-type question, as it is equivalent to determining the largest size of a convex geometric graph does not contain a $g$-angulation.   
We answer for the   ``forbidding configurations for $g$-angulations" question not only by determining their smallest size which is $n-g+1\/$, but rather giving a complete characterization of all those forbidding configurations of size at least $n-g+1\/$ and at most $n-1\/$.

\vspace{5mm}  We present the characterizations ${\cal F}_{n,g}(^*)\/$   and  ${\cal J}_{n,g}(^*_\beta)\/$ ($\beta\in \{1, 2, \ldots, 2g-3\}\/$) that forbidding $g\/$-angulations of $K_n\/$.

\vspace{5mm}  By Theorem \ref{fstarsg}, Propositions \ref{Jstar1}, \ref{Jstar2}, and \ref{Jstar-gamma}, we show that $K_n-{\cal J}_{n,g}(^*)\/$ and  $K_n-{\cal J}_{n,g}(^*_\beta)\/$ admits no $g\/$-angulation respectively.

\vspace{5mm}  If $|E(F_n)|\leq n-g+1\/$, we show by Theorem  \ref{n-2g} that $K_n-F_n\/$ admits a $g$-angulation if and only if $F_n\neq {\cal F}_{n,g}(^*)\/$.
If $n-g+2 \leq |E(F_n)|\leq n-1\/$, we show by Theorem  \ref{n-g+1+mu} that $K_n-F_n\/$ admits a $g\/$-angulation if and only if $F_n\neq {\cal J}_{n,g}(^*_\beta)\/$.

\vspace{5mm}
 For the case where $F_n\/$ has at least $n\/$ edges, it seems difficult to obtain a characterization on $F_n\/$ such that  $K_n -F_n\/$ admits a triangulation.
 For this we confine  our attention to seek for the possibility of arranging certain families of graphs $F_n\/$ as a convex geometric graph
  such that $K_n - F_n\/$ admits a $g$-angulation.
 If such a configuration exists for $F_n\/$, then we say that $F_n\/$ is {\em potentially $g$-angulable in $K_n\/$}.
 Potentially $g$-angulable graphs are considered in Section \ref{pt} where we
 (i) determine precisely the value of $n\/$ for which the $n\/$-cycle is potentially $g$-angulable in $K_n\/$ (Theorem \ref{fnconnected}), \ and
 (ii) characterize all $2\/$-regular graphs  which are potentially $g$-angulable in $K_n\/$ (Theorem \ref{disconnected}).

\vspace{5mm}
The potentially $g$-angulable  problem is extended to the regular case in Section \ref{regular} where we characterize all $3\/$-regular graphs which are potentially $4$-angulable in $K_n\/$ (Theorem \ref{4-angulation}).

\vspace{5mm}

Throughout, we shall adopt the following notations.  Unless otherwise stated, the vertices of a convex complete graph $K_n\/$ will be denoted by $v_0, v_1, \ldots, v_{n-1}\/$ in cyclic clockwise ordered. Also, unless otherwise specified, any operation on the subscript of $v_i\/$ is reduced modulo $n\/$.

 \vspace{1mm}
\begin{support} \label{earg}
Suppose $g \geq 3\/$ and $t \geq 2\/$ are natural numbers and assume that $n=g+t(g-2)\/$.   Let $F\/$ be a subgraph of a convex complete graph $K_n\/$. Assume that $F\/$ has  at most $n-1\/$ edges and having no boundary edge of $K_n\/$. Then $K_n -F\/$ has an edge of the form  $v_{j}v_{j+g-1}\/$.
\end{support}

\vspace{1mm}  \noindent
{\bf Proof:}   If  the lemma  is not true, then it implies that $v_{j-g+1}v_j, v_jv_{j+g-1}  \in E(F)\/$, and recursively, this implies that $F\/$ is a spanning subgraph of $K_n\/$ with minimum vertex-degree at least $2\/$. But this implies that $|E(F)| \geq n\/$, a contradiction.     \qed

\section{$g$-angulable graphs}     \label{ptg}

\begin{result}{\bf:} \label{atmostn-g}  Suppose  $g \geq 3\/$ is a natural number and let  $n=g+t(g-2)\/$ where  $t\geq 0$ is any natural number. Let  $F\/$ be a subgraph of the complete convex graph $K_n$. Suppose $F\/$ contains no boundary edge of $K_n\/$ and  $|E(F)| \leq n-g\/$. Then  $K_n-F\/$ admits   a  $g$-angulation. 
\end{result}

\vspace{5mm} \noindent
{\bf Proof:} We prove this by induction on $t\/$.

\vspace{5mm}
The result is trivially true if $t =0\/$. If  $t=1\/$, then $n =2(g-1)\/$ and $F\/$ has at most $g-2\/$ edges. Hence $K_n - F\/$ has a vertex, say $v_i\/$  which is  adjacent to every other vertex in $K_n - F\/$ and this means that $v_i v_{i +g-1}\/$ together with the boundary edges of $K_n\/$ form a $g\/$-angulation of $K_n - F\/$.

\vspace{5mm}
Assume that $t\geq 2\/$ and the result is true for all natural numbers $t'\/$ where $t' < t\/$. Since $g \geq 3\/$, by Lemma \ref{earg}, $K_n - F\/$ contains an edge of the form $v_iv_{i+g-1}\/$. By relabeling the vertices of $K_n - F\/$, if necessary we may assume that $i=0\/$.

\vspace{5mm}
If, for some $ j \in \{1, 2, \ldots, g-2\}\/$, $v_j\/$ is adjacent to every other vertex in $K_n - F\/$, then $v_jv_{j+kg-2k+1}\/$ where $ k =1, 2, \ldots, t\/$ together with the boundary  edges of $K_n\/$ form a $g\/$-angulation of $K_n -F\/$.

\vspace{5mm}
Hence we assume that, for each $j \in \{1, 2, \ldots, g-2\}\/$, $v_j\/$ is incident with at least one edge of $F\/$. Then the subgraph obtained from $K_n - F\/$ by deleting all the vertices $v_1, v_2, \ldots, v_{g-2}\/$ is a convex graph of the form $K_m - F'\/$ where $m = g +(t-1)(g-2)\/$ and $F' = F - \{v_1, v_2, \ldots, v_{g-2}\}\/$. Moreover, $F'\/$ has at most $m-g\/$ edges (since there are at least $g-2\/$ edges of $F\/$ incident to the vertices $v_1, v_2, \ldots, v_{g-2}\/$).

\vspace{5mm} By induction $K_m - F'\/$ admits a $g\/$-angulation $G'\/$. As such, $G' \cup \{v_0v_1v_2 \cdots v_{g-1} \}\/$ is a $g\/$-angulation for $K_n - F\/$.
This completes the proof.  $\qed$

\begin{note}{\bf:}
The result in Theorem \ref{atmostn-g} is tight with respect to $|E(F)|\/$ since deleting $n-g+1\/$ edges from $K_n\/$ does not guarantee that the resulting graph admits a $g\/$-angulation. This follows directly from the next result in the next section.
\end{note}

\section{Graphs with at most $n-g+1\/$ edges}  \label{n-1oelessg}

This section presents a characterization of  ${\cal F}_{n,g}(^*)\/$, that shares any possible $g$-angulation of $K_n$ with at least one edge.
We show that $G\/$  admits no $g$-angulation when $G=K_n-{\cal F}_{n,g}(^*)\/$, $G\/$ admits a $g$-angulation if and only if $F_n\neq {\cal F}_{n,g}(^*)\/$.

\vspace{1mm} Throughout, we let $d_G(v)\/$ denote the degree of $v\/$ in the graph $G\/$.

 \vspace{2mm}
\begin{define}   \label{fng}
 Suppose  $g \geq 3\/$ is a natural number and let  $n=g+t(g-2)\/$ where  $t\geq 1$ is any natural number.
Let ${\cal F}_{n,g}(^*)\/$ denote a convex geometric graph with no isolated vertices and having  $n-g+1\/$  edges  such that

\vspace{1mm} (i) for any $v_iv_j\/$ in  ${\cal F}_{n,g}(^*)\/$,  $|j-i| \equiv 1 \  ({\rm mod}\ (g-2))\/$,

\vspace{1mm} (ii) whenever $d(v_i) \geq 2\/$, then $v_j v_{j-g+1}\/$ is an edge of ${\cal F}_{n,g}(^*)\/$ for $j=i+1, \ldots , i+g-2\/$,

\vspace{1mm} (iii) the neighbor $v_i\/$ of any pendant vertex satisfies the condition that $v_j v_{j-g+1}\/$ is an edge of ${\cal F}_{n,g}(^*)\/$ for $j=i+1, \ldots , i+g-2\/$,

\vspace{1mm} (iv) for any two pendant vertices $v_r, v_s\/$ such that whenever $v_rv_i\/$ and $v_sv_j\/$ crosse each other in ${\cal F}_{n,g}(^*)\/$, then $|i-j|\leq g-2\/$.

\end{define}

\vspace{6mm} It is easy to see that in the definition,  conditions (ii) and (iii) imply that  all the pendant vertices of  ${\cal F}_{n,g}(^*)\/$ are in consecutive order. As such we may label the vertices of ${\cal F}_{n,g}(^*)\/$ so that $d(v_i) =1 \/$ if and only if $i = k, k+1, \ldots, n-1\/$ where $k\/$ is a  natural number with $g-2 \leq k \leq n-2(g-1)\/$ if $t \geq 2\/$, and $k=0\/$ if $t=1\/$.

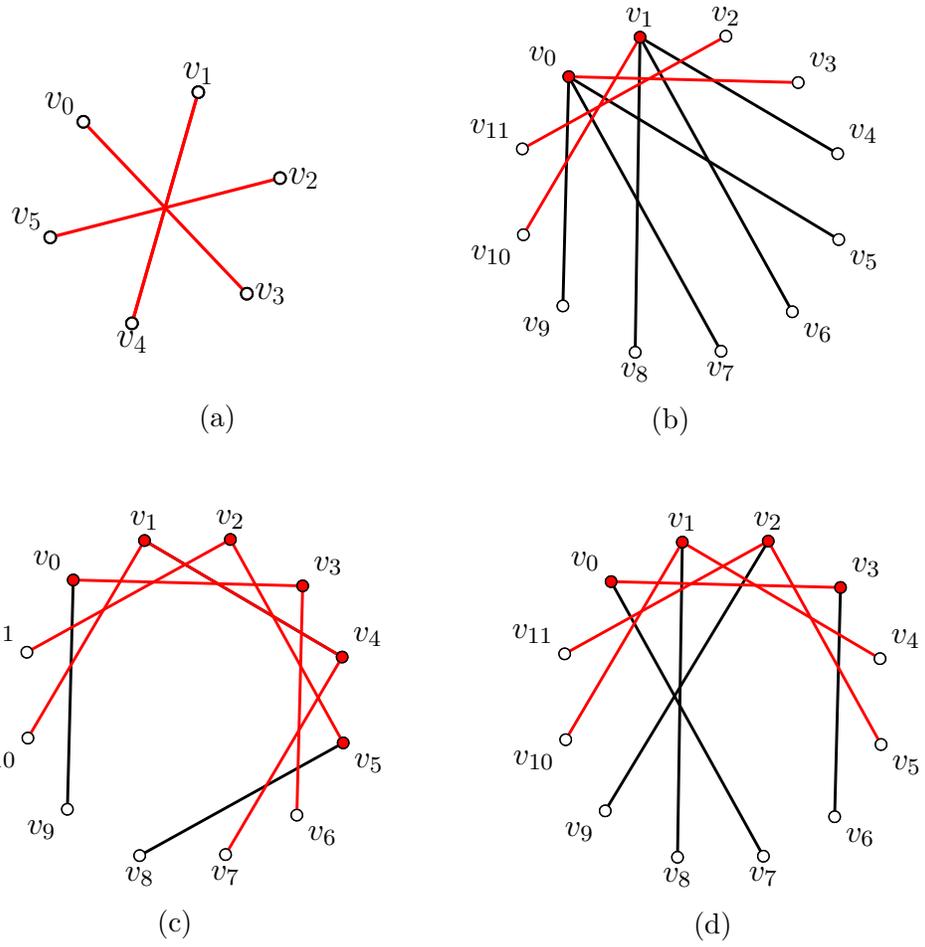
\begin{figure}[htp]
\begin{center}
\begin{minipage}{.4\textwidth}
\resizebox{4.5cm}{!}{
\begin{tikzpicture}[rotate=-16,.style={draw}]
\coordinate (center) at (0,0);

   \def\radius{2cm}
   \foreach \x in {0, 30,...,360} {
               }


      \coordinate (v0) at (-1.7,1.);\filldraw[black] (v0) circle(3pt); \node[above left] at (v0) {\LARGE${v_0}$};

      \coordinate (v1) at (0, 2);\filldraw[black] (v1) circle(3pt); \node[above] at (v1) {\LARGE${v_1}$};
       \coordinate (v2) at (1.7, 1);\filldraw[black] (v2) circle(3pt);\node[right] at (v2) {\LARGE${v_2}$};

       \coordinate (v3) at (1.7, -1);\filldraw[black] (v3) circle(3pt);\node[ right] at (v3) {\LARGE${v_3}$};
       \coordinate (v4) at (0, -2);\filldraw[black] (v4) circle(3pt);\node[below] at (v4) {\LARGE${v_4}$};

       \coordinate (v5) at (-1.7, -1);\filldraw[black] (v5) circle(3pt);\node[above left] at (v5) {\LARGE${v_5}$};

\filldraw[white] (0,3.5) circle(3pt);

     \draw [line width=1.5,red](v4)-- (v1);
     \draw [line width=1.5,red](v5)-- (v2);
     \draw [line width=1.5,red](v0) -- (v3);
     \draw [line width=1.5,red](v1) -- (v4);

\filldraw[black] (v0) circle(3pt);
\filldraw[black] (v1) circle(3pt);
\filldraw[black] (v2) circle(3pt);
\filldraw[black] (v3) circle(3pt);
\filldraw[black] (v4) circle(3pt);
\filldraw[black] (v5) circle(3pt);

 \filldraw[white] (v0) circle(2pt);
 \filldraw[white] (v1) circle(2pt);
 \filldraw[white] (v2) circle(2pt);
 \filldraw[white] (v3) circle(2pt);
 \filldraw[white] (v4) circle(2pt);
 \filldraw[white] (v5) circle(2pt);

\end{tikzpicture}
}

\vspace{4mm}

\centering (a)
\end{minipage}
\begin{minipage}{.4\textwidth}
\resizebox{6.cm}{!}{
\begin{tikzpicture}[rotate=-16,.style={draw}]
\coordinate (center) at (0,0);
   \def\radius{2cm}
   \foreach \x in {0, 30,...,360} {
             \filldraw[] (\x:2cm) circle(1pt);
               }


      \coordinate (v0) at (-1.7,1.);\filldraw[black] (v0) circle(2pt); \node[above left] at (v0) {${v_0}$};

      \coordinate (v1) at (-1., 1.7);\filldraw[black] (v1) circle(2pt); \node[above] at (v1) {${v_1}$};
       \coordinate (v2) at (0, 2);\filldraw[black] (v2) circle(2pt);\node[above] at (v2) {${v_2}$};

       \coordinate (v3) at (1., 1.7);\filldraw[black] (v3) circle(2pt);\node[above right] at (v3) {${v_3}$};
       \coordinate (v4) at (1.7, 1);\filldraw[black] (v4) circle(2pt);\node[above right] at (v4) {${v_4}$};

       \coordinate (v5) at (2, 0.);\filldraw[black] (v5) circle(2pt);\node[below right] at (v5) {${v_5}$};
       \coordinate (v6) at (1.7, -1);\filldraw[black] (v6) circle(2pt);\node[below right] at (v6) {${v_6}$};

       \coordinate (v7) at (1., -1.7);\filldraw[black] (v7) circle(2pt);\node[below] at (v7) {${v_7}$};
      \coordinate (v8) at (0, -2);\filldraw[black] (v8) circle(2pt);\node[below] at (v8) {${v_8}$};

       \coordinate (v9) at (-1., -1.7);\filldraw[black] (v9) circle(2pt);\node[below left] at (v9) {${v_{9}}$};
       \coordinate (v10) at (-1.7, -1);\filldraw[black] (v10) circle(2pt);\node[below left] at (v10) {${v_{10}}$};

       \coordinate (v11) at (-2, 0);\filldraw[black] (v11) circle(2pt);\node[above left] at (v11) {${v_{11}}$};

     \draw [line width=1,black](v4)-- (v1);
     \draw [line width=1,black](v5)-- (v0);
     \draw [line width=1,black](v6)-- (v1);
     \draw [line width=1,black](v7)-- (v0);
     \draw [line width=1,black](v8)-- (v1);
     \draw [line width=1,black](v9)-- (v0);

     \draw [line width=1,red](v10)-- (v1);
     \draw [line width=1,red](v11)-- (v2);
     \draw [line width=1,red](v0) -- (v3);

\filldraw[black] (v0) circle(2pt);
\filldraw[black] (v1) circle(2pt);
\filldraw[black] (v2) circle(2pt);
\filldraw[black] (v3) circle(2pt);
\filldraw[black] (v4) circle(2pt);
\filldraw[black] (v5) circle(2pt);
\filldraw[black] (v6) circle(2pt);
\filldraw[black] (v7) circle(2pt);
\filldraw[black] (v8) circle(2pt);
\filldraw[black] (v9) circle(2pt);
\filldraw[black] (v10) circle(2pt);
\filldraw[black] (v11) circle(2pt);

 \filldraw[red] (v0) circle(1.5pt);
 \filldraw[red] (v1) circle(1.5pt);
 \filldraw[white] (v2) circle(1.5pt);
 \filldraw[white] (v3) circle(1.5pt);
 \filldraw[white] (v4) circle(1.5pt);
 \filldraw[white] (v5) circle(1.5pt);
 \filldraw[white] (v6) circle(1.5pt);
 \filldraw[white] (v7) circle(1.5pt);
 \filldraw[white] (v8) circle(1.5pt);
 \filldraw[white] (v9) circle(1.5pt);
 \filldraw[white] (v10) circle(1.5pt);
 \filldraw[white] (v11) circle(1.5pt);

\end{tikzpicture}
}

\centering   (b)
\end{minipage}

\vspace{8mm}

\begin{minipage}{.4\textwidth}
\resizebox{6cm}{!}{
\begin{tikzpicture}[rotate=-16,.style={draw}]
\coordinate (center) at (0,0);
   \def\radius{2cm}
   \foreach \x in {0, 30,...,360} {
             \filldraw[] (\x:2cm) circle(1pt);
               }


      \coordinate (v0) at (-1.7,1.);\filldraw[black] (v0) circle(2pt); \node[above left] at (v0) {${v_0}$};

      \coordinate (v1) at (-1., 1.7);\filldraw[black] (v1) circle(2pt); \node[above] at (v1) {${v_1}$};
       \coordinate (v2) at (0, 2);\filldraw[black] (v2) circle(2pt);\node[above] at (v2) {${v_2}$};

       \coordinate (v3) at (1., 1.7);\filldraw[black] (v3) circle(2pt);\node[above right] at (v3) {${v_3}$};
       \coordinate (v4) at (1.7, 1);\filldraw[black] (v4) circle(2pt);\node[above right] at (v4) {${v_4}$};

       \coordinate (v5) at (2, 0.);\filldraw[black] (v5) circle(2pt);\node[below right] at (v5) {${v_5}$};
       \coordinate (v6) at (1.7, -1);\filldraw[black] (v6) circle(2pt);\node[below right] at (v6) {${v_6}$};

       \coordinate (v7) at (1., -1.7);\filldraw[black] (v7) circle(2pt);\node[below] at (v7) {${v_7}$};
      \coordinate (v8) at (0, -2);\filldraw[black] (v8) circle(2pt);\node[below] at (v8) {${v_8}$};

       \coordinate (v9) at (-1., -1.7);\filldraw[black] (v9) circle(2pt);\node[below left] at (v9) {${v_{9}}$};
       \coordinate (v10) at (-1.7, -1);\filldraw[black] (v10) circle(2pt);\node[below left] at (v10) {${v_{10}}$};

       \coordinate (v11) at (-2, 0);\filldraw[black] (v11) circle(2pt);\node[above left] at (v11) {${v_{11}}$};

     \draw [line width=1,black](v4)-- (v1);
     \draw [line width=1,black](v5)-- (v8);
     \draw [line width=1,black](v0)-- (v9);

     \draw [line width=1,red](v10)-- (v1);
     \draw [line width=1,red](v11)-- (v2);
     \draw [line width=1,red](v0) -- (v3);
     \draw [line width=1,red](v1) -- (v4);
     \draw [line width=1,red](v2) -- (v5);
     \draw [line width=1,red](v3) -- (v6);
     \draw [line width=1,red](v4) -- (v7);

\filldraw[black] (v0) circle(2pt);
\filldraw[black] (v1) circle(2pt);
\filldraw[black] (v2) circle(2pt);
\filldraw[black] (v3) circle(2pt);
\filldraw[black] (v4) circle(2pt);
\filldraw[black] (v5) circle(2pt);
\filldraw[black] (v6) circle(2pt);
\filldraw[black] (v7) circle(2pt);
\filldraw[black] (v8) circle(2pt);
\filldraw[black] (v9) circle(2pt);
\filldraw[black] (v10) circle(2pt);
\filldraw[black] (v11) circle(2pt);

 \filldraw[red] (v0) circle(1.5pt);
 \filldraw[red] (v1) circle(1.5pt);
 \filldraw[red] (v2) circle(1.5pt);
 \filldraw[red] (v3) circle(1.5pt);
 \filldraw[red] (v4) circle(1.5pt);
 \filldraw[red] (v5) circle(1.5pt);
 \filldraw[white] (v6) circle(1.5pt);
 \filldraw[white] (v7) circle(1.5pt);
 \filldraw[white] (v8) circle(1.5pt);
 \filldraw[white] (v9) circle(1.5pt);
 \filldraw[white] (v10) circle(1.5pt);
 \filldraw[white] (v11) circle(1.5pt);

\end{tikzpicture}

}

\centering (c)
\end{minipage}
\hspace{10mm}
\begin{minipage}{.4\textwidth}
\begin{center}
\resizebox{6cm}{!}{
\begin{tikzpicture}[rotate=-16,.style={draw}]
\coordinate (center) at (0,0);
   \def\radius{2cm}
   \foreach \x in {0, 30,...,360} {
             \filldraw[] (\x:2cm) circle(1pt);
               }


      \coordinate (v0) at (-1.7,1.);\filldraw[black] (v0) circle(2pt); \node[above left] at (v0) {${v_0}$};

      \coordinate (v1) at (-1., 1.7);\filldraw[black] (v1) circle(2pt); \node[above] at (v1) {${v_1}$};
       \coordinate (v2) at (0, 2);\filldraw[black] (v2) circle(2pt);\node[above] at (v2) {${v_2}$};

       \coordinate (v3) at (1., 1.7);\filldraw[black] (v3) circle(2pt);\node[above right] at (v3) {${v_3}$};
       \coordinate (v4) at (1.7, 1);\filldraw[black] (v4) circle(2pt);\node[above right] at (v4) {${v_4}$};

       \coordinate (v5) at (2, 0.);\filldraw[black] (v5) circle(2pt);\node[below right] at (v5) {${v_5}$};
       \coordinate (v6) at (1.7, -1);\filldraw[black] (v6) circle(2pt);\node[below right] at (v6) {${v_6}$};

       \coordinate (v7) at (1., -1.7);\filldraw[black] (v7) circle(2pt);\node[below] at (v7) {${v_7}$};
      \coordinate (v8) at (0, -2);\filldraw[black] (v8) circle(2pt);\node[below] at (v8) {${v_8}$};

       \coordinate (v9) at (-1., -1.7);\filldraw[black] (v9) circle(2pt);\node[below left] at (v9) {${v_{9}}$};
       \coordinate (v10) at (-1.7, -1);\filldraw[black] (v10) circle(2pt);\node[below left] at (v10) {${v_{10}}$};

       \coordinate (v11) at (-2, 0);\filldraw[black] (v11) circle(2pt);\node[above left] at (v11) {${v_{11}}$};

     \draw [line width=1,black](v3)-- (v6);
     \draw [line width=1,black](v2)-- (v9);
     \draw [line width=1,black](v1)-- (v8);
     \draw [line width=1,black](v0)-- (v7);

     \draw [line width=1,red](v10)-- (v1);
     \draw [line width=1,red](v11)-- (v2);
     \draw [line width=1,red](v0) -- (v3);
     \draw [line width=1,red](v1) --(v4);
     \draw [line width=1,red](v2) --(v5);

\filldraw[black] (v0) circle(2pt);
\filldraw[black] (v1) circle(2pt);
\filldraw[black] (v2) circle(2pt);
\filldraw[black] (v3) circle(2pt);
\filldraw[black] (v4) circle(2pt);
\filldraw[black] (v5) circle(2pt);
\filldraw[black] (v6) circle(2pt);
\filldraw[black] (v7) circle(2pt);
\filldraw[black] (v8) circle(2pt);
\filldraw[black] (v9) circle(2pt);
\filldraw[black] (v10) circle(2pt);
\filldraw[black] (v11) circle(2pt);

 \filldraw[red] (v0) circle(1.5pt);
 \filldraw[red] (v1) circle(1.5pt);
 \filldraw[red] (v2) circle(1.5pt);
 \filldraw[red] (v3) circle(1.5pt);
 \filldraw[white] (v4) circle(1.5pt);
 \filldraw[white] (v5) circle(1.5pt);
 \filldraw[white] (v6) circle(1.5pt);
 \filldraw[white] (v7) circle(1.5pt);
 \filldraw[white] (v8) circle(1.5pt);
 \filldraw[white] (v9) circle(1.5pt);
 \filldraw[white] (v10) circle(1.5pt);
 \filldraw[white] (v11) circle(1.5pt);

\end{tikzpicture}
}

\centering (d)

\end{center}
\end{minipage}
\caption{${\cal F}_{n,4}(^*)$}  \label{condition (iv)}
\end{center}
\end{figure}

\vspace{5mm} Examples of a geometric graph ${\cal F}_{n,g}(^*)\/$  are depicted in Figure\ref{condition (iv)} (a), (b), (c) and (d).
 In Figure\ref{condition (iv)}(a), ${\cal F}_{6,4}(^*)$ with $t=1\/$, in (b) ${\cal F}_{12,4}(^*)\/$ with minimum value of $k=g-2=2\/$, in (c)  ${\cal F}_{12,4}(^*)$ with maximum value of $k=n-2(g-2)=6$, and in (d) ${\cal F}_{12,4}(^*)$ satisfies condition (iv) of Definition\ref{fng}.

\vspace{1mm}
\begin{result} \label{fstarsg}  Suppose  $g \geq 3\/$ is a natural number and let  $n=g+t(g-2)\/$ where  $t\geq 1$ is any natural number.
   Then   $K_n - {\cal F}_{n,g}(^*)\/$ admits no $g$-angulation.

\end{result}

\vspace{1mm}  \noindent
{\bf Proof:}

We first observe that, when $t=1\/$, we have $n=2g-2\/$ and since ${\cal F}_{n,g}(^*)\/$   consists of $g-1\/$ edges of the form  $ v_jv_{j+g-1}\/$,  $j= 0, 1, \ldots, g-2\/$,   $K_n - {\cal F}_{n,g}(^*)\/$ admits no $g$-angulation.

\vspace{5mm} Hence assume that  $t \geq 2\/$.

\vspace{5mm}
Suppose $K_n - {\cal F}_{n,g}(^*)\/$ admit a $g\/$-angulation $G_n\/$.
Note that any  $g$-angulation on a set of points in convex position has a diagonal $v_qv_{q+g-1}\/$.

\vspace{5mm} Let $ F_m = {\cal F}_{n,g}(^*) - \{v_{q+1}, \ldots, v_{q+g-2}\}\/$ where $m = n-g+2\/$.

\vspace{5mm} By the definition of ${\cal F}_{n,g}(^*)\/$, we see that $d_{{\cal F}_{n,g}(^*)}(v_i)=1\/$  for each $i=q+1, q+2, \ldots, q+g-2$.
It is readily checked that  $F_m\/$ is of the form ${\cal F}_{m,g}(^*)\/$ where $m = g + (t-1)(g-2)\/$.  By induction, $K_m-F_m\/$ admits no  $g$-angulation, which is a contradiction with   $G_n$ is a  $g$-angulation for $K_n - {\cal  F}_{n,g}(^*)\/$.

\vspace{5mm} This completes the proof.  $\qed$

\vspace{5mm}
\begin{result}{\bf:}    \label{n-2g}
Suppose  $g \geq 3\/$ is a natural number and let  $n=g+t(g-2)\/$ where  $t\geq 1$ is any natural number. Suppose  $F_n\/$ is a subgraph of the  convex complete graph $K_n\/$ such that  $|E(F_n)| \leq n-g+1\/$  and $F_n\/$ contains no boundary edges  of $K_n\/$. Then $K_n - F_n\/$ admits a $g$-angulation unless
$F_n = {\cal F}_{n,g}(^*)\/$.

\end{result}

\vspace{5mm}  \noindent
{\bf Proof:}  In view of Theorem \ref{fstarsg},  we assume that $F_n \neq {\cal F}_{n,g}(^*)\/$.  Also, we may assume that $F_n\/$ contains no isolated vertices (otherwise   $K_n - F_n\/$ admits a $g$-angulation).

\vspace{5mm}
We prove the result by induction on $t\/$. The result is clearly true  for $t=1\/ $. Assume that $n=g+t(g-2)\/$ with $t\geq2$ and the result is true for all convex graphs $K_m - F_m\/$ where    $m= g+t'(g-2)\/$ with $t'\leq t-1$.

\vspace{5mm} By Lemma \ref{earg}, $K_n - F_n\/$ contains an edge of the form $v_{q}v_{q+g-1}\/$ which is not an edge of $F_n\/$.

\vspace{5mm}
Let $H\/$ be the convex graph obtained from $K_n - F_n\/$ by deleting the vertices  $v_{q+1}, \ldots, v_{q+g-2}$. Clearly, $H\/$ is of the form $K_m - F_m\/$ where $m = g + (t-1)(g-2)\/$.

\vspace{5mm} If $F_m \neq {\cal F}_{m,g}(^*)\/$, then  $K_m - F_m\/$ admits a $g$-angulation (by induction) and this implies that $K_n - F_n\/$ admits a $g$-angulation.
Hence assume that $F_m = {\cal F}_{m,g}(^*)\/$. As such, the  vertices  $v_{q+1}, \ldots, v_{q+g-2}$ are pendant vertices in $F_n\/$  (otherwise $|E(F_m)| \leq m-g\/$, a contradiction with assumption $F_m = {\cal F}_{m,g}(^*)\/$).

\vspace{5mm} There are two cases to consider. Let $S= \{q+1, \ldots, q+g-2\}\/$.

\vspace{5mm} {\em Case (1):     $m=2g-2\/$.}

\vspace{5mm} For each $i \in S\/$, $v_i\/$ is adjacent  either (i) to  $v_{i+g-1}\/$ or else (ii) to $v_{i-g+1}\/$.

\vspace{5mm} If  there exist $i, j \in S\/$ such that

\vspace{5mm}  \hspace{15mm}  $i < j\/$ and   $v_{i}\/$ is adjacent to  $v_{i+g-1}\/$ and $v_{j}\/$ is adjacent to  $v_{j-g+1}\/$  \hfill $(\ast)\/$

 \vspace{5mm} \noindent  then $v_{i} v_{i-g+1}\/$ and $v_{j} v_{j+g-1}\/$ are the diagonals of a required $g\/$-angulation.

 \vspace{5mm} Hence assume that no $i, j \in S\/$ satisfy the condition $(\ast)\/$.
 Let $i\/$ be the largest  integer in $S\/$ such that $v_i\/$ is adjacent to $v_{i-g+1}\/$. This implies that $v_j\/$ is adjacent to $v_{j-g+1}\/$ if $j <i\/$,  and $v_j\/$ is adjacent to $v_{j+g -1}\/$ otherwise. But this means that $F_n\/$ is ${\cal F}_{n,g}(^*)\/$, a contradiction.

 \vspace{5mm} On the other hand, if there exists  no $i \in S\/$ such that $v_i\/$ is adjacent to $v_{i-g+1}\/$, then $v_j\/$ is adjacent to $v_{j+g-1}\/$ for all $j \in S\/$. But again $F_n\/$ is  is ${\cal F}_{n,g}(^*)\/$, a contradiction.

\vspace{5mm} {\em Case (2):     $m \geq 3g-4\/$.}

\vspace{5mm} Suppose  there is an $i \in S\/$ such that $v_i\/$ is adjacent to a vertex $v_j\/$ in $F_m\/$ with  $|j-i| \not \equiv 1 \  ({\rm mod}\ (g-2))\/$. Then the diagonals
 $v_iv_{i+r(g-2)+1}\/$, $r = 1, 2,   \ldots, t\/$ together with the boundary edges of $K_n\/$ yield a $g\/$-angulation of  $K_n - F_n\/$.

\vspace{5mm} Hence assume that  for any $i\in S\/$, $v_i\/$ is adjacent to a vertex $v_j\/$ in $F_m\/$ such that  $|j-i| \equiv 1 \  ({\rm mod}\ (g-2))\/$.

\vspace{5mm}
  (i) Suppose  both $v_q$ and $v_{q+g-1}$ are non-pendent vertices in ${\cal F}_{m,g}(^*)\/$.

\vspace{5mm}  If $v_{q+1}\/$ is not adjacent to $v_{q-g+2}\/$ in $F_n\/$, then  obtained by deleting all the vertices $v_{q-g+3}, \ldots, v_{q}\/$ is of the form $K_{m} - F'_{m}\/$ where  $m = n-g+2\/$,  $|E(F'_{m})| \leq m -g\/$ (because $d_{F_n}(v_q) \geq 2\/$) and $F'_{m}\/$ contains no boundary edges of $K_{m}\/$. By Theorem \ref{atmostn-g}, $K_{m} - F'_{m}\/$ has a $g\/$-angulation which together $v_{q-g+1}v_{q-g+2} \cdots v_{q+1}\/$ form a $g\/$-angulation for $K_n - F_n\/$.

\vspace{5mm}   If $v_{q+1}\/$ is adjacent to $v_{q-g+2}\/$ in $F_n\/$, then $v_{q+1}\/$ is not adjacent to $v_{q+g-2}\/$ in $F_n\/$. In this case, we consider the convex graph obtained by deleting the vertices $v_{q+2}, v_{q+3}, \ldots, v_{q+g-1}\/$ and apply similar argument before to conclude that  $K_n - F_n\/$ admits  a $g\/$-angulation.

\vspace{5mm}
  (ii) Suppose  only one of  $v_q$ or $v_{q+g-1}$ is a non-pendent vertex in ${\cal F}_{m,g}(^*)\/$.

  \vspace{5mm} We can assume without loss of generality that  $v_q$ (since we can relabel the vertices of $K_n - F_n\/$). Then $v_{q+g-1}\/$ is a pendant vertex in $F_n\/$.  If for some $i \in S\/$, $v_i\/$ is not adjacent to $v_{i-g+1}\/$ in $F_n\/$, then by the method similar to case in (i), we see that $K_n - F_n\/$ admits a $g\/$-angulation. On the other hand, if $v_i\/$ is adjacent to $v_{i-g+1}\/$ in $F_n\/$ for all $i \in S\/$, then $F_n\/$ is
  ${\cal F}_{n,g}(^*)\/$, a contradiction.

\vspace{5mm}
  (iii) Suppose  both $v_q$ and $v_{q+g-1}$ are pendent vertices in ${\cal F}_{m,g}(^*)\/$.

 \vspace{5mm} If $i \in S\/$, we let $v_{s_i}$ be the neighbor of $v_i$ in $F_n\/$.

\vspace{5mm} Suppose  there exist $i, j \in S\/$ such that

\vspace{5mm}  \hspace{15mm}  $i < j\/$ and    $v_iv_{s_i}\/$ and $v_jv_{s_j}\/$ crosse in  $F_n\/$ with $s_j-s_i > g-2\/$  \hfill $(\star)\/$

\vspace{5mm} Let $H_1\/$ and $H_2\/$  denote the convex subgraphs of $K_n - F_n\/$ induced by the vertices
$v_{s_i+a},  \ldots, v_{s_j}, \ldots, v_i, \ldots, v_{j}\/$  and $v_j, \ldots, v_{s_i}, \ldots, v_{s_i+a} \/$ respectively
where $a\in\{1, \ldots, g-2\}\/$ with $|s_i+a-j|\equiv 1({\rm mod}\ (g-2))\/$.
Then $H_1\/$ and $H_2\/$ each admits a $g\/$-angulation (since $v_i$ and $v_j$  is adjacent to every other vertex of $H_1\/$ and $H_2\/$ respectively)  which together yields a $g\/$-angulation for $K_n - F_n\/$.

\vspace{5mm} Hence we assume that no $i, j\in S\/$ satisfy the condition $(\star)\/$.

\vspace{5mm} Suppose that  for some $i \in S\/$, $v_{s_i-j}v_{s_i-j+g-1}\notin E(F_n)$ for some $j\in \{1, \ldots, g-2\}\/$.

\vspace{5mm} Then the subgraph obtained by deleting $g-2\/$  vertices $v_{s_i-j+1}, \ldots, v_{s_i},  \ldots, \linebreak v_{s_i-j+g-2}\/$ (from $K_n - F_n\/$) is of the form $K_{m} - F'_{m}\/$ where $m=n-g+2$ with $|E(F'_{m})| \leq m -g\/$ (because $d_{F_n}(v_{s_i}) \geq 2\/$) and $F'_{m}\/$ contains no boundary edges of $K_{m}\/$. By Theorem \ref{atmostn-g}, $K_{m} - F'_{m}\/$ has a $g\/$-angulation which together with $v_{s_i-j} \cdots v_{s_i}  \cdots v_{s_i-j+g-1}\/$ form a $g\/$-angulation for $K_n - F_n\/$.

 \vspace{5mm}  Hence we assume that for any $i \in S\/$, $v_{s_i-j}v_{s_i-j+g-1}\/$ is an edge in  $F_n\/$ for any  $j\in \{1, \ldots, g-2\}\/$.

\vspace{5mm} Suppose there is a pendant vertex $v_r\/$ such that $r \not \in S\/$ and the edge  $v_rv_s\/$  (incident to $v_r\/$) crosses $v_iv_{s_i}\/$ ($i \in S\/$) with $|s- s_i| > g-2\/$, then again a $g\/$-angulation of $K_n - F_n\/$ can be constructed as in the previous case (where the condition $(\star)\/$ is satisfied).

\vspace{5mm} Hence,  for any pendant vertex $v_r\/$ where $r \not \in S\/$, $v_rv_s\/$ does not crosse $v_iv_{s_i}\/$ with  $|s- s_i| > g-2\/$ for any $i \in S\/$. But this implies that $F_n\/$ is  $ {\cal F}_{n,g}(^*)\/$, a contradiction.

\vspace{1mm} This completes the proof.  \qed

\section{ $n-g +1 + \mu\/$ edges}   \label{n-1g}

We present in this section a characterization of  ${\cal J}_{n,g}(^*_\beta)\/$ (where $\beta\in \{1, 2, \ldots, 2g-3\}\/$) of size at most $n-1\/$ that shares any possible $g$-angulation of $K_n\/$ by at least one edge.
We show that  $G=K_n-F_n\/$  
admits a $g$-angulation if and only if $F_n\neq {\cal J}_{n,g}(^*_\beta)\/$.

\vspace{2mm}
\begin{define}  \label{jmu1}
Suppose  $g \geq 3\/$ is a natural number and let  $n=g+t(g-2)\/$ where  $t\geq 1$ is any natural number.
 Let ${\cal J}_{n,g}(^* _1)\/$ denote a convex geometric graph with $n\/$ vertices and  $n-g+1+\mu\/$ edges
  such that ${\cal J}_{n,g}(^* _1)- \{e_1, e_2, \ldots, e_{\mu}\}={\cal F}_{n,g}(^*)\/$ for some  $\mu$ edges $e_1, e_2, \ldots, e_{\mu}$. Here   $1\leq\mu\leq g-2$.

\end{define}

\vspace{2mm}

\begin{propo} \label{Jstar1}
    $K_n - {\cal J}_{n,g}(^* _1)\/$ admits no $g$-angulation for any natural number $n=g+t(g-2)\/$ with  $t\geq 1$.
\end{propo}

\vspace{1mm}  \noindent {\bf Proof:} Suppose on the contrary that   $K_n - {\cal J}_{n,g}(^* _1)\/$  admits a $g$-angulation $G_n\/$. Since
  ${\cal J}_{n,g}(^* _1)$  contains  $\mu$ edges $\{e_1, e_2, \ldots, e_{\mu}\}$   such that ${\cal J}_{n,g}(^* _1)- \{e_1, e_2, \ldots, e_{\mu}\}={\cal F}_{n,g}(^*)\/$, it follows that   $K_n-{\cal F}_{n,g}(^*)\/$ admits $G_n\/$, which is a contradiction.

\vspace{1mm} This completes the proof.  \qed

\vspace{8mm} Suppose  $g \geq 3\/$ is a natural number and let  $n=g+t(g-2)\/$ where  $t\geq 0$ is any natural number.  Let $F\/$ denote a subgraph of $K_n\/$.

\begin{itemize}
  \item  A pair of vertices $\{v_r, v_s\}\/$ is called  an {\em  $\alpha\/$-pair \/} in $F\/$ if whenever pair of edges  $v_iv_r, v_jv_s\/$ (with $|r - i|\equiv 1({\rm mod}\ (g-2))\/$   and   $|s - j|\equiv 1({\rm mod}\ (g-2))\/$), then $v_iv_r\/$ crosses $v_jv_s\/$ with   $|i - j|> g-2\/$.

  \item A vertex $v_i\in V(F)$ is called  a {\em  $g$-angulable vertex \/} in $F\/$ if whenever edge  $v_iv_j\/$ (if there exists) in $F\/$   satisfies that  $|j-i| \not\equiv 1 \  ({\rm mod}\ (g-2))\/$.
\end{itemize}

\vspace{2mm}
\begin{define}  \label{jmu2}
Suppose  $g \geq 3\/$ is a natural number and let  $n=g+t(g-2)\/$ where $t\geq 2$ is any natural number.
 Let ${\cal J}_{n,g}(^* _2)\/$ denote a convex geometric graph with $n\/$ vertices and  $n-g+1+\mu\/$ edges (where $1 \leq \mu \leq g-2\/$) with no $g\/$-angulable vertex in $K_n-{\cal J}_{n,g}(^* _2)\/$
  such that

 (i) if $\{v_r, v_s\}\/$ is an   $\alpha\/$-pair of vertices in ${\cal J}_{n,g}(^* _2)\/$ (which occurs at most once), then $v_r, v_s\/$ are adjacent in  ${\cal J}_{n,g}(^* _2)\/$ with $|r - s|\equiv 1({\rm mod}\ (g-2))\/$, and

 (ii)  whenever $v_iv_{i+g-1}\not \in   E({\cal J}_{n,g}(^* _2))$ and $\{v_{i+1}, \ldots, v_{i+g-2}\}\/$  contains non-pendant vertex (which is a vertex of an $\alpha\/$-pair or a vertex in the neighbor of a vertex of an $\alpha\/$-pair in case ${\cal J}_{n,g}(^* _2)\/$ has an $\alpha\/$-pair),
 then ${\cal J}_{n,g}(^* _2) - \{v_{i+1}, \ldots,  v_{i+g-2}\}\/$ is either
  ${\cal F}_{m,g}(^*)\/$ or ${\cal J}_{m,g}(^* _1)\/$   (where $m = n-g+2\/$).
\end{define}

\vspace{5mm}
 \begin{note} Suppose $v_iv_{i+g-1}\not \in   E({\cal J}_{n,g}(^* _2))$ and
 either  all vertices in $\{v_{i+1}, \linebreak \ldots,  v_{i+g-2}\}\/$ are pendant (in case ${\cal J}_{n,g}(^* _2)\/$ has no $\alpha\/$-pair)
  or $\{v_{i+1}, \ldots,  v_{i+g-2}\}\/$ contains a non-pendant vertex but neither a vertex of an $\alpha\/$-pair nor a vertex in the neighbor of a vertex of an $\alpha\/$-pair (in case ${\cal J}_{n,g}(^* _2)\/$ has an $\alpha\/$-pair),
  then it is easy to see that ${\cal J}_{n,g}(^* _2) - \{v_{i+1}, \ldots, v_{i+g-2}\}\/$ is ${\cal J}_{m,g}(^* _2)$ where $m = n -g+2\/$.
  See  Figure \ref{alpha-pair of vertices  in Fm gives type-2}.
 \end{note}

\vspace{5mm} An example of  a geometric graph ${\cal J}_{n,g}(^* _2)\/$ with and without $\alpha\/$-pair is depicted in Figure \ref{example-type-2} (a) and (b) respectively.  In Figure \ref{example-type-2} (a), $\{v_6, v_{10}\}\/$ is the $\alpha\/$-pair in ${\cal J}_{14,5}(^* _2)\/$.

\begin{figure}[htp]
\begin{center}
\begin{minipage}{.48\textwidth}
\begin{center}
\resizebox{6.cm}{!}{
\begin{tikzpicture}
\coordinate (center) at (0,0);
  \def\radius{2.5cm}
   \foreach \x in {0,25.71428,...,360} {
             \filldraw[] (\x:2cm) circle(1pt);
               }


      \coordinate (v0) at (-1.26,1.58);\filldraw[black] (v0) circle(2.5pt);  \node[above left] at (v0) {${v_0}$};
       \coordinate (v1) at (-0.45,1.96);\filldraw[black] (v1) circle(2.5pt);       \node[above] at (v1) {${v_1}$};
       \coordinate (v2) at (0.45,1.96);\filldraw[black] (v2) circle(2.5pt);        \node[above] at (v2) {${v_2}$};
       \coordinate (v3) at (1.26,1.58);\filldraw[black] (v3) circle(2.5pt);\node[above right] at (v3) {${v_{3}}$};
       \coordinate (v4) at (1.8,0.87);\filldraw[black] (v4) circle(2.5pt);         \node[right] at (v4) {${v_4}$};
       \coordinate (v5) at (2,0);\filldraw[black] (v5) circle(2.5pt);              \node[right] at (v5) {${v_5}$};
       \coordinate (v6) at (1.8,-0.87);\filldraw[black] (v6) circle(2.5pt);\node[below right] at (v6) {${v_{6}}$};
       \coordinate (v7) at (1.26,-1.58);\filldraw[black] (v7) circle(2.5pt);       \node[below] at (v7) {${v_7}$};
       \coordinate (v8) at (0.45,-1.96);\filldraw[black] (v8) circle(2.5pt);       \node[below] at (v8) {${v_8}$};
       \coordinate (v9) at (-0.45,-1.96);\filldraw[black] (v9) circle(2.5pt);  \node[below left] at (v9) {${v_9}$};
        \coordinate (v10) at (-1.26,-1.58);\filldraw[black] (v10) circle(2.5pt);\node[below left] at (v10) {${v_{10}}$};
        \coordinate (v11) at (-1.8,-0.87);\filldraw[black] (v11) circle(2.5pt);       \node[left] at (v11) {${v_{11}}$};
        \coordinate (v12) at (-2,0);\filldraw[black] (v12) circle(2.5pt);           \node[ left] at (v12) {${v_{12}}$};
       \coordinate (v13) at (-1.8,0.87);\filldraw[black] (v13) circle(2.5pt);        \node[ left] at (v13) {${v_{13}}$};


       \draw [line width=1,black](v13) -- (v3);
    \draw [line width=1,black](v9) -- (v2);
    \draw [line width=1,black](v8) -- (v1);
    \draw [line width=1,black](v7) -- (v0);

       \draw [line width=1.2,blue](v3) --(v10) --(v6) -- (v13);
       \draw [line width=1,black](v11) -- (v1);
       \draw [line width=1,black](v12) -- (v2);

      \draw [line width=1,black](v5) -- (v1);
       \draw [line width=1,black](v4) -- (v0);

       \draw [line width=1,black](v3) -- (v7);
       \draw [line width=1,black](v5) -- (v9);

\filldraw[white] (v0) circle(1.5pt);
\filldraw[white] (v1) circle(1.5pt);
\filldraw[white] (v2) circle(1.5pt);
\filldraw[white] (v3) circle(1.5pt);
\filldraw[white] (v4) circle(1.5pt);
\filldraw[white] (v5) circle(1.5pt);
\filldraw[white] (v6) circle(1.5pt);
\filldraw[white] (v7) circle(1.5pt);
\filldraw[white] (v8) circle(1.5pt);
\filldraw[white] (v9) circle(1.5pt);
\filldraw[white] (v10) circle(1.5pt);
\filldraw[white] (v11) circle(1.5pt);
\filldraw[white] (v12) circle(1.5pt);
\filldraw[white] (v13) circle(1.5pt);

\end{tikzpicture}
}

\centering{(a)}

\end{center}
\end{minipage}
\begin{minipage}{.48\textwidth}
\begin{center}
\resizebox{6.cm}{!}{
\begin{tikzpicture}
\coordinate (center) at (0,0);
  \def\radius{2.5cm}
   \foreach \x in {0,25.71428,...,360} {
             \filldraw[] (\x:2cm) circle(1pt);
               }


       \coordinate (v0) at (-1.26,1.58);\filldraw[black] (v0) circle(2.5pt);\node[above left] at (v0) {${v_0}$};
       \coordinate (v1) at (-0.45,1.96);\filldraw[black] (v1) circle(2.5pt);\node[above] at (v1) {${v_1}$};
       \coordinate (v2) at (0.45,1.96);\filldraw[black] (v2) circle(2.5pt);\node[above] at (v2) {${v_2}$};
       \coordinate (v3) at (1.26,1.58);\filldraw[black] (v3) circle(2.5pt);\node[above right] at (v3) {${v_{3}}$};
       \coordinate (v4) at (1.8,0.87);\filldraw[black] (v4) circle(2.5pt);\node[right] at (v4) {${v_4}$};
       \coordinate (v5) at (2,0);\filldraw[black] (v5) circle(2.5pt);\node[right] at (v5) {${v_5}$};
       \coordinate (v6) at (1.8,-0.87);\filldraw[black] (v6) circle(2.5pt);\node[below right] at (v6) {${v_{6}}$};
       \coordinate (v7) at (1.26,-1.58);\filldraw[black] (v7) circle(2.5pt);\node[below] at (v7) {${v_7}$};
       \coordinate (v8) at (0.45,-1.96);\filldraw[black] (v8) circle(2.5pt);\node[below] at (v8) {${v_8}$};
       \coordinate (v9) at (-0.45,-1.96);\filldraw[black] (v9) circle(2.5pt);\node[below left] at (v9) {${v_9}$};
        \coordinate (v10) at (-1.26,-1.58);\filldraw[black] (v10) circle(2.5pt);\node[below left] at (v10) {${v_{10}}$};
        \coordinate (v11) at (-1.8,-0.87);\filldraw[black] (v11) circle(2.5pt);\node[below left] at (v11) {${v_{11}}$};
        \coordinate (v12) at (-2,0);\filldraw[black] (v12) circle(2.5pt);\node[ left] at (v12) {${v_{12}}$};
       \coordinate (v13) at (-1.8,0.87);\filldraw[black] (v13) circle(2.5pt);\node[ left] at (v13) {${v_{13}}$};

     \draw [line width=1.2,orange](v1)-- (v11);
     \draw [line width=1.2,orange](v1)-- (v5);
     \draw [line width=1.2,purple](v12)-- (v2)-- (v6);
     \draw [line width=1.2,blue](v13)-- (v3)-- (v7);

     \draw [line width=1,black] (v11)-- (v4)-- (v8);
     \draw [line width=1,black](v12)-- (v5)-- (v9);
     \draw [line width=1,black] (v6)-- (v13);
     \draw [line width=1,black](v10)-- (v0)-- (v7);

    \filldraw[black] (v0) circle(2.5pt);
    \filldraw[black] (v1) circle(2.5pt);
    \filldraw[black] (v2) circle(2.5pt);
   \filldraw[black] (v3) circle(2.5pt);
    \filldraw[black] (v4) circle(2.5pt);
   \filldraw[black] (v5) circle(2.5pt);
    \filldraw[black] (v6) circle(2.5pt);
    \filldraw[black] (v7) circle(2.5pt);
   \filldraw[black] (v8) circle(2.5pt);
    \filldraw[black] (v9) circle(2.5pt);
   \filldraw[black] (v10) circle(2.5pt);
    \filldraw[black] (v11) circle(2.5pt);
    \filldraw[black] (v12) circle(2.5pt);
    \filldraw[black] (v13) circle(2.5pt);

    \filldraw[orange] (v1) circle(1.5pt);
    \filldraw[purple] (v2) circle(1.5pt);
    \filldraw[blue] (v3) circle(1.5pt);

    \filldraw[white] (v4) circle(1.5pt);
    \filldraw[white] (v5) circle(1.5pt);
    \filldraw[white] (v6) circle(1.5pt);
    \filldraw[white] (v7) circle(1.5pt);
    \filldraw[white] (v8) circle(1.5pt);
    \filldraw[white] (v9) circle(1.5pt);

    \filldraw[white] (v10) circle(1.5pt);
    \filldraw[white] (v11) circle(1.5pt);
    \filldraw[white] (v12) circle(1.5pt);
    \filldraw[white] (v13) circle(1.5pt);
    \filldraw[white] (v0) circle(1.5pt);

\end{tikzpicture}
}

\centering (b)

\end{center}
\end{minipage}
\caption{${\cal J}_{14,5}(^*_2)\/$}   \label{example-type-2}

\end{center}
\end{figure}

\vspace{2mm}

\begin{propo} \label{Jstar2}
     $K_n - {\cal J}_{n,g}(^*_2)\/$ admits no $g$-angulation for any natural number $n=g+t(g-2)\/$ with $t\geq 2$.

\end{propo}

\vspace{1mm}  \noindent {\bf Proof:}We prove this by induction on $n\/$.

\vspace{5mm} Consider the case $t=2\/$. Here $n = 3g-4\/$.

\vspace{5mm} Lemma \ref{earg} asserts the existence of an edge $v_qv_{q+g-1}\/$ in $K_n - {\cal J}_{n,g}(^* _2)\/$. Choose $q\/$ such that $\{v_{q+1}, \ldots, v_{q+g-2}\}\/$ contains non-pendant vertex.  
By definition  $ {\cal J}_{n,g}(^* _2) - \{v_{q+1}, \ldots, v_{q+g-2}\} = F_m\/$ is either ${\cal F}_{m,g}(^*)\/$ or ${\cal J}_{m,g}(^* _1)\/$, $K_m - F_m\/$ admits no $g\/$-angulation by Theorem \ref{fstarsg} and Proposition \ref{Jstar1}. Here $m = n-g+2\/$.

\vspace{5mm}  Assume on the contrary that $K_n - {\cal J}_{n,g}(^* _2)\/$ admits a $g$-angulation $G_n\/$. Then $G_n \/$ does not contain the edge $v_qv_{q+g-1}\/$  (otherwise this implies that $K_m - F_m\/$ admits $g\/$-angulation, a contradiction).

 \vspace{5mm} Since $t=2\/$, $G_n\/$ has only two diagonal edges.  Clearly  at least one of these diagonal edges, say $e_1\/$ is incident to a vertex in $\{v_{q+1}, \ldots, v_{q+g-2}\}\/$. This is clearly not possible if all vertices in $\{v_{q+1}, \ldots, v_{q+g-2}\}\/$ are non-pendant. Hence assume that $e_1\/$ is a incident to a pendant vertex in $\{v_{q+1}, \ldots, v_{q+g-2}\}\/$. That is $e_1 = v_jv_r\/$ where $j \in \{q+1, \ldots, q+g-2\}\/$. This means that $v_jv_s\/$ is an edge of ${\cal J}_{n,g}(^* _2)\/$ (where $|r-s| =g-2\/$)  and $v_s\/$ is a non-pendant vertex in ${\cal J}_{n,g}(^* _2)\/$ (otherwise $v_s\/$ is an isolated vertex in  $F_m ={\cal J}_{n,g}(^* _2) - \{v_{q+1}, \ldots, v_{q+g-2}\}\/$, a contradiction).

 \vspace{5mm} Let $v_pv_{p+g-1}\/$ be the second diagonal edge of $G_n\/$ (where $ p > j\/$). It is easy to see that $v_s\/$ is one of the vertex in $\{v_{p+1}, \ldots, v_{p+g-2}\}\/$ (see for example Figure \ref{inetial step type 2}).    But this means that $v_j\/$ is  an isolated vertex in  ${\cal J}_{n,g}(^* _2) - \{v_{p+1}, \ldots, v_{p+g-2}\}\/$ (a contradiction).

\vspace{5mm} Now suppose $t \geq 3\/$.

\vspace{5mm} Assume on the contrary that $K_n - {\cal  J}_{n,g}(^*_2)$ admits a $g$-angulation $G_n\/$.
Note that any  $g$-angulation on a set of points in convex position has a diagonal $v_qv_{q+g-1}\/$. 

\vspace{5mm} Let $ F_m = {\cal J}_{n,g}(^* _2) - \{v_{q+1}, \ldots, v_{q+g-2}\}\/$ where $m = n-g+2\/$.
Then $G_n -\{v_{q+1}, \ldots, v_{q+g-2}\}\/$ is  a $g$-angulation  for $K_m-F_m\/$.

\vspace{5mm} {\em Case (1):    }  ${\cal J}_{n,g}(^* _2)\/$ contains no $\alpha\/$-pair.

\vspace{5mm} (i)  Suppose  $\{v_{q+1}, \ldots, v_{q+g-2}\}\/$ contains non-pendent vertex.
Then by the definition of ${\cal  J}_{n,g}(^*_2)\/$,   $F_m\/$ is either ${\cal F}_{m,g}(^*)\/$ or $ {\cal J}_{m,g}(^*_1)\/$.
By Theorem \ref{fstarsg} or Proposition \ref{Jstar1} respectively, $K_m-F_m\/$ admits no  $g$-angulation, a contradiction.

\vspace{5mm} (ii)  Suppose all vertices in  $\{v_{q+1}, \ldots, v_{q+g-2}\}\/$  are pendant vertices.
  Then $F_m={\cal J}_{m,g}(^*_2)\/$. By induction, $K_m-F_m\/$ admits no  $g$-angulation, again a contradiction that  $G_n$ is a  $g$-angulation for $K_n - {\cal  J}_{n,g}(^*_{2})\/$.

\vspace{5mm} {\em Case (2):    }  ${\cal J}_{n,g}(^* _2)\/$ contains an $\alpha\/$-pair.

\vspace{5mm}  If $\{v_{q+1}, \ldots, v_{q+g-2}\}\/$ contains a vertex of an $\alpha\/$-pair or a vertex in the neighbor of a vertex of an $\alpha\/$-pair, then (by definition of ${\cal  J}_{n,g}(^*_2)\/$,)   $F_m\/$ is either ${\cal F}_{m,g}(^*)\/$ or $ {\cal J}_{m,g}(^*_1)\/$.
If $\{v_{q+1}, \ldots, v_{q+g-2}\}\/$ contains neither a vertex of an $\alpha\/$-pair nor a vertex in the neighbor of a vertex of an $\alpha\/$-pair, then
   $F_m\/$ is ${\cal J}_{m,g}(^*_2)\/$.

\vspace{5mm} In any case, by Theorem \ref{fstarsg},  Proposition \ref{Jstar1} or by induction,  $K_m-F_m\/$ admits no  $g$-angulation,  a contradiction that  $G_n$ is a  $g$-angulation for $K_n - {\cal  J}_{n,g}(^*_{2})\/$.

 \vspace{5mm}    This completes the proof.  \qed

\begin{figure}[htp]
\begin{center}
\begin{minipage}{.49\textwidth}
\begin{center}
\resizebox{5.5cm}{!}{
\begin{tikzpicture}[rotate=2,.style={draw}]
\coordinate (center) at (0,0);
   \def\radius{1.cm}
  \foreach \x in {0, 45,...,360} {
            \filldraw[] (\x:1.cm) circle(1pt);
               }

       \coordinate (v1) at (.68,0.68);\filldraw[blue] (v1) circle(1pt);\node[right] at (v1) {\tiny${v_{1}}$};
       \coordinate (v2) at (1,0);\filldraw[black] (v2) circle(1pt);\node[right] at (v2) {\tiny${v_2}$};
       \coordinate (v3) at (.68,-0.68);\filldraw[black] (v3) circle(1pt);\node[right] at (v3) {\tiny${v_3}$};
       \coordinate (v4) at (0,-1);\filldraw[black] (v4) circle(1pt);\node[below] at (v4) {\tiny${{v_4}}$};
       \coordinate (v5) at (-.68,-0.68);\filldraw[blue] (v5) circle(1pt);\node[left] at (v5) {\tiny${v_5}$};
       \coordinate (v6) at (-1,0);\filldraw[blue] (v6) circle(1pt);\node[left] at (v6) {\tiny${v_6}$};
       \coordinate (v7) at (-.68,0.68);\filldraw[blue] (v7) circle(1pt);\node[left] at (v7) {\tiny${v_{7}}$};
       \coordinate (v0) at (0,1.);\filldraw[blue] (v0) circle(1pt);\node[above] at (v0) {\tiny${v_0}$};

   \draw [line width=0.7,black](v6) -- (v1) -- (v4);
   \draw [line width=0.7,black](v5) -- (v0) -- (v3);

    \draw [line width=0.7,black](v7) -- (v2);
    \draw [line width=0.7,black](v1) -- (v3);

\filldraw[black] (v0) circle(2.pt);
\filldraw[black] (v1) circle(2pt);
\filldraw[black] (v2) circle(2pt);
\filldraw[black] (v3) circle(2pt);
\filldraw[black] (v4) circle(2pt);
\filldraw[black] (v5) circle(2pt);
\filldraw[black] (v6) circle(2pt);
\filldraw[black] (v7) circle(2pt);

\filldraw[white] (v0) circle(1.3pt);
\filldraw[white] (v1) circle(1.3pt);
\filldraw[white] (v2) circle(1.3pt);
\filldraw[white] (v3) circle(1.3pt);
\filldraw[white] (v4) circle(1.3pt);
\filldraw[white] (v5) circle(1.3pt);
\filldraw[white] (v6) circle(1.3pt);
\filldraw[white] (v7) circle(1.3pt);

\end{tikzpicture}
}
\end{center}

\end{minipage}
\end{center}
\caption{${\cal J}_{8,4}(^*_2)$.}    \label{inetial step type 2}
\end{figure}


\vspace{2mm}
\begin{define}   \label{jmu-gamma}
Suppose  $g \geq 4\/$ and $\gamma \geq 3\/$ are natural numbers and let  $n=g+t(g-2)\/$ where   $t \geq 3\/$ is any natural number.
 Let ${\cal J}_{n,g}(^* _{\gamma})\/$  denote a convex geometric graph with $n\/$ vertices and  $n-g+1+\mu\/$ edges  (where $1 \leq \mu \leq g-2\/$) with no  $g$-angulable vertex and having only $\gamma-1\/$ $\alpha\/$-pairs of vertices    such that

   \vspace{5mm}   (i) whenever $\{v_r, v_s\}\/$ is an   $\alpha\/$-pair of vertices, then $v_r, v_s\/$ are adjacent in  ${\cal J}_{n,g}(^*_{\gamma})\/$ with $|r - s|\equiv 1({\rm mod}\ (g-2))\/$, and

  \vspace{5mm}   (ii)  whenever $v_iv_{i+g-1}\not \in   E({\cal J}_{n,g}(^* _{\gamma}))$ and $\{v_{i+1}, \ldots, v_{i+g-2}\}\/$  contains a vertex of an $\alpha\/$-pair or a vertex in the neighbor of a vertex of an $\alpha\/$-pair, then ${\cal J}_{n,g}(^* _{\gamma}) - \{v_{i+1}, \ldots,  v_{i+g-2}\}\/$ is either   ${\cal F}_{m,g}(^*)\/$ or ${\cal J}_{m,g}(^* _\beta)\/$   where $m = n-g+2\/$ and  $\beta \leq \gamma-1\/$.
\end{define}

\begin{figure}[htp]
\begin{center}
\begin{minipage}{.49\textwidth}
\begin{center}
\resizebox{6.cm}{!}{
\begin{tikzpicture}
\coordinate (center) at (0,0);
   \def\radius{3cm}
   \foreach \x in {0, 20,...,360} {
             \filldraw[] (\x:3cm) circle(1pt);
               }


      \coordinate (v0) at (-1.51, 2.61);\filldraw[blue] (v0) circle(1pt);\node[above left] at (v0) {\Large${v_0}$};

       \coordinate (v1) at (-.52, 2.94);\filldraw[blue] (v1) circle(1pt);\node[above left] at (v1) {\Large${v_1}$};
      \coordinate (v2) at (.52, 2.94);\filldraw[blue] (v2) circle(1pt);\node[above right] at (v2) {\Large${v_2}$};

       \coordinate (v3) at (1.51, 2.61);\filldraw[blue] (v3) circle(1pt);\node[above right] at (v3) {\Large${v_{3}}$};
       \coordinate (v4) at (2.3,1.93);\filldraw[blue] (v4) circle(1pt);\node[above right] at (v4) {\Large${v_4}$};

       \coordinate (v5) at (2.82, 1.02);\filldraw[blue] (v5) circle(1pt);\node[above right] at (v5) {\Large${v_5}$};
       \coordinate (v6) at (3,0);\filldraw[blue] (v6) circle(1pt);\node[right] at (v6) {\Large${v_{6}}$};

       \coordinate (v7) at (2.82, -1.02);\filldraw[blue] (v7) circle(1pt);\node[below right] at (v7) {\Large${v_7}$};

       \coordinate (v8) at (2.3,-1.93);\filldraw[blue] (v8) circle(1pt);\node[below right] at (v8) {\Large${v_8}$};

       \coordinate (v9) at (1.51, -2.61);\filldraw[blue] (v9) circle(1pt);\node[below right] at (v9) {\Large${v_{9}}$};
       \coordinate (v10) at (.52, -2.94);\filldraw[blue] (v10) circle(1pt);\node[below right] at (v10) {\Large${v_{10}}$};

       \coordinate (v11) at (-.52, -2.94);\filldraw[blue] (v11) circle(1pt);\node[below left] at (v11) {\Large${v_{11}}$};

       \coordinate (v12) at (-1.51, -2.61);\filldraw[blue] (v12) circle(1pt);\node[below left] at (v12) {\Large${v_{12}}$};

       \coordinate (v13) at (-2.3,-1.93);\filldraw[blue] (v13) circle(1pt);\node[below left] at (v13) {\Large${v_{13}}$};

       \coordinate (v14) at (-2.82, -1.02);\filldraw[blue] (v14) circle(1pt);\node[below left] at (v14) {\Large${v_{14}}$};

       \coordinate (v15) at (-3,0);\filldraw[blue] (v15) circle(1pt);\node[left] at (v15) {\Large${v_{15}}$};

       \coordinate (v16) at (-2.82, 1.02);\filldraw[blue] (v16) circle(1pt);\node[above left] at (v16) {\Large${v_{16}}$};

       \coordinate (v17) at (-2.3,1.93);\filldraw[blue] (v17) circle(1pt);\node[above left] at (v17) {\Large${v_{17}}$};

     \draw [line width=1.3,black](v17)-- (v2) -- (v5)--(v8);
     \draw [line width=1.3,black](v0)-- (v3) -- (v6);
     \draw [line width=1.3,black](v1)-- (v4) -- (v7);

     \draw [line width=1.3,black](v11)-- (v6);
     \draw [line width=1.3,black](v10)-- (v5);
     \draw [line width=1.3,black](v15)-- (v2);
     \draw [line width=1.3,black](v14)-- (v1);

     \draw [line width=1.7,blue](v0)-- (v13) -- (v16)--(v3);
     \draw [line width=1.7,orange](v4)-- (v9) -- (v12)--(v7);

    \filldraw[black] (v0) circle(3.5pt);
    \filldraw[black] (v1) circle(3.5pt);
    \filldraw[black] (v2) circle(3.5pt);
   \filldraw[black] (v3) circle(3.5pt);
    \filldraw[black] (v4) circle(3.5pt);
   \filldraw[black] (v5) circle(3.5pt);
    \filldraw[black] (v6) circle(3.5pt);
    \filldraw[black] (v7) circle(3.5pt);
   \filldraw[black] (v8) circle(3.5pt);
    \filldraw[black] (v9) circle(3.5pt);
   \filldraw[black] (v10) circle(3.5pt);
    \filldraw[black] (v11) circle(3.5pt);
    \filldraw[black] (v12) circle(3.5pt);
    \filldraw[black] (v13) circle(3.5pt);
    \filldraw[black] (v14) circle(3.5pt);
    \filldraw[black] (v15) circle(3.5pt);
    \filldraw[black] (v16) circle(3.5pt);
    \filldraw[black] (v17) circle(3.5pt);

    \filldraw[white] (v0) circle(2pt);
    \filldraw[white] (v1) circle(2pt);
    \filldraw[white] (v2) circle(2pt);
    \filldraw[white] (v3) circle(2pt);

    \filldraw[white] (v4) circle(2pt);
    \filldraw[white] (v5) circle(2pt);
    \filldraw[white] (v6) circle(2pt);
    \filldraw[white] (v7) circle(2pt);
    \filldraw[white] (v8) circle(2pt);
    \filldraw[white] (v9) circle(2pt);

    \filldraw[white] (v10) circle(2pt);
    \filldraw[white] (v11) circle(2pt);
    \filldraw[white] (v12) circle(2pt);
    \filldraw[white] (v13) circle(2pt);

    \filldraw[white] (v14) circle(2pt);
    \filldraw[white] (v15) circle(2pt);
    \filldraw[white] (v16) circle(2pt);
    \filldraw[white] (v17) circle(2pt);

\end{tikzpicture}
}

\centering{(a)  ${\cal J}_{18,4}(^*_3)\/$}
\end{center}
\end{minipage}
\begin{minipage}{.49\textwidth}
\begin{center}
\resizebox{6.cm}{!}{
\begin{tikzpicture}
\coordinate (center) at (0,0);
   \def\radius{3cm}
   \foreach \x in {0, 20,...,360} {
             \filldraw[] (\x:3cm) circle(1pt);
               }


      \coordinate (v0) at (-1.51, 2.61);\filldraw[blue] (v0) circle(1pt);\node[above left] at (v0) {\Large${v_0}$};

       \coordinate (v1) at (-.52, 2.94);\filldraw[blue] (v1) circle(1pt);\node[above left] at (v1) {\Large${v_1}$};
      \coordinate (v2) at (.52, 2.94);\filldraw[blue] (v2) circle(1pt);\node[above right] at (v2) {\Large${v_2}$};

       \coordinate (v3) at (1.51, 2.61);\filldraw[blue] (v3) circle(1pt);\node[above right] at (v3) {\Large${v_{3}}$};
       \coordinate (v4) at (2.3,1.93);\filldraw[blue] (v4) circle(1pt);\node[above right] at (v4) {\Large${v_4}$};

       \coordinate (v5) at (2.82, 1.02);\filldraw[blue] (v5) circle(1pt);\node[above right] at (v5) {\Large${v_5}$};
       \coordinate (v6) at (3,0);\filldraw[blue] (v6) circle(1pt);\node[right] at (v6) {\Large${v_{6}}$};

       \coordinate (v7) at (2.82, -1.02);\filldraw[blue] (v7) circle(1pt);\node[below right] at (v7) {\white${v_7}$};

       \coordinate (v8) at (2.3,-1.93);\filldraw[blue] (v8) circle(1pt);\node[below right] at (v8) {\white${v_8}$};

       \coordinate (v9) at (1.51, -2.61);\filldraw[blue] (v9) circle(1pt);\node[below right] at (v9) {\Large${v_{9}}$};
       \coordinate (v10) at (.52, -2.94);\filldraw[blue] (v10) circle(1pt);\node[below right] at (v10) {\Large${v_{10}}$};

       \coordinate (v11) at (-.52, -2.94);\filldraw[blue] (v11) circle(1pt);\node[below left] at (v11) {\Large${v_{11}}$};

       \coordinate (v12) at (-1.51, -2.61);\filldraw[blue] (v12) circle(1pt);\node[below left] at (v12) {\Large${v_{12}}$};

       \coordinate (v13) at (-2.3,-1.93);\filldraw[blue] (v13) circle(1pt);\node[below left] at (v13) {\Large${v_{13}}$};

       \coordinate (v14) at (-2.82, -1.02);\filldraw[blue] (v14) circle(1pt);\node[below left] at (v14) {\Large${v_{14}}$};

       \coordinate (v15) at (-3,0);\filldraw[blue] (v15) circle(1pt);\node[left] at (v15) {\Large${v_{15}}$};

       \coordinate (v16) at (-2.82, 1.02);\filldraw[blue] (v16) circle(1pt);\node[above left] at (v16) {\Large${v_{16}}$};

       \coordinate (v17) at (-2.3,1.93);\filldraw[blue] (v17) circle(1pt);\node[above left] at (v17) {\Large${v_{17}}$};

     \draw [line width=1.3,black](v17)-- (v2) -- (v5);
     \draw [line width=1.3,black](v0)-- (v3) -- (v6);
     \draw [line width=1.3,black](v1)-- (v4);

     \draw [line width=1.3,black](v11)-- (v6);
     \draw [line width=1.3,black](v10)-- (v5);
     \draw [line width=1.3,black](v15)-- (v2);
     \draw [line width=1.3,black](v14)-- (v1);

     \draw [line width=1.7,blue](v0)-- (v13) -- (v16)--(v3);
     \draw [line width=1.3,black](v4)-- (v9) -- (v12);

    \filldraw[black] (v0) circle(3.5pt);
    \filldraw[black] (v1) circle(3.5pt);
    \filldraw[black] (v2) circle(3.5pt);
   \filldraw[black] (v3) circle(3.5pt);
    \filldraw[black] (v4) circle(3.5pt);
   \filldraw[black] (v5) circle(3.5pt);
    \filldraw[black] (v6) circle(3.5pt);
    \filldraw[black] (v7) circle(3.5pt);
   \filldraw[black] (v8) circle(3.5pt);
    \filldraw[black] (v9) circle(3.5pt);
   \filldraw[black] (v10) circle(3.5pt);
    \filldraw[black] (v11) circle(3.5pt);
    \filldraw[black] (v12) circle(3.5pt);
    \filldraw[black] (v13) circle(3.5pt);
    \filldraw[black] (v14) circle(3.5pt);
    \filldraw[black] (v15) circle(3.5pt);
    \filldraw[black] (v16) circle(3.5pt);
    \filldraw[black] (v17) circle(3.5pt);

    \filldraw[white] (v0) circle(2pt);
    \filldraw[white] (v1) circle(2pt);
    \filldraw[white] (v2) circle(2pt);
    \filldraw[white] (v3) circle(2pt);

    \filldraw[white] (v4) circle(2pt);
    \filldraw[white] (v5) circle(2pt);
    \filldraw[white] (v6) circle(2pt);
    \filldraw[white] (v7) circle(4pt);
    \filldraw[white] (v8) circle(4pt);
    \filldraw[white] (v9) circle(2pt);

    \filldraw[white] (v10) circle(2pt);
    \filldraw[white] (v11) circle(2pt);
    \filldraw[white] (v12) circle(2pt);
    \filldraw[white] (v13) circle(2pt);

    \filldraw[white] (v14) circle(2pt);
    \filldraw[white] (v15) circle(2pt);
    \filldraw[white] (v16) circle(2pt);
    \filldraw[white] (v17) circle(2pt);

\end{tikzpicture}
}

\centering{ (b) ${\cal J}_{16,4}(^*_2)\/$}
\end{center}
\end{minipage}
\end{center}
\caption{${\cal J}_{n,g}(^*_3)\/$.}     \label{ex type 3}
\end{figure}

\vspace{5mm} An example of a convex geometric graph ${\cal J}_{n,g}(^*_3)\/$  is depicted in Figure \ref{ex type 3}.
 In Figure \ref{ex type 3} (a), $\{v_{9}, v_{12}\}$ and $\{v_{13}, v_{16}\}$ are two $\alpha\/$-pairs in  ${\cal J}_{18,4}(^*_3)\/$,
 in (b) $\{v_{13}, v_{16}\}$ is the $\alpha\/$-pair in the convex geometric graph ${\cal J}_{16,4}(^*_2)\/$ which obtained from  ${\cal J}_{18,4}(^*_3)\/$ by deleting two vertices $v_7\/$ and  $v_8\/$ (with $v_7\/$ is a neighbor of a vertex ($v_{12}\/$) of the $\alpha\/$-pair $\{v_9,  v_{12}\}\/$).
 For more instance see Figure \ref{example type 3}.

\vspace{2mm}

\begin{propo} \label{Jstar-gamma}
 $K_n - {\cal J}_{n,g}(^*_{\gamma})\/$ admits no $g$-angulation for any natural number $n=g+t(g-2)\/$ with $t\geq 3$.
\end{propo}

\vspace{1mm}  \noindent {\bf Proof:} We prove this by induction on $n\/$.

\vspace{5mm} Consider the case $t=3\/$. Here $n = 4g-6\/$.

\vspace{5mm} Assume on the contrary that $K_n - {\cal  J}_{n,g}(^*_\gamma)$ admits a $g$-angulation $G_n\/$.
Note that $G_n\/$  has a diagonal $v_qv_{q+g-1}\/$ for some $q\/$.

\vspace{5mm} Let $F_m={\cal J}_{n,g}(^*_\gamma)-\{v_{q+1}, \ldots, v_{q+g-2}\}\/$ where $m = n-g+2\/$. Then  $G_n-\{v_{q+1}, \ldots, v_{q+g-2}\}\/$ is a $g$-angulation of $K_m-F_m\/$.

\vspace{5mm} Let $\{v_r, v_s\}\/$ be an $\alpha$-pair in ${\cal J}_{n,g}(^*_\gamma)\/$ and let
$v_iv_r\/$ and $v_jv_s\/$ the edges incident to $v_r\/$ and $v_s\/$ respectively.     Here  $|i-r|\equiv 1({\rm mod}\ (g-2))\/$  and
 $|j-s|\equiv 1({\rm mod}\ (g-2))\/$.  Since $|i-j| > g-2\/$  and $|r-s| = g-1\/$.

\vspace{5mm}  Since $m = n-g+2\/$ then $t=2\/$. Hence, there is no any $\alpha$-pair in  $F_m\/$.

\vspace{5mm}  Thus, $\{v_{q+1}, \ldots, v_{q+g-2}\}\/$ contains a vertex of $v_i, v_j, v_r, v_s\/$.

\vspace{5mm} By definition of ${\cal J}_{n,g}(^*_\gamma)\/$, $F_m\/$  is either ${\cal F}_{m,g}(^*)\/$ or ${\cal J}_{m,g}(^* _\beta)\/$ with  $\beta=1, 2\/$.  \
 By Theorem \ref{fstarsg} or Proposition \ref{Jstar1} or else Proposition \ref{Jstar2}, $K_m-F_m\/$ admits no  $g$-angulation, a contradiction.

\vspace{5mm} Now suppose $t \geq 4\/$.

\vspace{5mm} Assume on the contrary that $K_n - {\cal  J}_{n,g}(^*_\gamma)$ admits a $g$-angulation $G_n\/$.
Then $G_n\/$ has a diagonal $v_qv_{q+g-1}\/$.

\vspace{5mm} Let $ F_m = {\cal J}_{n,g}(^* _\gamma) - \{v_{q+1}, \ldots, v_{q+g-2}\}\/$ where $m = n-g+2\/$.
Then $G_n -\{v_{q+1}, \ldots, v_{q+g-2}\}\/$ is  a $g$-angulation  for $K_m-F_m\/$.

\vspace{5mm}  If $\{v_{q+1}, \ldots, v_{q+g-2}\}\/$ contains a vertex of an $\alpha\/$-pair or a vertex in the neighbor of a vertex of an $\alpha\/$-pair, then (by definition of ${\cal  J}_{n,g}(^*_\gamma)\/$,)   $F_m\/$ is either ${\cal F}_{m,g}(^*)\/$ or $ {\cal J}_{m,g}(^*_\beta)\/$     where $m = n-g+2\/$ and  $\beta \leq \gamma-1\/$.
If $\{v_{q+1}, \ldots, v_{q+g-2}\}\/$ contains neither a vertex of an $\alpha\/$-pair nor a vertex in the neighbor of a vertex of an $\alpha\/$-pair, then
   $F_m\/$ is  ${\cal J}_{m,g}(^*_\gamma)\/$.

\vspace{5mm} In any case, by Theorem \ref{fstarsg},  Proposition \ref{Jstar1}, Proposition \ref{Jstar2} or by induction,  $K_m-F_m\/$ admits no  $g$-angulation,  a contradiction that  $G_n$ is a  $g$-angulation for $K_n - {\cal  J}_{n,g}(^*_{\gamma})\/$  (see for example Figure \ref{alpha-pair of vertices  in Fm gives type-2}).

\vspace{1mm}  This completes the proof.   \qed

\begin{figure}[htp]
\begin{center}
\begin{minipage}{.45\textwidth}
\begin{center}
\resizebox{6cm}{!}{
\begin{tikzpicture}
\coordinate (center) at (0,0);
   \def\radius{3cm}
   \foreach \x in {0, 20,...,360} {
             \filldraw[] (\x:3cm) circle(1pt);
               }


      \coordinate (v0) at (-1.51, 2.61);\filldraw[blue] (v0) circle(1pt); \node[above left] at (v0) {${v_0}$};

       \coordinate (v1) at (-.52, 2.94);\filldraw[blue] (v1) circle(1pt);           \node[above left] at (v1) {${v_1}$};
      \coordinate (v2) at (.52, 2.94);\filldraw[blue] (v2) circle(1pt);  \node[above right] at (v2) {${v_2}$};

       \coordinate (v3) at (1.51, 2.61);\filldraw[blue] (v3) circle(1pt);\node[above right] at (v3) {${v_{3}}$};
       \coordinate (v4) at (2.3,1.93);\filldraw[blue] (v4) circle(1pt);\node[above right] at (v4) {${v_4}$};

       \coordinate (v5) at (2.82, 1.02);\filldraw[blue] (v5) circle(1pt);\node[above right] at (v5) {${v_5}$};
       \coordinate (v6) at (3,0);\filldraw[blue] (v6) circle(1pt);                            \node[right] at (v6) {${v_{6}}$};

       \coordinate (v7) at (2.82, -1.02);\filldraw[blue] (v7) circle(1pt);\node[below right] at (v7) {${v_7}$};

       \coordinate (v8) at (2.3,-1.93);\filldraw[blue] (v8) circle(1pt);\node[below right] at (v8) {${v_8}$};

       \coordinate (v9) at (1.51, -2.61);\filldraw[blue] (v9) circle(1pt);\node[below right] at (v9) {${v_{9}}$};
       \coordinate (v10) at (.52, -2.94);\filldraw[blue] (v10) circle(1pt);               \node[below right] at (v10) {${v_{10}}$};

       \coordinate (v11) at (-.52, -2.94);\filldraw[blue] (v11) circle(1pt);\node[below left] at (v11) {${v_{11}}$};

       \coordinate (v12) at (-1.51, -2.61);\filldraw[blue] (v12) circle(1pt);\node[below left] at (v12) {${v_{12}}$};

       \coordinate (v13) at (-2.3,-1.93);\filldraw[blue] (v13) circle(1pt);\node[below left] at (v13) {${v_{13}}$};

       \coordinate (v14) at (-2.82, -1.02);\filldraw[blue] (v14) circle(1pt);\node[below left] at (v14) {${v_{14}}$};

       \coordinate (v15) at (-3,0);\filldraw[blue] (v15) circle(1pt);                  \node[left] at (v15) {${v_{15}}$};

       \coordinate (v16) at (-2.82, 1.02);\filldraw[blue] (v16) circle(1pt);\node[above left] at (v16) {${v_{16}}$};

       \coordinate (v17) at (-2.3,1.93);\filldraw[blue] (v17) circle(1pt);\node[above left] at (v17) {${v_{17}}$};

     \draw [line width=1.5,black] (v2) -- (v11);

    \draw [line width=1.,black](v6)-- (v11);

     \draw [line width=1.7,red](v17)-- (v4) -- (v13)--(v8) -- (v17);
     \draw [line width=1.7,blue](v16)-- (v3) -- (v12)--(v7) -- (v16);
     \draw [line width=2,green](v9)-- (v0) -- (v5)-- (v14)-- (v9);
     \draw [line width=1.8,orange](v6)--(v1)-- (v10)-- (v15)-- (v6);
     \draw [dashed, yellow](v6)--(v1)-- (v10)-- (v15)-- (v6);

  \draw[transform canvas={xshift=1.3pt}, blue] (v0)-- (v9);

  \draw[transform canvas={yshift= 1pt}, blue] (v5)-- (v14);

  \draw[transform canvas={yshift=1pt}, blue] (v0)-- (v5);

  \draw[transform canvas={yshift=1pt}, blue] (v9)-- (v14);

    \filldraw[black] (v0) circle(3.5pt);
    \filldraw[black] (v1) circle(3.5pt);
    \filldraw[black] (v2) circle(3.5pt);
   \filldraw[black] (v3) circle(3.5pt);
    \filldraw[black] (v4) circle(3.5pt);
   \filldraw[black] (v5) circle(3.5pt);
    \filldraw[black] (v6) circle(3.5pt);
    \filldraw[black] (v7) circle(3.5pt);
   \filldraw[black] (v8) circle(3.5pt);
    \filldraw[black] (v9) circle(3.5pt);
   \filldraw[black] (v10) circle(3.5pt);
    \filldraw[black] (v11) circle(3.5pt);
    \filldraw[black] (v12) circle(3.5pt);
    \filldraw[black] (v13) circle(3.5pt);
    \filldraw[black] (v14) circle(3.5pt);
    \filldraw[black] (v15) circle(3.5pt);
    \filldraw[black] (v16) circle(3.5pt);
    \filldraw[black] (v17) circle(3.5pt);

    \filldraw[green] (v0) circle(2pt);
    \filldraw[yellow] (v1) circle(2pt);
    \filldraw[white] (v2) circle(2pt);
    \filldraw[blue] (v3) circle(2pt);

    \filldraw[red] (v4) circle(2pt);
    \filldraw[green] (v5) circle(2pt);
    \filldraw[yellow] (v6) circle(2pt);
    \filldraw[blue] (v7) circle(2pt);
    \filldraw[red] (v8) circle(2pt);
    \filldraw[green] (v9) circle(2pt);

    \filldraw[yellow] (v10) circle(2pt);
    \filldraw[white] (v11) circle(2pt);
    \filldraw[blue] (v12) circle(2pt);
    \filldraw[red] (v13) circle(2pt);

    \filldraw[green] (v14) circle(2pt);
    \filldraw[yellow] (v15) circle(2pt);
    \filldraw[blue] (v16) circle(2pt);
    \filldraw[red] (v17) circle(2pt);

\end{tikzpicture}
}

\centering{(a) ${\cal J}_{18,6}(^*_8)\/$}

\end{center}
\end{minipage}
\hspace{4mm}
\begin{minipage}{.45\textwidth}
\begin{center}
\resizebox{6cm}{!}{
\begin{tikzpicture}
\coordinate (center) at (0,0);
   \def\radius{3cm}
   \foreach \x in {0, 20,...,360} {
             \filldraw[] (\x:3cm) circle(1pt);
               }


      \coordinate (v0) at (-1.51, 2.61);\filldraw[blue] (v0) circle(1pt);\node[above left] at (v0) {${v_0}$};

       \coordinate (v1) at (-.52, 2.94);\filldraw[blue] (v1) circle(1pt);\node[above left] at (v1) {${v_1}$};
      \coordinate (v2) at (.52, 2.94);\filldraw[blue] (v2) circle(1pt);\node[above right] at (v2) {${v_2}$};

       \coordinate (v3) at (1.51, 2.61);\filldraw[blue] (v3) circle(1pt);\node[above right] at (v3) {${v_{3}}$};
       \coordinate (v4) at (2.3,1.93);\filldraw[blue] (v4) circle(1pt);\node[above right] at (v4) {${v_4}$};

       \coordinate (v5) at (2.82, 1.02);\filldraw[blue] (v5) circle(1pt);\node[above right] at (v5) {${v_5}$};
       \coordinate (v6) at (3,0);\filldraw[blue] (v6) circle(1pt);\node[above right] at (v6) {${v_{6}}$};

       \coordinate (v7) at (2.82, -1.02);\filldraw[blue] (v7) circle(1pt);\node[below right] at (v7) {${v_7}$};

       \coordinate (v8) at (2.3,-1.93);\filldraw[blue] (v8) circle(1pt);\node[below right] at (v8) {${v_8}$};

       \coordinate (v9) at (1.51, -2.61);\filldraw[blue] (v9) circle(1pt);\node[below right] at (v9) {${v_{9}}$};
       \coordinate (v10) at (.52, -2.94);\filldraw[blue] (v10) circle(1pt);\node[below right] at (v10) {${v_{10}}$};

       \coordinate (v11) at (-.52, -2.94);\filldraw[blue] (v11) circle(1pt);\node[below left] at (v11) {${v_{11}}$};

       \coordinate (v12) at (-1.51, -2.61);\filldraw[blue] (v12) circle(1pt);\node[below left] at (v12) {${v_{12}}$};

       \coordinate (v13) at (-2.3,-1.93);\filldraw[blue] (v13) circle(1pt);\node[below left] at (v13) {${v_{13}}$};

       \coordinate (v14) at (-2.82, -1.02);\filldraw[blue] (v14) circle(1pt);\node[below left] at (v14) {${v_{14}}$};

       \coordinate (v15) at (-3,0);\filldraw[blue] (v15) circle(1pt);\node[left] at (v15) {${v_{15}}$};

       \coordinate (v16) at (-2.82, 1.02);\filldraw[blue] (v16) circle(1pt);\node[above left] at (v16) {${v_{16}}$};

       \coordinate (v17) at (-2.3,1.93);\filldraw[blue] (v17) circle(1pt);\node[above left] at (v17) {${v_{17}}$};

     \draw [line width=1.5,black](v2) -- (v11);

     \draw [line width=1.7,red](v17)-- (v4) -- (v13)--(v8) -- (v17);
     \draw [line width=1.7,blue](v16)-- (v3) -- (v12)--(v7) -- (v16);
     \draw [line width=1.9,green](v9)-- (v0) -- (v5)-- (v14)-- (v9);
     \draw [line width=1.8, orange](v6)--(v1)-- (v10)-- (v15)-- (v6);
     \draw [dashed, yellow](v6)--(v1)-- (v10)-- (v15)-- (v6);

  \draw[transform canvas={xshift=1.3pt}, blue] (v0)-- (v9);

  \draw[transform canvas={yshift= 1pt}, blue] (v5)-- (v14);

  \draw[transform canvas={yshift=1pt}, blue] (v0)-- (v5);

  \draw[transform canvas={yshift=1pt}, blue] (v9)-- (v14);

    \filldraw[black] (v0) circle(3.5pt);
    \filldraw[black] (v1) circle(3.5pt);
    \filldraw[black] (v2) circle(3.5pt);
   \filldraw[black] (v3) circle(3.5pt);
    \filldraw[black] (v4) circle(3.5pt);
   \filldraw[black] (v5) circle(3.5pt);
    \filldraw[black] (v6) circle(3.5pt);
    \filldraw[black] (v7) circle(3.5pt);
   \filldraw[black] (v8) circle(3.5pt);
    \filldraw[black] (v9) circle(3.5pt);
   \filldraw[black] (v10) circle(3.5pt);
    \filldraw[black] (v11) circle(3.5pt);
    \filldraw[black] (v12) circle(3.5pt);
    \filldraw[black] (v13) circle(3.5pt);
    \filldraw[black] (v14) circle(3.5pt);
    \filldraw[black] (v15) circle(3.5pt);
    \filldraw[black] (v16) circle(3.5pt);
    \filldraw[black] (v17) circle(3.5pt);

    \filldraw[green] (v0) circle(2pt);
    \filldraw[yellow] (v1) circle(2pt);
    \filldraw[white] (v2) circle(2pt);
    \filldraw[blue] (v3) circle(2pt);

    \filldraw[red] (v4) circle(2pt);
    \filldraw[green] (v5) circle(2pt);
    \filldraw[yellow] (v6) circle(2pt);
    \filldraw[blue] (v7) circle(2pt);
    \filldraw[red] (v8) circle(2pt);
    \filldraw[green] (v9) circle(2pt);

    \filldraw[yellow] (v10) circle(2pt);
    \filldraw[white] (v11) circle(2pt);
    \filldraw[blue] (v12) circle(2pt);
    \filldraw[red] (v13) circle(2pt);

    \filldraw[green] (v14) circle(2pt);
    \filldraw[yellow] (v15) circle(2pt);
    \filldraw[blue] (v16) circle(2pt);
    \filldraw[red] (v17) circle(2pt);

\end{tikzpicture}
}

\centering{(b) ${\cal J}_{18,6}(^*_9)\/$}

\end{center}
\end{minipage}
\caption{${\cal J}_{18,6}(^*_\gamma)\/$,  $t=3\/$,  and $\gamma\in\{8, 9\}\/$}     \label{J18,6 gamma}

\end{center}
\end{figure}
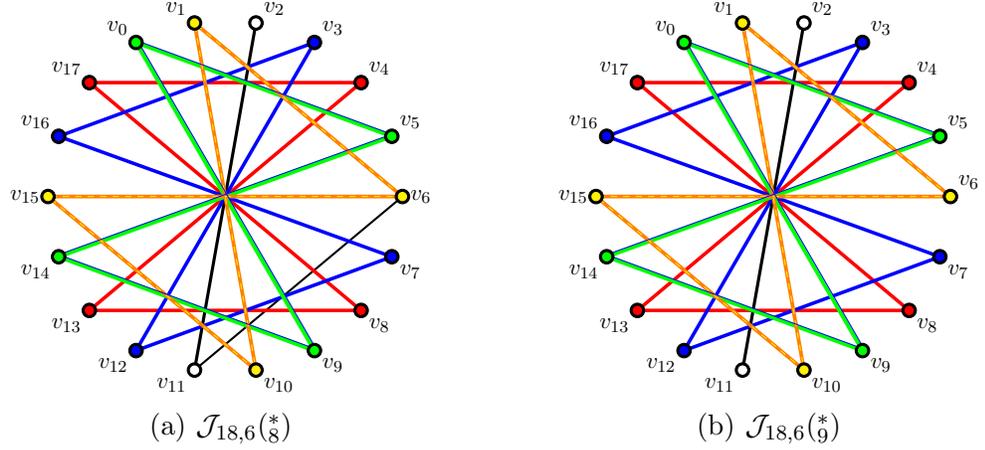

\vspace{5mm} An example of a convex geometric graph ${\cal J}_{n,g}(^*_\gamma)\/$  is depicted in Figure \ref{J18,6 gamma}.
 In Figure \ref{J18,6 gamma} (a), ${\cal J}_{18,6}(^*_8)\/$ has $\{v_{10},v_{15}\}\/$, $\{v_{10-\ell},v_{15-\ell}\}\/$ and $\{v_{6-\ell},v_{1-\ell}\}\/$  as $\alpha$-pairs where $\ell=1, 2, 3\/$,
 in Figure \ref{J18,6 gamma} (b) ${\cal J}_{18,6}(^*_9)\/$ has $\{v_{7+\ell},v_{12+\ell}\}\/$ and $\{v_{16+\ell},v_{3+\ell}\}\/$ as $\alpha$-pairs where $\ell=0, 1, 2, 3\/$.


\vspace{1mm}
\begin{result}   \label{n-g+1+mu}
Suppose  $g \geq 3\/$ is a natural number and let  $n=g+t(g-2)\/$ where  $t\geq 1$ is any natural number.    Suppose $F_n\/$ is a subgraph of the  convex complete graph $K_n\/$ such that  $n-g+2 \leq |E(F_n)| \leq n-1\/$  and $F_n\/$ contains no boundary edges  of $K_n\/$. Then $K_n - F_n\/$ admits a $g\/$-angulation    unless $F_n \/$ is  ${\cal J}_{n,g}(^*_\beta)\/$ for some  $\beta\in\{1, 2,  \ldots, 2g-3\}\/$.
\end{result}

\vspace{3mm}  \noindent
{\bf Proof:} In view of Propositions  \ref{Jstar1}, \ref{Jstar2} and \ref{Jstar-gamma}, we assume that $F_n \neq {\cal J}_{n,g}(^*_\beta)\/$ for any $\beta \geq 1\/$.

\vspace{5mm} The case $g=3\/$ has been treated in \cite{alta:refer}. Hence we assume that $g \geq 4\/$.

\vspace{5mm} We prove the result by induction on $t\/$.

\vspace{2mm} Suppose $t=1\/$. If $K_n - F_n\/$ contains an edge of the form $v_jv_{j+g-1}\/$, then $v_jv_{j+g-1}\/$ together with the boundary edges of $K_n\/$ is a $g\/$-angulation of  $K_n - F_n\/$.

\vspace{5mm}
On the other hand, if $v_jv_{j+g-1}\/$ is not an edge of $K_n - F_n\/$ for any $j =0 , 1, \ldots, g-2\/$, (that is, $v_jv_{j+g-1}\/$ is an edge of $F_n\/$),    then $F_n\/$ contains ${\cal F}_{n,g}(^*)\/$ as a subgraph.  This contradicts the assumption that $F_n \neq {\cal J}_{n,g}(^*_1)\/$.

\vspace{5mm} Now, assume that $t\geq2$  and the result is true for all convex graphs $K_m - F_m\/$ where $m= g+t'(g-2)\/$, $t'\leq t-1$.

\vspace{5mm} By Lemma \ref{earg}, $K_n - F_n\/$ contains an edge of the form $v_{j}v_{j+g-1}\/$ which is not an edge of $F_n\/$.
By relabeling if necessary, we may take $j=0$.

\vspace{5mm}  Now delete the set of vertices $v_1, \dots, v_{g-2}$ from    $K_n - F_n\/$.   Let $K_{m} - F_{m}\/$ denote the resulting convex graph. Here  $m= g+(t-1)(g-2)\/$.

\vspace{5mm} If $|E(F_{m})|\leq m-g$, then  $K_{m} - F_{m}\/$ admits a $g$-angulation by  Theorem \ref{atmostn-g}. Clearly, this $g\/$-angulation gives rise to   a  $g$-angulation for $K_n - F_n\/$.

\vspace{5mm} Hence we assume that $|E(F_{m})|\geq m-g+1$.

\vspace{5mm} Suppose $F_{m}\/$ is neither ${\cal F}_{m,g}(^*)\/$ nor ${\cal J}_{m,g}(^*_\beta)\/$ for any $\beta\/$. By Theorem \ref{n-2g} or by induction
$K_{m} - F_{m}\/$ admits a $g$-angulation. Again  this $g\/$-angulation gives rise to a $g$-angulation for $K_n - F_n\/$.

\vspace{5mm}  Hence assume that  $F_{m}\/$ is either ${\cal F}_{m,g}(^*)\/$ or ${\cal J}_{m,g}(^*_\beta)\/$ for some $\beta\in\{1, 2, \ldots 2g-3\}\/$.

\vspace{5mm}  Suppose that $\{v_1, \dots, v_{g-2}\}\/$ contains no vertex of an $\alpha\/$-pair in $F_n\/$.

\vspace{5mm} If $F_{m}={\cal F}_{m,g}(^*)\/$, then   $F_{n}\/$ is either ${\cal F}_{n,g}(^*_1)\/$ or ${\cal J}_{n,g}(^*_2)\/$. Either case is a contradiction.

\vspace{5mm} Hence assume that   $F_{m}= {\cal J}_{m,g}(^*_\beta)\/$.  If $\{v_1, \dots, v_{g-2}\}\/$ contains a neighbor of a vertex of an $\alpha\/$-pair in $F_n\/$, then  $F_{n}={\cal F}_{n,g}(^*_{\gamma})\/$  ($\gamma = \beta +1\/$).    If  $\{v_1, \dots, v_{g-2}\}\/$ contains no neighbor of a vertex of an   $\alpha\/$-pair in $F_n\/$, then   $F_{n}={\cal F}_{n,g}(^*_\beta)\/$. Either case is  a contradiction. \

\vspace{5mm} Now assume that  $\{v_r, v_s\}\/$ is an $\alpha\/$-pair in $F_n\/$ and   $v_r \in \{v_1, \dots, v_{g-2}\}\/$.

\vspace{5mm} If $v_r,  v_s\/$ are adjacent and  $|r - s|\equiv 1({\rm mod}\ (g-2))\/$ then $F_{n}= {\cal J}_{n,g}(^*_\beta)\/$ for some $\beta\in \{2, \ldots, 2g-3\}\/$,  a contradiction.  Hence either $v_r, v_s\/$ are non-adjacent or else $|r - s|\not\equiv 1({\rm mod}\ (g-2))\/$ .

\vspace{5mm} Assume without loss of generality that $r < s\/$. Let $i\/$ be  the largest integer such that $v_i\/$ is adjacent to $v_r\/$ in $F_n\/$,  and $j\/$ the smallest integer such that $v_j\/$ is adjacent to $v_s\/$ in $F_n\/$.

 \vspace{5mm} Let $k\/$ be an integer such that $ i <k < j\/$ and $|k-s| \equiv 1({\rm mod}\ (g-2))\/$. Then $v_r\/$ (respectively $v_s\/$) is a $g\/$-angulable vertex in the convex subgraph induced by $v_k, v_{k+1}, \ldots, v_s\/$ (respectively $v_s, v_{s+1}, \ldots, v_k\/$). This yields  a $g\/$-angulation for $K_n - F_n\/$.

\vspace{1mm} This completes the proof.    \qed

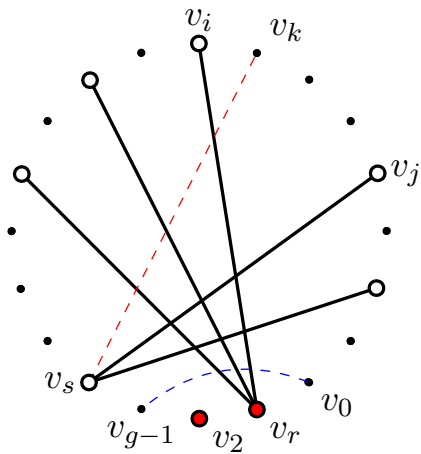
\begin{figure}[htp]
\begin{center}
\begin{minipage}{.9\textwidth}
\begin{center}
\resizebox{6cm}{!}{
\begin{tikzpicture}[rotate=-36,.style={draw}]
\coordinate (center) at (0,0);
   \def\radius{2cm}
   \foreach \x in {0, 18,...,360} {
             \filldraw[] (\x:2cm) circle(1pt);
               }


      \coordinate (v0) at (-1.62, 1.18);\filldraw[black] (v0) circle(1pt);

       \coordinate (v1) at (-1.18, 1.62);\filldraw[black] (v1) circle(1pt); \node[above] at (v1) {${v_{i}}$};
       \coordinate (v2) at (-.62, 1.9);\filldraw[black] (v2) circle(1pt);\node[above right] at (v2) {${v_{k}}$};

       \coordinate (v3) at (0, 2);\filldraw[black] (v3) circle(1pt);
       \coordinate (v4) at (.62, 1.9);\filldraw[black] (v4) circle(1pt);

       \coordinate (v5) at (1.18, 1.62);\filldraw[black] (v5) circle(1pt); \node[ right] at (v5) {${v_{j}}$};
       \coordinate (v6) at (1.62, 1.18);\filldraw[black] (v6) circle(1pt);

       \coordinate (v7) at (1.89, .62);\filldraw[black] (v7) circle(1pt);
      \coordinate (v8) at (2, 0);\filldraw[black] (v8) circle(1pt);

       \coordinate (v9) at (1.89, -.62);\filldraw[black] (v9) circle(1pt);\node[below right] at (v9) {${v_{0}}$};
       \coordinate (v10) at (1.62, -1.18);\filldraw[black] (v10) circle(1pt);\node[below right] at (v10) {${v_{r}}$};

       \coordinate (v11) at (1.18, -1.62);\filldraw[black] (v11) circle(1pt);\node[below right] at (v11) {${v_{2}}$};

       \coordinate (v12) at (.62, -1.9);\filldraw[black] (v12) circle(1pt);\node[below] at (v12) {${v_{g-1}}$};

       \coordinate (v13) at (0, -2);\filldraw[black] (v13) circle(1pt);\node[left] at (v13) {${v_{s}}$};

       \coordinate (v14) at (-.62, -1.9);\filldraw[black] (v14) circle(1pt);

       \coordinate (v15) at (-1.18, -1.62);\filldraw[black] (v15) circle(1pt);

       \coordinate (v16) at (-1.62, -1.18);\filldraw[black] (v16) circle(1pt);

       \coordinate (v17) at (-1.89, -.62);\filldraw[black] (v17) circle(1pt);

       \coordinate (v18) at (-2, 0);\filldraw[black] (v18) circle(1pt);

       \coordinate (v19) at (-1.89, .62);\filldraw[black] (v19) circle(1pt);

     \draw [line width=1,black](v10)-- (v17);
     \draw [line width=1,black](v10)-- (v19);
     \draw [line width=1,black](v10)-- (v1);

     \draw [line width=1,black](v13)-- (v5);
     \draw [line width=1,black](v13)-- (v7);

     \draw [dashed,red](v13)-- (v2);
     \draw [dashed,blue](v9) to [bend right] (v12);

    \filldraw[black] (v1) circle(2.5pt);
   \filldraw[black] (v5) circle(2.5pt);
    \filldraw[black] (v7) circle(2.5pt);
   \filldraw[black] (v10) circle(2.5pt);
   \filldraw[black] (v11) circle(2.5pt);
    \filldraw[black] (v13) circle(2.5pt);
    \filldraw[black] (v17) circle(2.5pt);
    \filldraw[black] (v19) circle(2.5pt);

    \filldraw[white] (v19) circle(1.5pt);
    \filldraw[white] (v1) circle(1.5pt);

    \filldraw[white] (v5) circle(1.5pt);
    \filldraw[white] (v7) circle(1.5pt);

    \filldraw[red] (v10) circle(1.5pt);
    \filldraw[red] (v11) circle(1.5pt);
    \filldraw[white] (v13) circle(1.5pt);
    \filldraw[white] (v17) circle(1.5pt);

\end{tikzpicture}
}

\end{center}
\end{minipage}
\caption{$g=4$, and $|k-s|=9 \equiv 1 \ ({\rm mod}\ (g-2))\/$}   \label{big crossing}
\end{center}
\end{figure}

\section{Potentially $g$-angulable graphs}     \label{pt}

We now look at the possibility of placing a graph $F\/$ with  $n\/$ vertices and $n\/$ edges in  the convex complete graph $K_n\/$ so that $K_n - F\/$ admits a $g$-angulation. We shall confine our attention to the case where $F_n\/$ is a $2\/$-regular graph.

\vspace{5mm}
\begin{define}{\bf:}  \label{potentiallyg}
Let $K_n\/$ be a convex complete graph with $n\/$ vertices. $F\/$ is said to be   potentially $g$-angulable  if there exists a configuration of $F\/$ in $K_n\/$ such that $K_n - F\/$ admits a $g$-angulation.
\end{define}

\vspace{5mm}
\begin{result}{\bf:}   \label{fnconnected}
 Suppose $F_n\/$ is an $n\/$-cycle and $g\geq 4$ is a natural number such that $n=g+t(g-2)$.  Then $F_n\/$ is potentially $g$-angulable  if and only if $n \geq 5\/$.
\end{result}

\vspace{2mm}  \noindent
{\bf Proof:}  It is easy to see that $K_n - F_n\/$ admits no $g$-angulation if $n \leq 4\/$.

\vspace{5mm} For $n=6\/$ suppose  $F\/$ is a  $6\/$-cycle. Let $F\/$ be of the form $v_0v_2v_4v_1v_5v_3v_0\/$, then $K_6 - F\/$ admits a $g$-angulation for $g\in\{4, 6\}$, when $g=4$ the diagonal is $v_2v_5\/$.

\vspace{5mm} For the rest of the proof, we assume that $n \geq 5\/$, where $n\neq 6$.

\vspace{5mm} When $n \/$ is odd, let $F_n\/$ takes the form \[ v_0 v_{2} v_{4}  \ldots v_{n-3} v_{n-1} v_1 v_3 v_5 \ldots v_{n-4} v_{n-2} v_0.  \]

\vspace{5mm} When  $n \/$ is even, let $F_n\/$ takes the form \[ v_0 v_{2} v_{4} \ldots v_{n-4} v_{n-2} v_1 v_3 v_5 \ldots v_{n-5} v_{n-1} v_{n-3} v_0.  \]

 \vspace{5mm} In both cases, the edges $v_2v_{3+i(g-2)}\/$, $i=1, \ldots, t\/$ together with the boundary  edges $v_0v_1 v_2 \ldots v_{n-1}v_0\/$ is a $g$-angulation of $K_n - F_n\/$.

 \vspace{5mm} This completes the proof. \qed

 \vspace{5mm}
 \begin{result}{\bf:} \label{disconnected}
Let $F_n\/$ be a $2\/$-regular graph with $n\/$ vertices and $g\geq 4$ is a natural number such that $n=g+t(g-2)$. Then $F_n\/$ is potentially $g$-angulable if and only if $n \geq 5\/$.
\end{result}

 \vspace{5mm} \noindent
{\bf Proof:} If $F_n\/$ is connected, the result is true by Theorem \ref{fnconnected}. Hence we assume that $F_n\/$ is a union of disjoint cycles.

 \vspace{5mm}  Let $C$ be a smallest cycle in $F_n$ and let $v_xv_y$ and $v_yv_z$ are two edges in $C$.
Consider $F^*\/$ to be a union of disjoint cycles in $F_n-C$ and let $|V(C)|=p$.

 \vspace{5mm}    (i) If $F^*$ is a $4$-cycle, Take $F^* = \{v_1v_3v_2v_4v_1\}\/$.
    Insert $v_x$, $v_y$ and $v_z$ of $C$ into edges $v_{4}v_1, v_1v_2, v_2v_3$ respectively. Then $K_7-F_7$ admits $7$-angulation with $t=0$.
     In case that $p=4$, insert the fourth vertex $v_w$ into the edge $v_3v_4$. Then $K_8-F_8$ admits $4$-angulation with two diagonal edges $v_yv_3$ and $v_yv_4$ and admits $5$-angulation with a diagonal edge $v_yv_w$.

 \vspace{5mm}    (ii) If $F^*$ is not a $4$-cycle,  then place $F^*$ on $K^*_{n-p}$ so that $F^*$ contains no boundary edge of  $K^*_{n-p}$.
 \vspace{5mm}   Let $v_1, v_2, \ldots, v_{n-p}$ denotes the vertices of  $K^*_{n-p}$.
Insert $v_x$, $v_y$ and $v_z$ of $C$ into edges $v_1v_2, v_2v_3, v_3v_4$ respectively.
In case that $p>3$ insert the rest of vertices of $C$, which are $p-3$ vertices, into the edges of the path $v_4v_5 \cdots v_{n-p}v_1$ and place $C$ on $K_{n}$ so that $C$ contains no boundary edge of  $K_{n}$.

 \vspace{5mm}   Relabel the vertices of $K_n$ to be $u_0, u_1, u_2, \ldots, u_{n-1}$ with $v_y=u_0$.
Hence, we have $u_0u_{1+i(g-2)}$, $i=1, 2, \dots, t$ together with the edges $u_0u_1 u_2u_3 \ldots u_{n-1}u_0\/$ is a $g$-angulation of $K_n - F_n\/$ (since $u_0$ is adjacent only to $u_{n-2}$ and $u_{2}$).

 \vspace{5mm}   This completes the proof.   $\qed$

\section{Regular graphs}  \label{regular}

In view of the results in the preceding section, it is natural to ask which regular graph is potentially  $g$-angulable in $K_n\/$.

 \vspace{5mm}
{\bf Problem:} \
Let $r \geq 3\/$  and  $g \geq 3\/$   be two natural numbers and let $G\/$ be an $r\/$-regular graph with $n\/$ vertices   where $n = g + (g-2)t\/$. It is true that there is a natural number  $n_0(r, g)\/$ such that when $n \geq n_0(r, g)\/$, then $G\/$ is potentially $g$-angulable in the convex complete graph $K_n\/$?

 \vspace{5mm}   We believe that the above problem is true. However we do not have a complete  answer for this even when restricted to the case $r=3\/$.  Nevertheless  we offer the following special case of a $3\/$-regular graph which is well-known in the literature.

 \vspace{5mm}
Suppose $n\/$ and $k\/$ are two integers such that $1 \leq k \leq n-1\/$ and $n \geq 5\/$. The {\em generalized Petersen graph\/} $P(n,k)\/$ is defined to have vertex-set
  $\{a_i, b_i : i = 0, 1, \ldots, n-1\}\/$ and edge-set $E_1 \cup E_2 \cup E_3\/$ where   $E_1 = \{ a_ia_{i+1} : i = 0, 1, \ldots, n-1\}\/$,
$E_2 = \{ b_ib_{i+k} : i = 0, 1, \ldots, n-1\}\/$   and $E_3 = \{ a_ib_i : i = 0, 1, \ldots, n-1\}\/$     with subscripts reduced modulo $n \/$. Edges in $E_3\/$ are called the spokes of $P(n,k)\/$.

\begin{propo}{\bf:}  \label{petersen}
Suppose $ 1 \leq k < n/2\/$ and $g\geq 4$ is a natural number such that $2n=g+t(g-2)$. Then the generalized Petersen graph $P(n, k)\/$ is potentially $g$-angulable in the convex complete graph $K_{2n}\/$  where  $n \geq 5\/$.
\end{propo}

 \vspace{5mm}    \noindent
{\bf Proof:}  Let the vertices of $K_{2n}\/$ be denoted  $v_1, v_2, v_3, \ldots, v_{2n}\/$.  We shall pack $P(n,k)\/$ on $K_{2n}\/$ so that $K_{2n} - P(n,k)\/$ admits a $g$-angulation.

 \vspace{5mm}   {\em Case (1): \ $k=1\/$}

 \vspace{5mm}   $P(n,1)\/$ consists of two $n\/$-cycles $C = a_0a_1a_2 \cdots c_{n-1}a_0\/$ and $C' = b_0b_1b_2 \cdots  b_{n-1}b_0\/$ together with the edges $a_ib_i\/$, $i=0, 1, 2, \ldots, n-1\/$.

 \vspace{5mm}    Place  $C\/$ on $K_{2n}\/$ so that $C\/$ takes the form   $v_2v_{4}v_6 \cdots v_{2n-4}v_{2n-2}v_{2n}v_2\/$
              and that $C'\/$ takes the form   $v_1v_{3}v_5v_{4}v_7 \cdots v_{2n-3}v_{2n-1}v_1\/$.

 \vspace{5mm}   {\em Case (2): \ $1 < k < n/2\/$}

 \vspace{5mm}   $P(n,k)\/$ consists of two $n\/$-cycles $C = a_0a_1a_2 \cdots c_{n-1}a_0\/$ and $C' = \{ b_ib_{i+k}, \ i =0, 1, 2,  \ldots, n-1\}\/$ together with the edges $a_ib_i\/$, $i=0, 1, 2, \ldots, n-1\/$.

 \vspace{5mm}   Place $C\/$ and $C'\/$ on $K_{2n}\/$ so that  $C\/$ takes the form    $v_2v_4v_6 \ldots v_{2n-2}v_{2n}v_2\/$
             and $C'\/$ takes the form  $\{ v_iv_{i+2k} \ : \ i =1, 3,  \ldots, 2n-1\}  \/$. The operations on the subscripts are reduced modulo $2n\/$.

 \vspace{5mm}   In both cases: Let the spokes take the form  \ $v_iv_{i+3}\/$, \ $i =0, 2, 4, \ldots, 2n-2\/$. Here also the operations on the subscripts are reduced modulo $2n\/$.

 \vspace{5mm}   In both cases: if $t=1$, then $K_n - F_n\/$ has a $g$-angulation whose diagonal is $v_1v_g$.
If $t\geq2$,  consider the subgraph $H\/$ induced by the sets of vertices $\{v_g, v_{g+1}, \ldots, \linebreak v_{2n-(g-1)}\} \cup \{ v_0, v_1\}\/$.
 Since the vertex $v_0\/$ is not adjacent to every vertex in $H\/$ and $|V(H)|=2n-2(g-2)$, $K_{2n-2(g-2)}- H\/$ admits a $g$-angulation $G\/$. Then $G \cup \{v_1v_2 \cdots v_g \}\cup \{ v_{2n-(g-1)}v_{2n-(g-1)+1} 
  \ldots v_{2n-1}v_0\}\/$ is a $g$-angulation for $K_{2n} - P(n,k)\/$.

 \vspace{5mm}   This completes the proof.  $\qed$

\begin{figure}[htp]
\begin{minipage}{.32\textwidth}
\begin{center}
\resizebox{5cm}{!}{
\begin{tikzpicture}
\coordinate (center) at (0,0);
  \def\radius{2.5cm}
   \foreach \x in {0, 51.4285714,...,360} {
             \filldraw[] (\x:2.cm) circle(1pt);

             \filldraw[] (\x:1.2cm) circle(1pt);

               }


      \coordinate (v0) at (-0.27,1.16);\filldraw[black] (v0) circle(2.5pt);\node[above right] at (v0) {\red${b_6}$};
       \coordinate (v1) at (-0.45,1.96);\filldraw[black] (v1) circle(2.5pt);\node[above] at (v1) {\blue${a_6}$};
       \coordinate (v2) at (0.75,.95);\filldraw[black] (v2) circle(2.5pt);\node[ right] at (v2) {\red${b_0}$};
       \coordinate (v3) at (1.26,1.58);\filldraw[black] (v3) circle(2.5pt);\node[above right] at (v3) {\blue${a_0}$};
       \coordinate (v4) at (1.2,0.);\filldraw[black] (v4) circle(2.5pt);\node[below] at (v4) {\red${b_1}$};
       \coordinate (v5) at (2,0);\filldraw[black] (v5) circle(2.5pt);\node[right] at (v5) {\blue${a_1}$};
       \coordinate (v6) at (.75,-0.95);\filldraw[black] (v6) circle(2.5pt);\node[below] at (v6) {\red${b_2}$};
       \coordinate (v7) at (1.26,-1.58);\filldraw[black] (v7) circle(2.5pt);\node[below right] at (v7) {\blue${a_2}$};
       \coordinate (v8) at (-0.27,-1.16);\filldraw[black] (v8) circle(2.5pt);\node[left] at (v8) {\red $b_3$};
       \coordinate (v9) at (-0.45,-1.96);\filldraw[black] (v9) circle(2.5pt);\node[below left] at (v9) {\blue${a_3}$};
       \coordinate (v10) at (-1.08,-.53);\filldraw[black] (v10) circle(2.5pt);\node[above left] at (v10) {\red${b_4}$};
        \coordinate (v11) at (-1.8,-0.87);\filldraw[black] (v11) circle(2.5pt);\node[left] at (v11) {\blue${a_4}$};
       \coordinate (v12) at (-1.08,.53);\filldraw[black] (v12) circle(2.5pt);\node[above] at (v12) {\red${b_5}$};
       \coordinate (v13) at (-1.8,0.87);\filldraw[black] (v13) circle(2.5pt);\node[above left] at (v13) {\blue${a_5}$};


       \draw [line width=1.2,blue](v1) -- (v3) -- (v5) -- (v7) -- (v9) -- (v11) -- (v13) -- (v1) ;

      \draw [line width=1.2,purple](v0) -- (v6) -- (v12) -- (v4) -- (v10) -- (v2) -- (v8) -- (v0) ;

       \draw [line width=1,black](v0) -- (v1);
       \draw [line width=1,black](v2) -- (v3);
       \draw [line width=1,black](v4) -- (v5);
       \draw [line width=1,black](v6) -- (v7);
       \draw [line width=1,black](v8) -- (v9);
       \draw [line width=1,black](v10) -- (v11);
       \draw [line width=1,black](v12) -- (v13);

\filldraw[black] (v0) circle(2.8pt);
\filldraw[black] (v1) circle(2.8pt);
\filldraw[black] (v2) circle(2.8pt);
\filldraw[black] (v3) circle(2.8pt);
\filldraw[black] (v4) circle(2.8pt);
\filldraw[black] (v5) circle(2.8pt);
\filldraw[black] (v6) circle(2.8pt);
\filldraw[black] (v7) circle(2.8pt);
\filldraw[black] (v8) circle(2.8pt);
\filldraw[black] (v9) circle(2.8pt);
\filldraw[black] (v10) circle(2.8pt);
\filldraw[black] (v11) circle(2.8pt);
\filldraw[black] (v12) circle(2.8pt);
\filldraw[black] (v13) circle(2.8pt);

\filldraw[red] (v0) circle(1.5pt);
\filldraw[blue] (v1) circle(1.5pt);
\filldraw[red] (v2) circle(1.5pt);
\filldraw[blue] (v3) circle(1.5pt);
\filldraw[red] (v4) circle(1.5pt);
\filldraw[blue] (v5) circle(1.5pt);
\filldraw[red] (v6) circle(1.5pt);
\filldraw[blue] (v7) circle(1.5pt);
\filldraw[red] (v8) circle(1.5pt);
\filldraw[blue] (v9) circle(1.5pt);
\filldraw[red] (v10) circle(1.5pt);
\filldraw[blue] (v11) circle(1.5pt);
\filldraw[red] (v12) circle(1.5pt);
\filldraw[blue] (v13) circle(1.5pt);

\end{tikzpicture}
}
\end{center}

\end{minipage}
\begin{minipage}{.33\textwidth}
\begin{center}
\resizebox{5.5cm}{!}{
\begin{tikzpicture}[rotate=-14,.style={draw}]
\coordinate (center) at (0,0);
  \def\radius{2.5cm}
   \foreach \x in {0,25.71428,...,360} {
             \filldraw[] (\x:2cm) circle(1pt);
               }


       \coordinate (v0) at (-1.26,1.58);\filldraw[black] (v0) circle(2.5pt);\node[above left] at (v0) {\blue${a_0}$};
       \coordinate (v1) at (-0.45,1.96);\filldraw[black] (v1) circle(2.5pt);\node[above] at (v1) {\red${b_6}$};
       \coordinate (v2) at (0.45,1.96);\filldraw[black] (v2) circle(2.5pt);\node[above right] at (v2) {\blue${a_1}$};
       \coordinate (v3) at (1.26,1.58);\filldraw[black] (v3) circle(2.5pt);\node[above right] at (v3) {\red${b_0}$};
       \coordinate (v4) at (1.8,0.87);\filldraw[black] (v4) circle(2.5pt);\node[right] at (v4) {\blue${a_2}$};
       \coordinate (v5) at (2,0);\filldraw[black] (v5) circle(2.5pt);\node[right] at (v5) {\red${b_1}$};
       \coordinate (v6) at (1.8,-0.87);\filldraw[black] (v6) circle(2.5pt);\node[below right] at (v6) {\blue${a_3}$};
       \coordinate (v7) at (1.26,-1.58);\filldraw[black] (v7) circle(2.5pt);\node[below right] at (v7) {\red${b_2}$};
       \coordinate (v8) at (0.45,-1.96);\filldraw[black] (v8) circle(2.5pt);\node[below] at (v8) {\blue $a_4$};
       \coordinate (v9) at (-0.45,-1.96);\filldraw[black] (v9) circle(2.5pt);\node[below left] at (v9) {\red${b_3}$};
        \coordinate (v10) at (-1.26,-1.58);\filldraw[black] (v10) circle(2.5pt);\node[below left] at (v10) {\blue${a_5}$};
        \coordinate (v11) at (-1.8,-0.87);\filldraw[black] (v11) circle(2.5pt);\node[left] at (v11) {\red${b_4}$};
        \coordinate (v12) at (-2,0);\filldraw[black] (v12) circle(2.5pt);\node[ left] at (v12) {\blue${a_6}$};
       \coordinate (v13) at (-1.8,0.87);\filldraw[black] (v13) circle(2.5pt);\node[above left] at (v13) {\red${b_5}$};


       \draw [line width=1.2,blue](v0) to [bend right] (v2) to [bend right] (v4) to [bend right] (v6) to [bend right] (v8) to [bend right] (v10) to [bend right] (v12) to [bend right] (v0) ;

       \draw [line width=1.2,purple](v1) -- (v7) -- (v13) -- (v5) -- (v11) -- (v3) -- (v9) -- (v1) ;

       \draw [line width=1,black](v0)  to [bend right] (v3);
       \draw [line width=1,black](v2)  to [bend right] (v5);
       \draw [line width=1,black](v4)  to [bend right] (v7);
       \draw [line width=1,black](v6)  to [bend right] (v9);
       \draw [line width=1,black](v8)  to [bend right] (v11);
       \draw [line width=1,black](v10)  to [bend right] (v13);
       \draw [line width=1,black](v12)  to [bend right] (v1);

\filldraw[black] (v0) circle(2.8pt);
\filldraw[black] (v1) circle(2.8pt);
\filldraw[black] (v2) circle(2.8pt);
\filldraw[black] (v3) circle(2.8pt);
\filldraw[black] (v4) circle(2.8pt);
\filldraw[black] (v5) circle(2.8pt);
\filldraw[black] (v6) circle(2.8pt);
\filldraw[black] (v7) circle(2.8pt);
\filldraw[black] (v8) circle(2.8pt);
\filldraw[black] (v9) circle(2.8pt);
\filldraw[black] (v10) circle(2.8pt);
\filldraw[black] (v11) circle(2.8pt);
\filldraw[black] (v12) circle(2.8pt);
\filldraw[black] (v13) circle(2.8pt);

\filldraw[blue] (v0) circle(1.5pt);
\filldraw[red] (v1) circle(1.5pt);
\filldraw[blue] (v2) circle(1.5pt);
\filldraw[red] (v3) circle(1.5pt);
\filldraw[blue] (v4) circle(1.5pt);
\filldraw[red] (v5) circle(1.5pt);
\filldraw[blue] (v6) circle(1.5pt);
\filldraw[red] (v7) circle(1.5pt);
\filldraw[blue] (v8) circle(1.5pt);
\filldraw[red] (v9) circle(1.5pt);
\filldraw[blue] (v10) circle(1.5pt);
\filldraw[red] (v11) circle(1.5pt);
\filldraw[blue] (v12) circle(1.5pt);
\filldraw[red] (v13) circle(1.5pt);

\end{tikzpicture}
}
\end{center}


\end{minipage}
\hspace{1mm}
\begin{minipage}{.32\textwidth}
\begin{center}
\resizebox{5.5cm}{!}{
\begin{tikzpicture}[rotate=-14,.style={draw}]
\coordinate (center) at (0,0);
  \def\radius{2.5cm}
   \foreach \x in {0,25.71428,...,360} {
             \filldraw[] (\x:2cm) circle(1pt);
               }


       \coordinate (v0) at (-1.26,1.58);\filldraw[black] (v0) circle(2.5pt);\node[above left] at (v0) {${v_0}$};
       \coordinate (v1) at (-0.45,1.96);\filldraw[black] (v1) circle(2.5pt);\node[above] at (v1) {${v_1}$};
       \coordinate (v2) at (0.45,1.96);\filldraw[black] (v2) circle(2.5pt);\node[above right] at (v2) {${v_2}$};
       \coordinate (v3) at (1.26,1.58);\filldraw[black] (v3) circle(2.5pt);\node[above right] at (v3) {${v_{3}}$};
       \coordinate (v4) at (1.8,0.87);\filldraw[black] (v4) circle(2.5pt);\node[right] at (v4) {${v_4}$};
       \coordinate (v5) at (2,0);\filldraw[black] (v5) circle(2.5pt);\node[right] at (v5) {${v_5}$};
       \coordinate (v6) at (1.8,-0.87);\filldraw[black] (v6) circle(2.5pt);\node[below right] at (v6) {${v_{6}}$};
       \coordinate (v7) at (1.26,-1.58);\filldraw[black] (v7) circle(2.5pt);\node[below right] at (v7) {${v_7}$};
       \coordinate (v8) at (0.45,-1.96);\filldraw[black] (v8) circle(2.5pt);\node[below] at (v8) {${v_8}$};
       \coordinate (v9) at (-0.45,-1.96);\filldraw[black] (v9) circle(2.5pt);\node[below left] at (v9) {${v_9}$};
        \coordinate (v10) at (-1.26,-1.58);\filldraw[black] (v10) circle(2.5pt);\node[below left] at (v10) {${v_{10}}$};
        \coordinate (v11) at (-1.8,-0.87);\filldraw[black] (v11) circle(2.5pt);\node[left] at (v11) {${v_{11}}$};
        \coordinate (v12) at (-2,0);\filldraw[black] (v12) circle(2.5pt);\node[ left] at (v12) {${v_{12}}$};
       \coordinate (v13) at (-1.8,0.87);\filldraw[black] (v13) circle(2.5pt);\node[above left] at (v13) {${v_{13}}$};



       \draw [line width=1.2,blue](v0) to [bend right] (v2) to [bend right] (v4) to [bend right] (v6) to [bend right] (v8) to [bend right] (v10) to [bend right] (v12) to [bend right] (v0) ;

       \draw [line width=1.2,purple](v1) -- (v7) -- (v13) -- (v5) -- (v11) -- (v3) -- (v9) -- (v1) ;

       \draw [line width=1,black](v0)  to [bend right] (v3);
       \draw [line width=1,black](v2)  to [bend right] (v5);
       \draw [line width=1,black](v4)  to [bend right] (v7);
       \draw [line width=1,black](v6)  to [bend right] (v9);
       \draw [line width=1,black](v8)  to [bend right] (v11);
       \draw [line width=1,black](v10)  to [bend right] (v13);
       \draw [line width=1,black](v12)  to [bend right] (v1);

\draw [decorate,decoration={snake, segment length=1.3mm, amplitude=.3mm}, line width=.03,red] (v0) --  (v1) to [bend right] (v4) --  (v5)  -- (v6)  -- (v7)  -- (v8)  -- (v9)  -- (v10)  -- (v11) to [bend right] (v0);

\filldraw[black] (v0) circle(2.8pt);
\filldraw[black] (v1) circle(2.8pt);
\filldraw[black] (v2) circle(2.8pt);
\filldraw[black] (v3) circle(2.8pt);
\filldraw[black] (v4) circle(2.8pt);
\filldraw[black] (v5) circle(2.8pt);
\filldraw[black] (v6) circle(2.8pt);
\filldraw[black] (v7) circle(2.8pt);
\filldraw[black] (v8) circle(2.8pt);
\filldraw[black] (v9) circle(2.8pt);
\filldraw[black] (v10) circle(2.8pt);
\filldraw[black] (v11) circle(2.8pt);
\filldraw[black] (v12) circle(2.8pt);
\filldraw[black] (v13) circle(2.8pt);

\filldraw[blue] (v0) circle(1.5pt);
\filldraw[red] (v1) circle(1.5pt);
\filldraw[blue] (v2) circle(1.5pt);
\filldraw[red] (v3) circle(1.5pt);
\filldraw[blue] (v4) circle(1.5pt);
\filldraw[red] (v5) circle(1.5pt);
\filldraw[blue] (v6) circle(1.5pt);
\filldraw[red] (v7) circle(1.5pt);
\filldraw[blue] (v8) circle(1.5pt);
\filldraw[red] (v9) circle(1.5pt);
\filldraw[blue] (v10) circle(1.5pt);
\filldraw[red] (v11) circle(1.5pt);
\filldraw[blue] (v12) circle(1.5pt);
\filldraw[red] (v13) circle(1.5pt);

\end{tikzpicture}
}
\end{center}


\end{minipage}
\caption{$P(7,3)$}    \label{petersen}

\end{figure}

\vspace{1mm}  In this part, we shall prove that $3\/$-regular graph with $n =4 + 2t\/$ vertices where $t \geq 2\/$ is potentially $4$-angulable.

\begin{support}{\bf:} \label{4labeling}
Let $G\/$ be a $3\/$-regular graph with $n =4 + 2t\/$ vertices where $t \geq 2\/$.
Suppose $V(G)$ can be labeled as  $v_1, v_2, \ldots, v_{n}$ such that

 (i) $v_1v_a, v_1v_{b}, v_1v_{c} \in E(G)$  where $a,b, c\/$ are distinct odd integers (different from $1\/$,    and

 (ii)   $|i-j|\notin \{1, n-1\}$ whenever $v_iv_j \in E(G)$.

 \noindent  Then $G\/$ is potentially $4$-angulable.

\end{support}

 \vspace{5mm}    \noindent
{\bf Proof:} Place the vertices  $v_1, v_2, \ldots, v_{n}\/$ of $G\/$ in convex position and put them in clock wise order. Then
 $v_1v_{4+2k}\/$, $k=0, 1, \ldots, t-1$ together with all boundary edges of $K_n\/$ is a $4$-angulation of $K_n-G$.      $\qed$

\vspace{7mm}
In the next lemma, we shall show that all $3\/$-regular graphs with at least $8\/$ vertices admit a labeling as described in Lemma \ref{4labeling} unless it is the cube $Q_3\/$ (on $8\/$ vertices) in which case it has a labeling as shown in Figure \ref{cube Q3}.   If $Q_3\/$, with the given labeling, is placed on convex  position (with clockwise order), then $v_2v_5, v_1v_6\/$ together with the boundary edges of $K_8\/$ is a $4$-angulation of $K_8-Q_3$.  This implies that all connected $3\/$-regular graphs are potentially $4\/$-angulable.

\vspace{5mm} To facilitate the proof of Theorem \ref{4-angulation}, we shall need to consider $3\/$-regular graphs where double edges are allowed. By a {\em $2\/$-cycle\/} in a graph, written $uvu\/$, we mean two edges of the graph of the form $uv\/$ and $vu\/$.

\vspace{5mm}
Let $G\/$ be a $3\/$-regular graph on $n\/$ vertices and let $e=xy$ be an edge in $G$. 
Suppose  $N_{G}(x)=\{x_1, x_2, y\}$ and $N_{G}(y)=\{y_1, y_2, x\}$.
Let $G_e\/$ be the  graph obtained from $G-e\/$ by replacing the paths $x_1xx_2\/$ and $y_1yy_2\/$ with the edges  $x_1x_2\/$ and  $y_1y_2\/$ respectively.   Then $G_e\/$ is a $3\/$-regular graph with $n-2\/$ vertices.

\vspace{5mm} In the case that $G\/$ has no $2\/$-cycle, it is easy to see that $G\/$ contains an edge $e\/$ such that $G_e\/$ has at most one $2\/$-cycle if $n \geq 6\/$.

 \vspace{5mm}
\begin{support}{\bf:} \label{labeling}
Let $G\/$ be a connected $3\/$-regular  graph with $n\geq 8\/$ vertices having at most one $2\/$-cycle.

(i) Suppose $G\/$ has no multiple edges. Then $V(G)\/$ admits a labeling as described in Lemma \ref{4labeling} unless $G\/$ is the cube which has labeling as shown in Figure \ref{cube Q3}.

(ii) Suppose $G\/$ has a $2\/$-cycle $uvu\/$.  Then $V(G)\/$ can be labeled as $v_1, v_2, \ldots, v_n\/$ such that $u=v_1, v =v_a\/$, $v_1v_b \in E(G)\/$ where $a, b\/$ are distinct odd integers different from $1\/$, and  that  $|i-j|\notin \{1, n-1\}$ whenever $v_iv_j \in E(G)$.
\end{support}

 \vspace{5mm}    \noindent
{\bf Proof:}
We prove this lemma by induction on $n\/$.

\begin{figure}[htp]
\begin{center}

\begin{minipage}{.3\textwidth}            

\resizebox{4.5cm}{!}{
\begin{tikzpicture}
\coordinate (center) at (0,0);
   \def\radius{2.2cm}
   \foreach \x in {0, 60,...,360} {
               }

\draw [line width=1.5, white] plot [smooth,  tension=.9] coordinates {(0,4) (-2,3) (-2,0) (0,-1)};

       \coordinate (v1) at (-1.2, 0);\filldraw[black] (v1) circle(4pt);    \node[left] at (v1)    {\LARGE${v_2}$};
       \coordinate (v2) at (1.2,0);\filldraw[black] (v2) circle(4pt);      \node[right] at (v2)   {\LARGE${v_{6}}$};
       \coordinate (v3) at (0.2, 1);\filldraw[black] (v3) circle(4pt);    \node[above] at (0.3, 1)    {\LARGE${v_8}$};
       \coordinate (v4) at (-0.2, 2);\filldraw[black] (v4) circle(4pt);     \node[below] at (-0.3, 2)   {\LARGE${v_5}$};
       \coordinate (v5) at (-1.2, 3);\filldraw[black] (v5) circle(4pt);    \node[left] at (v5)    {\LARGE${v_7}$};
       \coordinate (v6) at (1.2, 3);\filldraw[black] (v6) circle(4pt);     \node[right] at (v6)   {\LARGE${v_3}$};

      \coordinate (v9) at (0.2, 4);  \filldraw[black] (v9) circle(4pt);     \node[right] at (v9)     {\LARGE${v_1}$};
      \coordinate (v10) at (-0.2,-1);    \filldraw[black] (v10) circle(4pt);     \node[right] at (v10) {\LARGE${v_{4}}$};

     \draw [line width=1.5,black](v10) -- (v1) -- (v5) -- (v9) -- (v6) -- (v2) -- (v10);
   \draw [line width=1.5,black] (v6) -- (v4)-- (v5);
     \draw [line width=1.5,black](v1) -- (v3) -- (v2);
     \draw [line width=1.5,black](v10) -- (v3) ;
     \draw [line width=1.5,black](v4) -- (v9);

       \filldraw[white] (-1.2, 0) circle(2pt);
       \filldraw[white] (1.2,0) circle(2pt);
       \filldraw[white] (0.2, 1) circle(2pt);
       \filldraw[white] (-0.2, 2) circle(2pt);
       \filldraw[white] (-1.2, 3) circle(2pt);
       \filldraw[white] (1.2, 3) circle(2pt);

       \filldraw[white] (0.2, 4) circle(2pt);
       \filldraw[white] (-0.2,-1) circle(2pt);

\end{tikzpicture}
}

\centering  1

\end{minipage}
\begin{minipage}{.3\textwidth}              
\resizebox{4.5cm}{!}{
\begin{tikzpicture}
\coordinate (center) at (0,0);
   \def\radius{2.2cm}
   \foreach \x in {0, 60,...,360} {
               }


       \coordinate (v1) at (-1, 0);\filldraw[black] (v1) circle(4pt);    \node[left] at (v1)    {\LARGE${v_4}$};
       \coordinate (v2) at (1,0);\filldraw[black] (v2) circle(4pt);      \node[right] at (v2)   {\LARGE${v_{6}}$};
       \coordinate (v3) at (-1, 1.5);\filldraw[black] (v3) circle(4pt);    \node[left] at (v3)    {\LARGE${v_8}$};
       \coordinate (v4) at (1, 1.5);\filldraw[black] (v4) circle(4pt);     \node[right] at (v4)   {\LARGE${v_3}$};
       \coordinate (v5) at (-1, 3);\filldraw[black] (v5) circle(4pt);    \node[left] at (v5)    {\LARGE${v_5}$};
       \coordinate (v6) at (1, 3);\filldraw[black] (v6) circle(4pt);     \node[right] at (v6)   {\LARGE${v_1}$};

      \coordinate (v9) at (0, 4);  \filldraw[black] (v9) circle(4pt);     \node[right] at (v9)     {\LARGE${v_7}$};
      \coordinate (v10) at (0,-1);    \filldraw[black] (v10) circle(4pt);     \node[right] at (v10) {\LARGE${v_{2}}$};

     \draw [line width=1.5,black](v10) -- (v1) -- (v3) -- (v5) -- (v9) -- (v6) -- (v4) -- (v2) -- (v10);
     \draw [line width=1.5,black](v3) -- (v4);
     \draw [line width=1.5,black](v5) -- (v6);
     \draw [line width=1.5,black](v1) -- (v2);

\draw [line width=1.5, black] plot [smooth,  tension=.9] coordinates {(0,4) (-2,3) (-2,0) (0,-1)};

       \filldraw[white] (-1, 0) circle(2pt);
       \filldraw[white] (1, 0) circle(2pt);
       \filldraw[white] (-1, 1.5) circle(2pt);
       \filldraw[white] (1, 1.5) circle(2pt);
       \filldraw[white] (-1, 3) circle(2pt);
       \filldraw[white] (1, 3) circle(2pt);
       \filldraw[white] (1, 4) circle(2pt);

       \filldraw[white] (0, 4) circle(2pt);
       \filldraw[white] (0,-1) circle(2pt);

\end{tikzpicture}
}

\centering  2

\end{minipage}
\begin{minipage}{.3\textwidth}              
\resizebox{4.5cm}{!}{
\begin{tikzpicture}
\coordinate (center) at (0,0);
   \def\radius{2.2cm}
   \foreach \x in {0, 60,...,360} {
               }


       \coordinate (v1) at (-1, 0);\filldraw[black] (v1) circle(4pt);    \node[left] at (v1)    {\LARGE${v_6}$};
       \coordinate (v2) at (1,0);\filldraw[black] (v2) circle(4pt);      \node[right] at (v2)   {\LARGE${v_{8}}$};
       \coordinate (v3) at (-1, 1.5);\filldraw[black] (v3) circle(4pt);    \node[left] at (v3)    {\LARGE${v_4}$};
       \coordinate (v4) at (1, 1.5);\filldraw[black] (v4) circle(4pt);     \node[right] at (v4)   {\LARGE${v_3}$};
       \coordinate (v5) at (-1, 3);\filldraw[black] (v5) circle(4pt);    \node[left] at (v5)    {\LARGE${v_7}$};
       \coordinate (v6) at (1, 3);\filldraw[black] (v6) circle(4pt);     \node[right] at (v6)   {\LARGE${v_1}$};

      \coordinate (v9) at (0, 4);  \filldraw[black] (v9) circle(4pt);     \node[right] at (v9)     {\LARGE${v_5}$};
      \coordinate (v10) at (0,-1);    \filldraw[black] (v10) circle(4pt);     \node[right] at (v10) {\LARGE${v_{2}}$};

     \draw [line width=1.5,black](v10) -- (v1) -- (v3) -- (v5) -- (v9) -- (v6) -- (v4) -- (v2) -- (v10);
     \draw [line width=1.5,black](v1) -- (v4);
     \draw [line width=1.5,black](v5) -- (v6);
     \draw [line width=1.5,black](v3) -- (v2);

\draw [line width=1.5, black] plot [smooth,  tension=.9] coordinates {(0,4) (-2,3) (-2,0) (0,-1)};

        \filldraw[white] (0,-1) circle(2pt);

       \filldraw[white] (-1, 0) circle(2pt);
       \filldraw[white] (1, 0) circle(2pt);

       \filldraw[white] (-1, 1.5) circle(2pt);
       \filldraw[white] (1, 1.5) circle(2pt);

       \filldraw[white] (-1, 3) circle(2pt);
       \filldraw[white] (1, 3) circle(2pt);

       \filldraw[white] (1, 4) circle(2pt);
       \filldraw[white] (0, 4) circle(2pt);

\end{tikzpicture}
}

\centering 3

\end{minipage}
\begin{minipage}{.3\textwidth}              

\resizebox{4.5cm}{!}{
\begin{tikzpicture}
\coordinate (center) at (0,0);
   \def\radius{2.2cm}
   \foreach \x in {0, 60,...,360} {
               }


       \coordinate (v1) at (-1, 0);\filldraw[black] (v1) circle(4pt);\node[left] at (-1.2, 0)   {\LARGE${v_3}$};
       \coordinate (v2) at (1,0);\filldraw[black] (v2) circle(4pt);\node[right] at (1.2,0)    {\LARGE${v_{8}}$};
       \coordinate (v3) at (-1, 1.5);\filldraw[black] (v3) circle(4pt);\node[left] at (v3)   {\LARGE${v_5}$};
       \coordinate (v4) at (1, 1.5);\filldraw[black] (v4) circle(4pt);\node[right] at (v4)   {\LARGE${v_1}$};
       \coordinate (v5) at (-1, 3);\filldraw[black] (v5) circle(4pt);\node[left] at (v5)   {\LARGE${v_2}$};
       \coordinate (v6) at (1, 3);\filldraw[black] (v6) circle(4pt);\node[right] at (v6)   {\LARGE${v_7}$};
       \coordinate (v7) at (-1, 4.5);\filldraw[black] (v7) circle(4pt);\node[left] at (-1.2, 4.5)   {\LARGE${v_6}$};
       \coordinate (v8) at (1, 4.5);\filldraw[black] (v8) circle(4pt);\node[right] at (1.2, 4.5)   {\LARGE${v_{4}}$};

     \draw [line width=1.5,black](v3) -- (v5) -- (v7) -- (v8) -- (v6) -- (v4) -- (v1) -- (v2) -- (v3);
     \draw [line width=1.5,black](v3) -- (v4);
     \draw [line width=1.5,black](v5) -- (v6);

\draw [line width=1.5, black] plot [smooth,  tension=.9] coordinates {(1,0) (2,1) (2,3.5) (1, 4.5)};
\draw [line width=1.5, black] plot [smooth,  tension=.9] coordinates {(-1,0) (-2,1) (-2,3.5) (-1, 4.5)};

         \filldraw[white] (-1, 0) circle(2pt);
       \filldraw[white] (1, 0) circle(2pt);
       \filldraw[white] (-1, 1.5) circle(2pt);
       \filldraw[white] (1, 1.5) circle(2pt);
       \filldraw[white] (-1, 3) circle(2pt);
       \filldraw[white] (1, 3) circle(2pt);
       \filldraw[white] (-1, 4.5) circle(2pt);
       \filldraw[white] (1, 4.5) circle(2pt);

       \filldraw[white] (0, 5.7) circle(2pt);
       \filldraw[white] (0, -1.2) circle(2pt);

\end{tikzpicture}
}

\centering 4

\end{minipage}                               
\begin{minipage}{.3\textwidth}              

\resizebox{4.cm}{!}{
\begin{tikzpicture}
\coordinate (center) at (0,0);
   \def\radius{2.2cm}
   \foreach \x in {0, 60,...,360} {
               }


       \coordinate (v1) at (-1, 0);\filldraw[black] (v1) circle(4pt);\node[left] at (v1)   {\LARGE${v_4}$};
       \coordinate (v2) at (1,0);\filldraw[black] (v2) circle(4pt);\node[right] at (v2)    {\LARGE${v_{2}}$};
       \coordinate (v3) at (-1, 1);\filldraw[black] (v3) circle(4pt);\node[left] at (v3)   {\LARGE${v_7}$};
       \coordinate (v4) at (1, 1);\filldraw[black] (v4) circle(4pt);\node[right] at (v4)   {\LARGE${v_8}$};
       \coordinate (v5) at (-1, 2);\filldraw[black] (v5) circle(4pt);\node[left] at (v5)   {\LARGE${v_5}$};
       \coordinate (v6) at (1, 2);\filldraw[black] (v6) circle(4pt);\node[right] at (v6)   {\LARGE${v_6}$};
       \coordinate (v7) at (-1, 3);\filldraw[black] (v7) circle(4pt);\node[left] at (v7)   {\LARGE${v_1}$};
       \coordinate (v8) at (1, 3);\filldraw[black] (v8) circle(4pt);\node[right] at (v8)   {\LARGE${v_{3}}$};

     \draw [line width=1.5,black](v8) -- (v6) -- (v4) -- (v2) -- (v1) -- (v3) -- (v5) -- (v7);
     \draw [line width=1.5,black](v1) -- (v6);
     \draw [line width=1.5,black](v2) -- (v3);
     \draw [line width=1.5,black](v4) -- (v5);

     \draw [line width=1.5,black](v7) to[bend right] (v8);
     \draw [line width=1.5,black](v8) to[bend right] (v7);

       \filldraw[white] (-1, 0) circle(2pt);
       \filldraw[white] (1, 0) circle(2pt);
       \filldraw[white] (-1, 1) circle(2pt);
       \filldraw[white] (1, 1) circle(2pt);
       \filldraw[white] (-1, 2) circle(2pt);
       \filldraw[white] (1, 2) circle(2pt);
       \filldraw[white] (-1, 3) circle(2pt);
       \filldraw[white] (1, 3) circle(2pt);

       \filldraw[white] (0, 5) circle(2pt);
       \filldraw[white] (0, -1) circle(2pt);

\end{tikzpicture}
}

\centering  5

\end{minipage}
\begin{minipage}{.3\textwidth}              
\resizebox{4cm}{!}{
\begin{tikzpicture}
\coordinate (center) at (0,0);
   \def\radius{2.2cm}
   \foreach \x in {0, 60,...,360} {
               }


       \coordinate (v1) at (-1, 0);\filldraw[black] (v1) circle(4pt);\node[left] at (v1)   {\LARGE${v_4}$};
       \coordinate (v2) at (1,0);\filldraw[black] (v2) circle(4pt);\node[right] at (v2)    {\LARGE${v_{7}}$};
       \coordinate (v3) at (-1, 1);\filldraw[black] (v3) circle(4pt);\node[left] at (v3)   {\LARGE${v_8}$};
       \coordinate (v4) at (1, 1);\filldraw[black] (v4) circle(4pt);\node[right] at (v4)   {\LARGE${v_2}$};
       \coordinate (v5) at (-1, 2);\filldraw[black] (v5) circle(4pt);\node[left] at (v5)   {\LARGE${v_5}$};
       \coordinate (v6) at (1, 2);\filldraw[black] (v6) circle(4pt);\node[right] at (v6)   {\LARGE${v_6}$};
       \coordinate (v7) at (-1, 3);\filldraw[black] (v7) circle(4pt);\node[left] at (v7)   {\LARGE${v_1}$};
       \coordinate (v8) at (1, 3);\filldraw[black] (v8) circle(4pt);\node[right] at (v8)   {\LARGE${v_{3}}$};

     \draw [line width=1.5,black](v8) -- (v6) -- (v4) -- (v2) -- (v1) -- (v3) -- (v5) -- (v7);
     \draw [line width=1.5,black](v1) -- (v6);
     \draw [line width=1.5,black](v4) -- (v3);
     \draw [line width=1.5,black](v2) -- (v5);

     \draw [line width=1.5,black](v7) to[bend right] (v8);
     \draw [line width=1.5,black](v8) to[bend right] (v7);

       \filldraw[white] (-1, 0) circle(2pt);
       \filldraw[white] (1, 0) circle(2pt);
       \filldraw[white] (-1, 1) circle(2pt);
       \filldraw[white] (1, 1) circle(2pt);
       \filldraw[white] (-1, 2) circle(2pt);
       \filldraw[white] (1, 2) circle(2pt);
       \filldraw[white] (-1, 3) circle(2pt);
       \filldraw[white] (1, 3) circle(2pt);

       \filldraw[white] (0, 5) circle(2pt);
       \filldraw[white] (0, -1) circle(2pt);

\end{tikzpicture}
}

\centering 6

\end{minipage}
\begin{minipage}{.3\textwidth}              
\resizebox{4cm}{!}{
\begin{tikzpicture}
\coordinate (center) at (0,0);
   \def\radius{2.2cm}
   \foreach \x in {0, 60,...,360} {
               }


       \coordinate (v1) at (-1, 0);\filldraw[black] (v1) circle(4pt);\node[left] at (v1)   {\LARGE${v_8}$};
       \coordinate (v2) at (1,0);\filldraw[black] (v2) circle(4pt);\node[right] at (v2)    {\LARGE${v_{7}}$};
       \coordinate (v3) at (-1, 1);\filldraw[black] (v3) circle(4pt);\node[left] at (v3)   {\LARGE${v_2}$};
       \coordinate (v4) at (1, 1);\filldraw[black] (v4) circle(4pt);\node[right] at (v4)   {\LARGE${v_4}$};
       \coordinate (v5) at (-1, 2);\filldraw[black] (v5) circle(4pt);\node[left] at (v5)   {\LARGE${v_5}$};
       \coordinate (v6) at (1, 2);\filldraw[black] (v6) circle(4pt);\node[right] at (v6)   {\LARGE${v_6}$};
       \coordinate (v7) at (-1, 3);\filldraw[black] (v7) circle(4pt);\node[left] at (v7)   {\LARGE${v_1}$};
       \coordinate (v8) at (1, 3);\filldraw[black] (v8) circle(4pt);\node[right] at (v8)   {\LARGE${v_{3}}$};

     \draw [line width=1.5,black](v8) -- (v6) -- (v4) -- (v2) -- (v3) -- (v5) -- (v7);
     \draw [line width=1.5,black](v1) -- (v3);
     \draw [line width=1.5,black](v1) -- (v4);
     \draw [line width=1.5,black](v1) -- (v6);
     \draw [line width=1.5,black](v2) -- (v5);

     \draw [line width=1.5,black](v7) to[bend right] (v8);
     \draw [line width=1.5,black](v8) to[bend right] (v7);

       \filldraw[white] (-1, 0) circle(2pt);
       \filldraw[white] (1, 0) circle(2pt);
       \filldraw[white] (-1, 1) circle(2pt);
       \filldraw[white] (1, 1) circle(2pt);
       \filldraw[white] (-1, 2) circle(2pt);
       \filldraw[white] (1, 2) circle(2pt);
       \filldraw[white] (-1, 3) circle(2pt);
       \filldraw[white] (1, 3) circle(2pt);

       \filldraw[white] (0, 5) circle(2pt);
       \filldraw[white] (0, -1) circle(2pt);

\end{tikzpicture}
}

\centering 7

\end{minipage}
\begin{minipage}{.3\textwidth}              

\resizebox{4cm}{!}{
\begin{tikzpicture}
\coordinate (center) at (0,0);
   \def\radius{2.2cm}
   \foreach \x in {0, 60,...,360} {
               }


       \coordinate (v1) at (-1, 0);\filldraw[black] (v1) circle(4pt);\node[left] at (v1)   {\LARGE${v_8}$};
       \coordinate (v2) at (1,0);\filldraw[black] (v2) circle(4pt);\node[right] at (v2)    {\LARGE${v_{6}}$};
       \coordinate (v3) at (-1, 1);\filldraw[black] (v3) circle(4pt);\node[left] at (v3)   {\LARGE${v_2}$};
       \coordinate (v4) at (1, 1);\filldraw[black] (v4) circle(4pt);\node[right] at (v4)   {\LARGE${v_4}$};
       \coordinate (v5) at (-1, 2);\filldraw[black] (v5) circle(4pt);\node[left] at (v5)   {\LARGE${v_5}$};
       \coordinate (v6) at (1, 2);\filldraw[black] (v6) circle(4pt);\node[right] at (v6)   {\LARGE${v_7}$};
       \coordinate (v7) at (-1, 3);\filldraw[black] (v7) circle(4pt);\node[left] at (v7)   {\LARGE${v_1}$};
       \coordinate (v8) at (1, 3);\filldraw[black] (v8) circle(4pt);\node[right] at (v8)   {\LARGE${v_{3}}$};

     \draw [line width=1.5,black](v8) -- (v6) -- (v4) -- (v2) --(v1) -- (v3) -- (v5) -- (v7);
     \draw [line width=1.5,black](v2) -- (v3);
     \draw [line width=1.5,black](v1) -- (v4);
     \draw [line width=1.5,black](v6) -- (v5);

     \draw [line width=1.5,black](v7) to[bend right] (v8);
     \draw [line width=1.5,black](v8) to[bend right] (v7);

       \filldraw[white] (-1, 0) circle(2pt);
       \filldraw[white] (1, 0) circle(2pt);
       \filldraw[white] (-1, 1) circle(2pt);
       \filldraw[white] (1, 1) circle(2pt);
       \filldraw[white] (-1, 2) circle(2pt);
       \filldraw[white] (1, 2) circle(2pt);
       \filldraw[white] (-1, 3) circle(2pt);
       \filldraw[white] (1, 3) circle(2pt);

       \filldraw[white] (0, 5) circle(2pt);
       \filldraw[white] (0, -1) circle(2pt);

\end{tikzpicture}
}

\centering 8

\end{minipage}
\begin{minipage}{.3\textwidth}              

\resizebox{4cm}{!}{
\begin{tikzpicture}
\coordinate (center) at (0,0);
   \def\radius{2.2cm}
   \foreach \x in {0, 60,...,360} {
               }


       \coordinate (v1) at (-1, 0);\filldraw[black] (v1) circle(4pt);\node[left] at (v1)   {\LARGE${v_8}$};
       \coordinate (v2) at (1,0);\filldraw[black] (v2) circle(4pt);\node[right] at (v2)    {\LARGE${v_{6}}$};
       \coordinate (v3) at (-1, 1);\filldraw[black] (v3) circle(4pt);\node[left] at (v3)   {\LARGE${v_2}$};
       \coordinate (v4) at (1, 1);\filldraw[black] (v4) circle(4pt);\node[right] at (v4)   {\LARGE${v_4}$};
       \coordinate (v5) at (-1, 2);\filldraw[black] (v5) circle(4pt);\node[left] at (v5)   {\LARGE${v_5}$};
       \coordinate (v6) at (1, 2);\filldraw[black] (v6) circle(4pt);\node[right] at (v6)   {\LARGE${v_7}$};
       \coordinate (v7) at (-1, 3);\filldraw[black] (v7) circle(4pt);\node[left] at (v7)   {\LARGE${v_1}$};
       \coordinate (v8) at (1, 3);\filldraw[black] (v8) circle(4pt);\node[right] at (v8)   {\LARGE${v_{3}}$};

     \draw [line width=1.5,black](v6) -- (v4) -- (v2) --(v1) -- (v3) -- (v6);
     \draw [line width=1.5,black](v5) -- (v6);
     \draw [line width=1.5,black](v5) -- (v7);
     \draw [line width=1.5,black](v5) -- (v8);

     \draw [line width=1.5,black](v1) -- (v4);
     \draw [line width=1.5,black](v2) -- (v3);

     \draw [line width=1.5,black](v7) to[bend right] (v8);
     \draw [line width=1.5,black](v8) to[bend right] (v7);

       \filldraw[white] (-1, 0) circle(2pt);
       \filldraw[white] (1, 0) circle(2pt);
       \filldraw[white] (-1, 1) circle(2pt);
       \filldraw[white] (1, 1) circle(2pt);
       \filldraw[white] (-1, 2) circle(2pt);
       \filldraw[white] (1, 2) circle(2pt);
       \filldraw[white] (-1, 3) circle(2pt);
       \filldraw[white] (1, 3) circle(2pt);

       \filldraw[white] (0, 5) circle(2pt);
       \filldraw[white] (0, -1) circle(2pt);

\end{tikzpicture}
}

\centering 9

\end{minipage}
\caption{Cubic graphs on 8 vertices}    \label{cubic graphs on 8 vertices}

\end{center}
\end{figure}
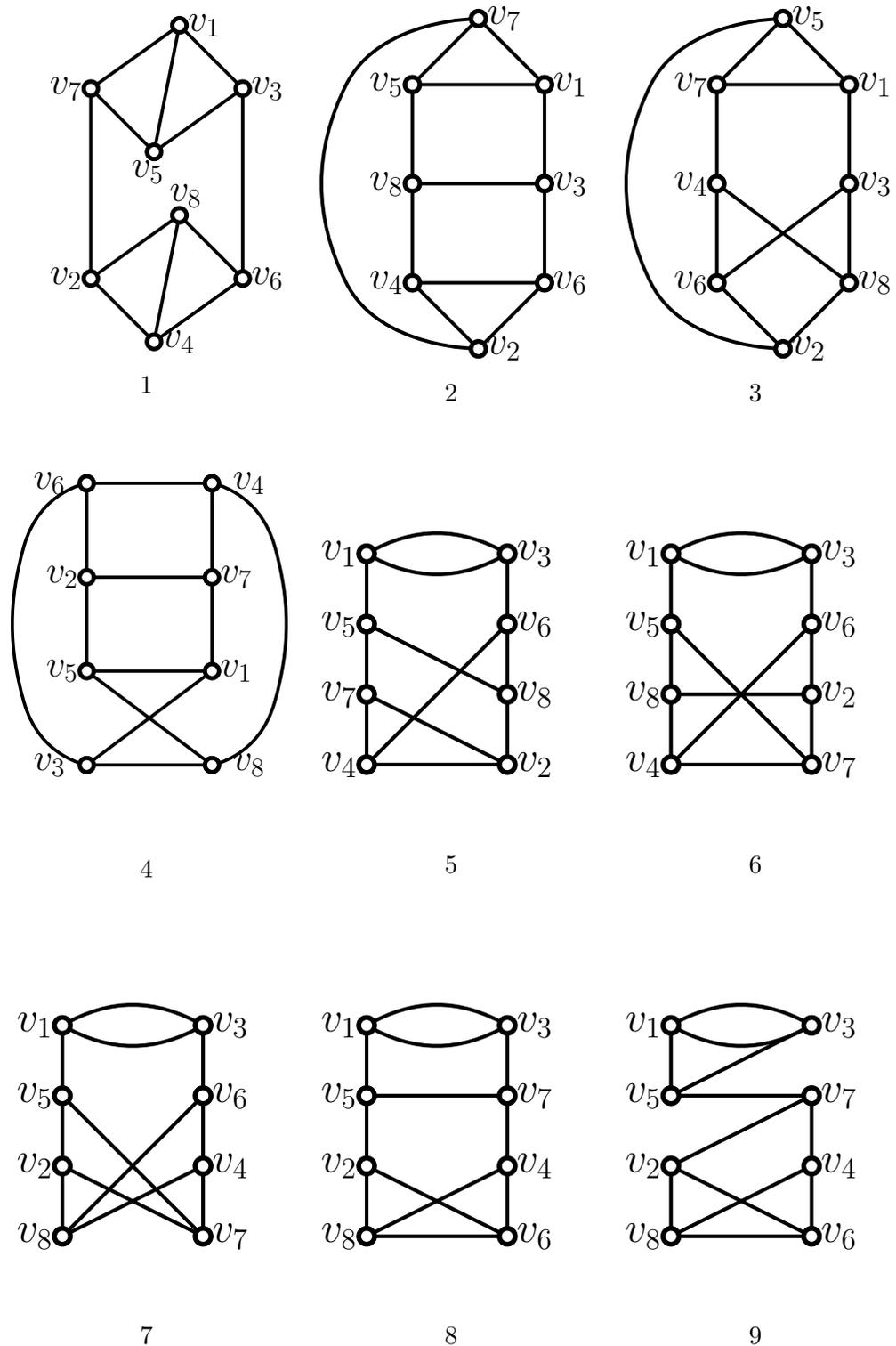

 \vspace{5mm}  For  $n=8\/$, we have checked that each cubic graph on $8\/$ vertices, except $Q_3\/$ cube admits a labeling on its vertices that satisfies the conditions (i) and (ii) (see Figure \ref{cubic graphs on 8 vertices}).

 \vspace{5mm}   Let $G\/$ be a $3\/$-regular  graph on $n\/$ vertices where $n \geq 10\/$.

\vspace{3mm}{\bf{Case (1)}    } $G$   has no $2\/$-cycle.

 \vspace{5mm}   Let $e=xy$ be chosen such that $G_e\/$ has at most one $2\/$-cycle.

 \vspace{5mm}   (1.1) $G_e$   has no $2\/$-cycle.

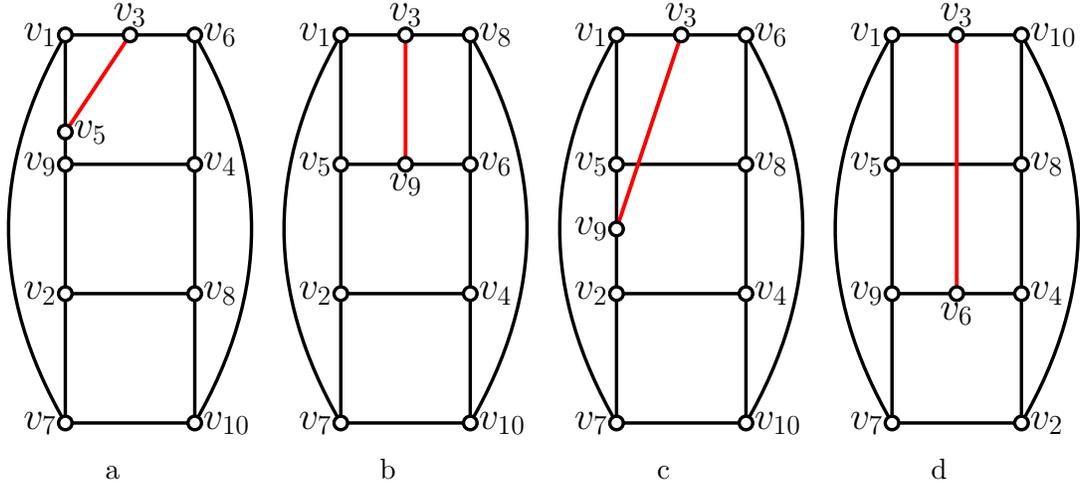
\begin{figure}[htp]
\begin{center}
\begin{minipage}{.24\textwidth}

\resizebox{4.cm}{!}{
\begin{tikzpicture}
\coordinate (center) at (0,0);
   \def\radius{2.2cm}
   \foreach \x in {0, 60,...,360} {
               }


       \coordinate (v1) at (-1, 0);\filldraw[black] (v1) circle(3.5pt);\node[left] at (v1)   {\LARGE${v_7}$};
       \coordinate (v2) at (1,0);\filldraw[black] (v2) circle(3.5pt);\node[right] at (v2)    {\LARGE${v_{10}}$};
       \coordinate (v3) at (-1, 2);\filldraw[black] (v3) circle(3.5pt);\node[left] at (v3)   {\LARGE${v_2}$};
       \coordinate (v4) at (1, 2);\filldraw[black] (v4) circle(3.5pt);\node[right] at (v4)   {\LARGE${v_8}$};
       \coordinate (v5) at (-1, 4);\filldraw[black] (v5) circle(3.5pt);\node[left] at (v5)   {\LARGE${v_9}$};
       \coordinate (v6) at (1, 4);\filldraw[black] (v6) circle(3.5pt);\node[right] at (v6)   {\LARGE${v_4}$};
       \coordinate (v7) at (-1, 6);\filldraw[black] (v7) circle(3.5pt);\node[left] at (v7)   {\LARGE${v_1}$};
       \coordinate (v8) at (1, 6);\filldraw[black] (v8) circle(3.5pt);\node[right] at (v8)   {\LARGE${v_6}$};

     \draw [line width=1.5,black](v1) -- (v3) -- (v5) -- (v7) -- (v8) -- (v6) -- (v4) -- (v2) -- (v1);
     \draw [line width=1.5,black](v3) -- (v4);
     \draw [line width=1.5,black](v5) -- (v6);
     \draw [line width=1.5,black](v2) to[bend right] (v8);
     \draw [line width=1.5,black](v7) to[bend right] (v1);

       \filldraw[white] (-1, 0) circle(2pt);
       \filldraw[white] (1, 0) circle(2pt);
       \filldraw[white] (-1, 2) circle(2pt);
       \filldraw[white] (1, 2) circle(2pt);
       \filldraw[white] (-1, 4) circle(2pt);
       \filldraw[white] (1, 4) circle(2pt);
       \filldraw[white] (-1, 6) circle(2pt);
       \filldraw[white] (1, 6) circle(2pt);

      \coordinate (v9) at (0, 6);
      \coordinate (v10) at (-1, 4.5);
     \draw [line width=1.7,red](v9) -- (v10);

      \filldraw[black] (v9) circle(3.5pt);\node[above] at (v9)     {\LARGE${v_3}$};
     \filldraw[black] (v10) circle(3.5pt);\node[right] at (v10) {\LARGE${v_{5}}$};

       \filldraw[white] (0, 6) circle(2pt);
       \filldraw[white] (-1, 4.5) circle(2pt);

\end{tikzpicture}
}

\centering  a

\end{minipage}
\begin{minipage}{.24\textwidth}
\resizebox{4cm}{!}{
\begin{tikzpicture}
\coordinate (center) at (0,0);
   \def\radius{2.2cm}
   \foreach \x in {0, 60,...,360} {
               }


       \coordinate (v1) at (-1, 0);\filldraw[black] (v1) circle(3.5pt);\node[left] at (v1)   {\LARGE${v_7}$};
       \coordinate (v2) at (1,0);\filldraw[black] (v2) circle(3.5pt);\node[right] at (v2)    {\LARGE${v_{10}}$};
       \coordinate (v3) at (-1, 2);\filldraw[black] (v3) circle(3.5pt);\node[left] at (v3)   {\LARGE${v_2}$};
       \coordinate (v4) at (1, 2);\filldraw[black] (v4) circle(3.5pt);\node[right] at (v4)   {\LARGE${v_4}$};
       \coordinate (v5) at (-1, 4);\filldraw[black] (v5) circle(3.5pt);\node[left] at (v5)   {\LARGE${v_5}$};
       \coordinate (v6) at (1, 4);\filldraw[black] (v6) circle(3.5pt);\node[right] at (v6)   {\LARGE${v_6}$};
       \coordinate (v7) at (-1, 6);\filldraw[black] (v7) circle(3.5pt);\node[left] at (v7)   {\LARGE${v_1}$};
       \coordinate (v8) at (1, 6);\filldraw[black] (v8) circle(3.5pt);\node[right] at (v8)   {\LARGE${v_8}$};

     \draw [line width=1.5,black](v1) -- (v3) -- (v5) -- (v7) -- (v8) -- (v6) -- (v4) -- (v2) -- (v1);
     \draw [line width=1.5,black](v3) -- (v4);
     \draw [line width=1.5,black](v5) -- (v6);
     \draw [line width=1.5,black](v2) to[bend right] (v8);
     \draw [line width=1.5,black](v7) to[bend right] (v1);

       \filldraw[white] (-1, 0) circle(2pt);
       \filldraw[white] (1, 0) circle(2pt);
       \filldraw[white] (-1, 2) circle(2pt);
       \filldraw[white] (1, 2) circle(2pt);
       \filldraw[white] (-1, 4) circle(2pt);
       \filldraw[white] (1, 4) circle(2pt);
       \filldraw[white] (-1, 6) circle(2pt);
       \filldraw[white] (1, 6) circle(2pt);

      \coordinate (v9) at (0, 6);
      \coordinate (v10) at (0, 4);
     \draw [line width=1.7,red](v9) -- (v10);

      \filldraw[black] (v9) circle(3.5pt);\node[above] at (v9)     {\LARGE${v_3}$};
     \filldraw[black] (v10) circle(3.5pt);\node[below] at (v10) {\LARGE${v_9}$};

       \filldraw[white] (0, 6) circle(2pt);
       \filldraw[white] (0, 4) circle(2pt);

\end{tikzpicture}
}

\centering  b

\end{minipage}
\begin{minipage}{.24\textwidth}
\resizebox{4cm}{!}{
\begin{tikzpicture}
\coordinate (center) at (0,0);
   \def\radius{2.2cm}
   \foreach \x in {0, 60,...,360} {
               }


       \coordinate (v1) at (-1, 0);\filldraw[black] (v1) circle(3.5pt);\node[left] at (v1)   {\LARGE${v_7}$};
       \coordinate (v2) at (1,0);\filldraw[black] (v2) circle(3.5pt);\node[right] at (v2)    {\LARGE${v_{10}}$};
       \coordinate (v3) at (-1, 2);\filldraw[black] (v3) circle(3.5pt);\node[left] at (v3)   {\LARGE${v_2}$};
       \coordinate (v4) at (1, 2);\filldraw[black] (v4) circle(3.5pt);\node[right] at (v4)   {\LARGE${v_4}$};
       \coordinate (v5) at (-1, 4);\filldraw[black] (v5) circle(3.5pt);\node[left] at (v5)   {\LARGE${v_5}$};
       \coordinate (v6) at (1, 4);\filldraw[black] (v6) circle(3.5pt);\node[right] at (v6)   {\LARGE${v_8}$};
       \coordinate (v7) at (-1, 6);\filldraw[black] (v7) circle(3.5pt);\node[left] at (v7)   {\LARGE${v_1}$};
       \coordinate (v8) at (1, 6);\filldraw[black] (v8) circle(3.5pt);\node[right] at (v8)   {\LARGE${v_6}$};

     \draw [line width=1.5,black](v1) -- (v3) -- (v5) -- (v7) -- (v8) -- (v6) -- (v4) -- (v2) -- (v1);
     \draw [line width=1.5,black](v3) -- (v4);
     \draw [line width=1.5,black](v5) -- (v6);
     \draw [line width=1.5,black](v2) to[bend right] (v8);
     \draw [line width=1.5,black](v7) to[bend right] (v1);

       \filldraw[white] (-1, 0) circle(2pt);
       \filldraw[white] (1, 0) circle(2pt);
       \filldraw[white] (-1, 2) circle(2pt);
       \filldraw[white] (1, 2) circle(2pt);
       \filldraw[white] (-1, 4) circle(2pt);
       \filldraw[white] (1, 4) circle(2pt);
       \filldraw[white] (-1, 6) circle(2pt);
       \filldraw[white] (1, 6) circle(2pt);

      \coordinate (v9) at (0, 6);
      \coordinate (v10) at (-1, 3);
     \draw [line width=1.7,red](v9) -- (v10);

      \filldraw[black] (v9) circle(3.5pt);\node[above] at (v9)     {\LARGE${v_3}$};
     \filldraw[black] (v10) circle(3.5pt);\node[left] at (v10) {\LARGE${v_{9}}$};

       \filldraw[white] (0, 6) circle(2pt);
       \filldraw[white] (-1, 3) circle(2pt);

\end{tikzpicture}
}

\centering c

\end{minipage}
\begin{minipage}{.24\textwidth}

\resizebox{4cm}{!}{
\begin{tikzpicture}
\coordinate (center) at (0,0);
   \def\radius{2.2cm}
   \foreach \x in {0, 60,...,360} {
               }


       \coordinate (v1) at (-1, 0);\filldraw[black] (v1) circle(3.5pt);\node[left] at (v1)   {\LARGE${v_7}$};
       \coordinate (v2) at (1,0);\filldraw[black] (v2) circle(3.5pt);\node[right] at (v2)    {\LARGE${v_{2}}$};
       \coordinate (v3) at (-1, 2);\filldraw[black] (v3) circle(3.5pt);\node[left] at (v3)   {\LARGE${v_9}$};
       \coordinate (v4) at (1, 2);\filldraw[black] (v4) circle(3.5pt);\node[right] at (v4)   {\LARGE${v_4}$};
       \coordinate (v5) at (-1, 4);\filldraw[black] (v5) circle(3.5pt);\node[left] at (v5)   {\LARGE${v_5}$};
       \coordinate (v6) at (1, 4);\filldraw[black] (v6) circle(3.5pt);\node[right] at (v6)   {\LARGE${v_8}$};
       \coordinate (v7) at (-1, 6);\filldraw[black] (v7) circle(3.5pt);\node[left] at (v7)   {\LARGE${v_1}$};
       \coordinate (v8) at (1, 6);\filldraw[black] (v8) circle(3.5pt);\node[right] at (v8)   {\LARGE${v_{10}}$};

     \draw [line width=1.5,black](v1) -- (v3) -- (v5) -- (v7) -- (v8) -- (v6) -- (v4) -- (v2) -- (v1);
     \draw [line width=1.5,black](v3) -- (v4);
     \draw [line width=1.5,black](v5) -- (v6);
     \draw [line width=1.5,black](v2) to[bend right] (v8);
     \draw [line width=1.5,black](v7) to[bend right] (v1);

       \filldraw[white] (-1, 0) circle(2pt);
       \filldraw[white] (1, 0) circle(2pt);
       \filldraw[white] (-1, 2) circle(2pt);
       \filldraw[white] (1, 2) circle(2pt);
       \filldraw[white] (-1, 4) circle(2pt);
       \filldraw[white] (1, 4) circle(2pt);
       \filldraw[white] (-1, 6) circle(2pt);
       \filldraw[white] (1, 6) circle(2pt);

      \coordinate (v9) at (0, 6);
      \coordinate (v10) at (0, 2);
     \draw [line width=1.7,red](v9) -- (v10);

      \filldraw[black] (v9) circle(3.5pt);\node[above] at (v9)     {\LARGE${v_3}$};
     \filldraw[black] (v10) circle(3.5pt);\node[below] at (v10) {\LARGE${v_{6}}$};

       \filldraw[white] (0, 6) circle(2pt);
       \filldraw[white] (0, 2) circle(2pt);

\end{tikzpicture}
}

\centering d

\end{minipage}
\caption{Cubic on 10 vertices}    \label{cubic on 10 vertices}
\end{center}
\end{figure}

 \vspace{5mm}   If  $G_e \/$ is the cube $Q_3\/$, then $G\/$ is any one of the four  $3\/$-regular graphs depicted in Figure \ref{cubic on 10 vertices}. Each of these graphs has a labeling that satisfies condition (i) of the lemma.



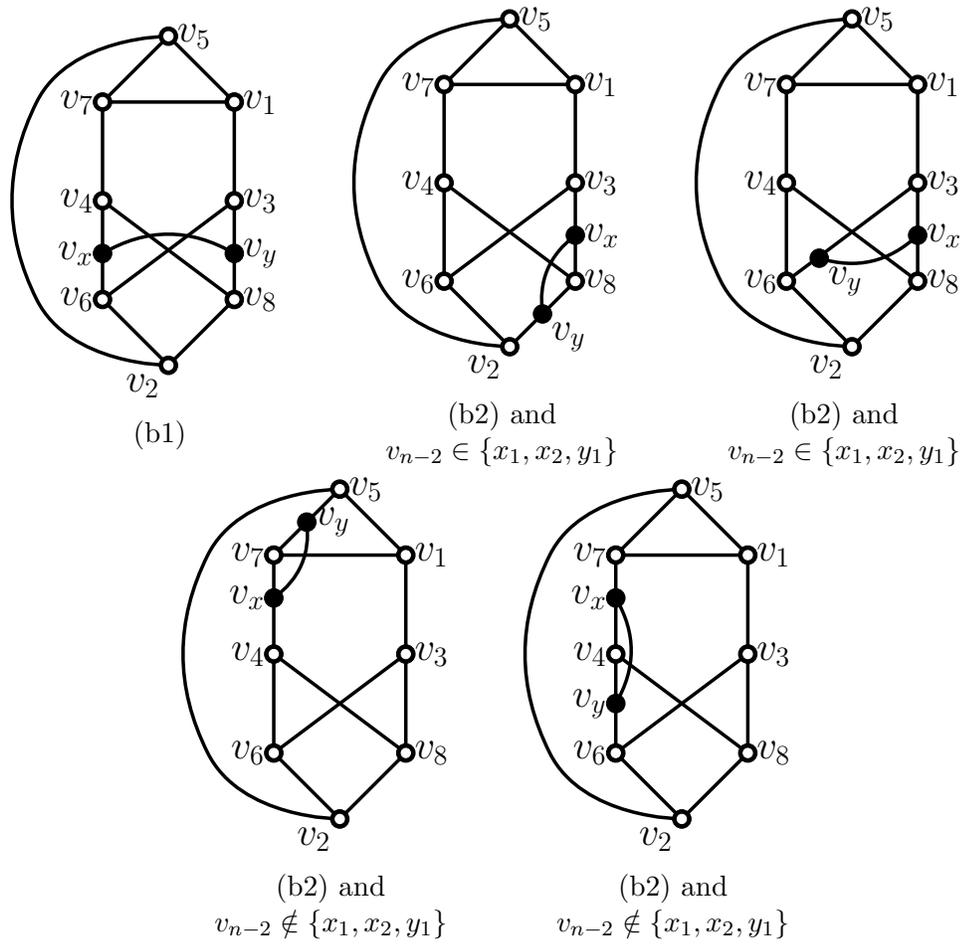
\begin{figure}[htp]
\begin{center}
\begin{minipage}{.3\textwidth}              
\resizebox{4.cm}{!}{
\begin{tikzpicture}
\coordinate (center) at (0,0);
   \def\radius{2.2cm}
   \foreach \x in {0, 60,...,360} {
               }


       \coordinate (v1) at (-1, 0);\filldraw[black] (v1) circle(4pt);    \node[left] at (v1)    {\LARGE${v_6}$};
       \coordinate (v2) at (1,0);\filldraw[black] (v2) circle(4pt);      \node[right] at (v2)   {\LARGE${v_{8}}$};
       \coordinate (v3) at (-1, 1.5);\filldraw[black] (v3) circle(4pt);    \node[left] at (v3)    {\LARGE${v_4}$};
       \coordinate (v4) at (1, 1.5);\filldraw[black] (v4) circle(4pt);     \node[right] at (v4)   {\LARGE${v_3}$};
       \coordinate (v5) at (-1, 3);\filldraw[black] (v5) circle(4pt);    \node[left] at (v5)    {\LARGE${v_7}$};
       \coordinate (v6) at (1, 3);\filldraw[black] (v6) circle(4pt);     \node[right] at (v6)   {\LARGE${v_1}$};

      \coordinate (v9) at (0, 4);  \filldraw[black] (v9) circle(4pt);     \node[right] at (v9)     {\LARGE${v_5}$};
      \coordinate (v10) at (0,-1);    \filldraw[black] (v10) circle(4pt);     \node[below left] at (v10) {\LARGE${v_{2}}$};

      \coordinate (v11) at (1, .7);    \filldraw[black] (v11) circle(4pt);     \node[right] at (v11) {\LARGE${v_{y}}$};
      \coordinate (v12) at (-1, .7);    \filldraw[black] (v12) circle(4pt);     \node[left] at (v12) {\LARGE${v_{x}}$};

     \draw [line width=1.5,black](v10) -- (v1) -- (v3) -- (v5) -- (v9) -- (v6) -- (v4) -- (v2) -- (v10);
     \draw [line width=1.5,black](v1) -- (v4);
     \draw [line width=1.5,black](v5) -- (v6);
     \draw [line width=1.5,black](v3) -- (v2);

\draw [line width=1.5, black] plot [smooth,  tension=.9] coordinates {(0,4) (-2,3) (-2,0) (0,-1)};

     \draw [line width=1.5,black](v11) to[bend right] (v12);

        \filldraw[white] (0,-1) circle(2pt);

       \filldraw[white] (-1, 0) circle(2pt);
       \filldraw[white] (1, 0) circle(2pt);

       \filldraw[white] (-1, 1.5) circle(2pt);
       \filldraw[white] (1, 1.5) circle(2pt);

       \filldraw[white] (-1, 3) circle(2pt);
       \filldraw[white] (1, 3) circle(2pt);

       \filldraw[white] (1, 4) circle(2pt);
       \filldraw[white] (0, 4) circle(2pt);

\end{tikzpicture}
}

\centering  (b1)

\end{minipage}
\begin{minipage}{.3\textwidth}              
\resizebox{4.cm}{!}{
\begin{tikzpicture}
\coordinate (center) at (0,0);
   \def\radius{2.2cm}
   \foreach \x in {0, 60,...,360} {
               }


       \coordinate (v1) at (-1, 0);\filldraw[black] (v1) circle(4pt);    \node[left] at (v1)    {\LARGE${v_6}$};
       \coordinate (v2) at (1,0);\filldraw[black] (v2) circle(4pt);      \node[right] at (v2)   {\LARGE${v_{8}}$};
       \coordinate (v3) at (-1, 1.5);\filldraw[black] (v3) circle(4pt);    \node[left] at (v3)    {\LARGE${v_4}$};
       \coordinate (v4) at (1, 1.5);\filldraw[black] (v4) circle(4pt);     \node[right] at (v4)   {\LARGE${v_3}$};
       \coordinate (v5) at (-1, 3);\filldraw[black] (v5) circle(4pt);    \node[left] at (v5)    {\LARGE${v_7}$};
       \coordinate (v6) at (1, 3);\filldraw[black] (v6) circle(4pt);     \node[right] at (v6)   {\LARGE${v_1}$};

      \coordinate (v9) at (0, 4);  \filldraw[black] (v9) circle(4pt);     \node[right] at (v9)     {\LARGE${v_5}$};
      \coordinate (v10) at (0,-1);    \filldraw[black] (v10) circle(4pt);     \node[below left] at (v10) {\LARGE${v_{2}}$};

      \coordinate (v11) at (1, .7);    \filldraw[black] (v11) circle(4pt);     \node[right] at (v11) {\LARGE${v_{x}}$};
      \coordinate (v12) at (0.5,-0.5);    \filldraw[black] (v12) circle(4pt);     \node[below right] at (v12) {\LARGE${v_{y}}$};

     \draw [line width=1.5,black](v10) -- (v1) -- (v3) -- (v5) -- (v9) -- (v6) -- (v4) -- (v2) -- (v10);
     \draw [line width=1.5,black](v1) -- (v4);
     \draw [line width=1.5,black](v5) -- (v6);
     \draw [line width=1.5,black](v3) -- (v2);

\draw [line width=1.5, black] plot [smooth,  tension=.9] coordinates {(0,4) (-2,3) (-2,0) (0,-1)};

     \draw [line width=1.5,black](v11) to[bend right] (v12);

        \filldraw[white] (0,-1) circle(2pt);

       \filldraw[white] (-1, 0) circle(2pt);
       \filldraw[white] (1, 0) circle(2pt);

       \filldraw[white] (-1, 1.5) circle(2pt);
       \filldraw[white] (1, 1.5) circle(2pt);

       \filldraw[white] (-1, 3) circle(2pt);
       \filldraw[white] (1, 3) circle(2pt);

       \filldraw[white] (1, 4) circle(2pt);
       \filldraw[white] (0, 4) circle(2pt);

\end{tikzpicture}
}

\centering  (b2)  and  $v_{n-2}\in \{x_1, x_2, y_1\}$

\end{minipage}
\begin{minipage}{.3\textwidth}              
\resizebox{4.cm}{!}{
\begin{tikzpicture}
\coordinate (center) at (0,0);
   \def\radius{2.2cm}
   \foreach \x in {0, 60,...,360} {
               }


       \coordinate (v1) at (-1, 0);\filldraw[black] (v1) circle(4pt);    \node[left] at (v1)    {\LARGE${v_6}$};
       \coordinate (v2) at (1,0);\filldraw[black] (v2) circle(4pt);      \node[right] at (v2)   {\LARGE${v_{8}}$};
       \coordinate (v3) at (-1, 1.5);\filldraw[black] (v3) circle(4pt);    \node[left] at (v3)    {\LARGE${v_4}$};
       \coordinate (v4) at (1, 1.5);\filldraw[black] (v4) circle(4pt);     \node[right] at (v4)   {\LARGE${v_3}$};
       \coordinate (v5) at (-1, 3);\filldraw[black] (v5) circle(4pt);    \node[left] at (v5)    {\LARGE${v_7}$};
       \coordinate (v6) at (1, 3);\filldraw[black] (v6) circle(4pt);     \node[right] at (v6)   {\LARGE${v_1}$};

      \coordinate (v9) at (0, 4);  \filldraw[black] (v9) circle(4pt);     \node[right] at (v9)     {\LARGE${v_5}$};
      \coordinate (v10) at (0,-1);    \filldraw[black] (v10) circle(4pt);     \node[below left] at (v10) {\LARGE${v_{2}}$};

      \coordinate (v11) at (1, .7);    \filldraw[black] (v11) circle(4pt);     \node[right] at (v11) {\LARGE${v_{x}}$};
      \coordinate (v12) at (-0.5,0.35);    \filldraw[black] (v12) circle(4pt);     \node[below right] at (v12) {\LARGE${v_{y}}$};

     \draw [line width=1.5,black](v10) -- (v1) -- (v3) -- (v5) -- (v9) -- (v6) -- (v4) -- (v2) -- (v10);
     \draw [line width=1.5,black](v1) -- (v4);
     \draw [line width=1.5,black](v5) -- (v6);
     \draw [line width=1.5,black](v3) -- (v2);

\draw [line width=1.5, black] plot [smooth,  tension=.9] coordinates {(0,4) (-2,3) (-2,0) (0,-1)};

     \draw [line width=1.5,black](v11) to[bend left] (v12);

        \filldraw[white] (0,-1) circle(2pt);

       \filldraw[white] (-1, 0) circle(2pt);
       \filldraw[white] (1, 0) circle(2pt);

       \filldraw[white] (-1, 1.5) circle(2pt);
       \filldraw[white] (1, 1.5) circle(2pt);

       \filldraw[white] (-1, 3) circle(2pt);
       \filldraw[white] (1, 3) circle(2pt);

       \filldraw[white] (1, 4) circle(2pt);
       \filldraw[white] (0, 4) circle(2pt);

\end{tikzpicture}
}

\centering  (b2) and $v_{n-2}\in \{x_1, x_2, y_1\}$

\end{minipage}
\begin{minipage}{.3\textwidth}              
\resizebox{4.cm}{!}{
\begin{tikzpicture}
\coordinate (center) at (0,0);
   \def\radius{2.2cm}
   \foreach \x in {0, 60,...,360} {
               }


       \coordinate (v1) at (-1, 0);\filldraw[black] (v1) circle(4pt);    \node[left] at (v1)    {\LARGE${v_6}$};
       \coordinate (v2) at (1,0);\filldraw[black] (v2) circle(4pt);      \node[right] at (v2)   {\LARGE${v_{8}}$};
       \coordinate (v3) at (-1, 1.5);\filldraw[black] (v3) circle(4pt);    \node[left] at (v3)    {\LARGE${v_4}$};
       \coordinate (v4) at (1, 1.5);\filldraw[black] (v4) circle(4pt);     \node[right] at (v4)   {\LARGE${v_3}$};
       \coordinate (v5) at (-1, 3);\filldraw[black] (v5) circle(4pt);    \node[left] at (v5)    {\LARGE${v_7}$};
       \coordinate (v6) at (1, 3);\filldraw[black] (v6) circle(4pt);     \node[right] at (v6)   {\LARGE${v_1}$};

      \coordinate (v9) at (0, 4);  \filldraw[black] (v9) circle(4pt);     \node[right] at (v9)     {\LARGE${v_5}$};
      \coordinate (v10) at (0,-1);    \filldraw[black] (v10) circle(4pt);     \node[below left] at (v10) {\LARGE${v_{2}}$};

      \coordinate (v11) at (-.5, 3.5);    \filldraw[black] (v11) circle(4pt);     \node[right] at (v11) {\LARGE${v_{y}}$};
      \coordinate (v12) at (-1,2.35);    \filldraw[black] (v12) circle(4pt);     \node[left] at (v12) {\LARGE${v_{x}}$};

     \draw [line width=1.5,black](v10) -- (v1) -- (v3) -- (v5) -- (v9) -- (v6) -- (v4) -- (v2) -- (v10);
     \draw [line width=1.5,black](v1) -- (v4);
     \draw [line width=1.5,black](v5) -- (v6);
     \draw [line width=1.5,black](v3) -- (v2);

\draw [line width=1.5, black] plot [smooth,  tension=.9] coordinates {(0,4) (-2,3) (-2,0) (0,-1)};

     \draw [line width=1.5,black](v11) to[bend left] (v12);

        \filldraw[white] (0,-1) circle(2pt);

       \filldraw[white] (-1, 0) circle(2pt);
       \filldraw[white] (1, 0) circle(2pt);

       \filldraw[white] (-1, 1.5) circle(2pt);
       \filldraw[white] (1, 1.5) circle(2pt);

       \filldraw[white] (-1, 3) circle(2pt);
       \filldraw[white] (1, 3) circle(2pt);

       \filldraw[white] (1, 4) circle(2pt);
       \filldraw[white] (0, 4) circle(2pt);

\end{tikzpicture}
}

\centering  (b2) and $v_{n-2}\notin \{x_1, x_2, y_1\}$

\end{minipage}
\begin{minipage}{.3\textwidth}              
\resizebox{4.cm}{!}{
\begin{tikzpicture}
\coordinate (center) at (0,0);
   \def\radius{2.2cm}
   \foreach \x in {0, 60,...,360} {
               }


       \coordinate (v1) at (-1, 0);\filldraw[black] (v1) circle(4pt);    \node[left] at (v1)    {\LARGE${v_6}$};
       \coordinate (v2) at (1,0);\filldraw[black] (v2) circle(4pt);      \node[right] at (v2)   {\LARGE${v_{8}}$};
       \coordinate (v3) at (-1, 1.5);\filldraw[black] (v3) circle(4pt);    \node[left] at (v3)    {\LARGE${v_4}$};
       \coordinate (v4) at (1, 1.5);\filldraw[black] (v4) circle(4pt);     \node[right] at (v4)   {\LARGE${v_3}$};
       \coordinate (v5) at (-1, 3);\filldraw[black] (v5) circle(4pt);    \node[left] at (v5)    {\LARGE${v_7}$};
       \coordinate (v6) at (1, 3);\filldraw[black] (v6) circle(4pt);     \node[right] at (v6)   {\LARGE${v_1}$};

      \coordinate (v9) at (0, 4);  \filldraw[black] (v9) circle(4pt);     \node[right] at (v9)     {\LARGE${v_5}$};
      \coordinate (v10) at (0,-1);    \filldraw[black] (v10) circle(4pt);     \node[below left] at (v10) {\LARGE${v_{2}}$};

      \coordinate (v11) at (-1, .75);    \filldraw[black] (v11) circle(4pt);     \node[left] at (v11) {\LARGE${v_{y}}$};
      \coordinate (v12) at (-1, 2.35);    \filldraw[black] (v12) circle(4pt);     \node[left] at (v12) {\LARGE${v_{x}}$};

     \draw [line width=1.5,black](v10) -- (v1) -- (v3) -- (v5) -- (v9) -- (v6) -- (v4) -- (v2) -- (v10);
     \draw [line width=1.5,black](v1) -- (v4);
     \draw [line width=1.5,black](v5) -- (v6);
     \draw [line width=1.5,black](v3) -- (v2);

\draw [line width=1.5, black] plot [smooth,  tension=.9] coordinates {(0,4) (-2,3) (-2,0) (0,-1)};

     \draw [line width=1.5,black](v11) to[bend right] (v12);

        \filldraw[white] (0,-1) circle(2pt);

       \filldraw[white] (-1, 0) circle(2pt);
       \filldraw[white] (1, 0) circle(2pt);

       \filldraw[white] (-1, 1.5) circle(2pt);
       \filldraw[white] (1, 1.5) circle(2pt);

       \filldraw[white] (-1, 3) circle(2pt);
       \filldraw[white] (1, 3) circle(2pt);

       \filldraw[white] (1, 4) circle(2pt);
       \filldraw[white] (0, 4) circle(2pt);

\end{tikzpicture}
}

\centering  (b2) and $v_{n-2}\notin \{x_1, x_2, y_1\}$

\end{minipage}
\end{center}
\caption{$G\/$ has no $2\/$-cycle and $G_e\/$ has no $2\/$-cycle (b)}    \label{Case(1) of labeling (1.1)(b)}
\end{figure}

\vspace{4mm} (b) Suppose $v_1 \not \in \{x_1, x_2, y_1, y_2\}\/$.     See Figure \ref{Case(1) of labeling (1.1)(b)}.

 \vspace{5mm}   (b1) Suppose  $\{x_1, x_2\} \cap \{y_1, y_2\} = \emptyset\/$.

 \vspace{5mm}     Let $v \in \{x_1, x_2, y_1, y_2\}\/$ be such that $v\/$ is not adjacent to $v_{n-2}\/$ in $G_e\/$. Without loss of generality, assume that $v=x_1\/$.
In $G_e\/$, change the label for $v\/$ from $v_i\/$ to $v_{n-1}\/$.  Extend this labeling on $G_e\/$ to a labeling of $G\/$ by assigning $x\/$ with $v_i\/$ and $y\/$ with $v_n\/$.

 \vspace{5mm}   (b2) Suppose  $\{x_1, x_2\} \cap \{y_1, y_2\} \neq \emptyset\/$.

 \vspace{5mm}   Assume that $x_2 = y_2\/$. Then $y_1 \neq x_1\/$.

 \vspace{5mm}   Suppose that $v_{n-2}\in \{x_1, x_2, y_1\}$.

 \vspace{5mm}   When $v_{n-2}=x_2$, then in $G_e\/$, change the label of $x_1$ from $v_{i}\/$ to $v_{n-1}\/$.
  In $G\/$, assign $v_i , v_{n}\/$ to $x, y\/$ respectively.

 \vspace{5mm}    When $v_{n-2}\neq x_2$, then assume without loss of generality that $v_{n-2} =x_1\/$.
In $G_e\/$, change the label for $y_1\/$ from $v_i\/$ to $v_{n}\/$ and change the label for $x_2\/$ from $v_j\/$ to $v_{n-1}\/$.  Extend this labeling on $G_e\/$ to a labeling of $G\/$ by assigning $y\/$ with $v_i\/$ and $x\/$ with $v_j\/$.

  \vspace{5mm}   Hence assume that $v_{n-2} \not\in \{x_1, x_2, y_1\}$.

  \vspace{5mm}   Suppose $x_2=v_{n-3}$.  Let $v \in \{x_1, y_1\}\/$, say $v=y_1\/$ be such that $v\/$ is not labeled with $v_{3}\/$. Then  in $G_e\/$, change the label for $v_{2}\/$ to $v_{n}\/$.
 In $G\/$, assign $v_{n-1}, v_{2}\/$ to $x, y\/$ respectively.

 \vspace{5mm}   Suppose $ x_2\neq v_{n-3}$. Let $v \in \{x_1, y_1\}\/$, say $v=y_1\/$ be such that $v\/$ is not labeled with $v_{n-3}\/$.  Then in $G_e\/$, change the labels for $v_{n-2}\/$ to $v_{n-1}\/$.   In $G\/$,  assign $v_{n}, v_{n-2}\/$ to $x, y\/$ respectively.

 \vspace{5mm}   (1.2) $G_e$   has a $2\/$-cycle.   See Figure \ref{Case(1) of labeling (1.2)}.

 \vspace{5mm}   Here we may assume that  $x_1x_2x_1\/$ is the $2\/$-cycle in $G_e\/$.  By induction, $G_e\/$ has a labeling $v_1, v_2, \ldots, v_{n-2}\/$ that satisfies condition (ii) of the lemma.   There are two cases to consider.

  \vspace{5mm}   (a)   $\{x_1, x_2\} \cap \{y_1, y_2\} = \emptyset\/$.

  \vspace{5mm}   We may assume without loss of generality  that  $x_1 = v_1\/$ and $x_2 = v_a\/$.

   \vspace{5mm}
 Let $v \in \{y_1, y_2\}\/$, say $v=y_1\/$ be such that $v\/$ is not labeled with $v_{n-2}\/$. In $G_e\/$, change the label for $y_1\/$ from $v_i\/$ to $v_n\/$. Extend this labeling of $G_e\/$ to a labeling of $G\/$ by assigning the label $v_{n-1}\/$ and $v_i\/$ to $x\/$ and $y\/$ respectively.

  \vspace{5mm}   (b)   $\{x_1, x_2\} \cap \{y_1, y_2\} \neq \emptyset\/$.

    \vspace{5mm}   (b1)   $\{x_1, x_2\} = \{y_1, y_2\} \/$.

 \vspace{5mm}   We may assume that $x_1 = v_1, x_2 = v_a\/$ and $x_3 =v_b\/$ where $x_3\/$ is the other neighbor of $x_1\/$ in $G_e\/$. To obtain a required labeling for $G\/$, we first change the label of a vertex in $G_e\/$ from $v_c\/$ to $v_n\/$ where $c \/$ is an odd integer and $c \not \in \{1, a, b\}\/$, and then assign $v_{n-1}, v_c\/$ to $x, y\/$ respectively.

    \vspace{5mm}   (b2)   $x_2 =y_2\/$  and  $x_1 \neq y_1\/$.

    \vspace{5mm}    Then either  (b2.1) $x_2 = v_1\/$ and $x_1 = v_a\/$  or (b2.2) $x_1 = v_1\/$ and $x_2 = v_a\/$.

    \vspace{5mm}   (b2.1) In this case, $y_1 = v_b\/$.  In $G_e\/$, change the label of $y_1\/$ to $v_n\/$. Then extend this new labeling on $G_e\/$ to a required labeling of $G\/$ by assigning $v_{n-1}, v_b\/$ to $x, y\/$ respectively.

 \vspace{5mm}   (b2.2)  In this case, we consider the label on $y_1\/$.

 \vspace{5mm}   Suppose  $y_1 \neq v_{n-2}\/$.   First change the label of a vertex in $G_e\/$ from $v_c\/$ to $v_n\/$ where $c \/$ is an odd integer and $c \not \in \{1, a, b\}\/$, and then assign $v_{n-1}, v_c\/$ to $x, y\/$ respectively.

 \vspace{5mm}   Suppose  $y_1 =  v_{n-2}\/$. Then in $G_e\/$, if $a=3\/$, then change $v_6\/$ to $v_n\/$;  if $a \neq 3\/$, then change $v_2\/$ to $v_n\/$. In any case, in $G\/$, label $x\/$ with $v_{n-1}\/$, and label $y\/$ with   $v_6\/$ and $v_2\/$ respectively.
This gives a required labeling for $G\/$.

\vspace{3mm}{\bf{Case (2)}    } $G$   has a $2\/$-cycle.    See Figure \ref{Case(2) of labeling (2.1)}.

 \vspace{5mm}    Let $xyx\/$ be the $2\/$-cycle of $G\/$. Also, let $x_1\/$ (respectively $y_1\/$) be the other neighbor of $x\/$ (respectively $y\/$).

  \vspace{5mm}   Here let $e = xy\/$. Then $G_e\/$ is a connected $3\/$-regular graph on $n-2\/$ vertices having at most one $2\/$-cycle.

 \vspace{5mm}   (2.1) $G_e$   has no $2\/$-cycle.   In this case, $x_1\/$ and $y_1\/$ are not adjacent in $G\/$.

 \vspace{5mm}   If  $G_e \/$ is the cube $Q_3\/$, then $G\/$ the $3\/$-regular graph depicted in Figure \ref{cubic G on 10 vertices has multiple edges}. which has  a labeling that satisfies condition (ii) of the lemma.


\begin{figure}[htp]
\begin{center}
\begin{minipage}{.37\textwidth}
\resizebox{4cm}{!}{
\begin{tikzpicture}
\coordinate (center) at (0,0);
   \def\radius{2.2cm}
   \foreach \x in {0, 60,...,360} {
               }


  \draw [line width=1.5, black] plot [smooth,  tension=.9] coordinates {(-1.5, 3) (-2.2, 2.4)  (-2.2, .6) (-1.5, 0) };

  \draw [line width=1.5, black] plot [smooth,  tension=.9] coordinates {(1.5, 3) (2.2, 2.4)  (2.2, .6) (1.5, 0) };

       \coordinate (v1) at (-1.5, 0);\filldraw[black] (v1) circle(4pt);\node[below left] at (v1)   {\LARGE${v_2}$};
       \coordinate (v2) at (1.5,0);\filldraw[black] (v2) circle(4pt);\node[below right] at (v2)    {\LARGE${v_{7}}$};
       \coordinate (v3) at (-1.5, 1);\filldraw[black] (v3) circle(4pt);\node[left] at (v3)   {\LARGE${v_8}$};
       \coordinate (v4) at (1.5, 1);\filldraw[black] (v4) circle(4pt);\node[right] at (v4)   {\LARGE${v_5}$};
       \coordinate (v5) at (-1.5, 2);\filldraw[black] (v5) circle(4pt);\node[left] at (v5)   {\LARGE${v_6}$};
       \coordinate (v6) at (1.5, 2);\filldraw[black] (v6) circle(4pt);\node[right] at (v6)   {\LARGE${v_3}$};
       \coordinate (v7) at (-1.5, 3);\filldraw[black] (v7) circle(4pt);\node[above left] at (v7)   {\LARGE${v_4}$};
       \coordinate (v8) at (1.5, 3);\filldraw[black] (v8) circle(4pt);\node[above right] at (v8)   {\LARGE${v_{1}}$};

       \coordinate (v9) at (-1, 3.7);\filldraw[white] (v9) circle(4pt);\node[left] at (v9)   {\white \LARGE${v_{9}}$};
       \coordinate (v10) at (1, 3,7);\filldraw[white] (v10) circle(4pt);\node[right] at (v10)   {\white \LARGE${v_{1}}$};

     \draw [line width=1.5,black](v1) -- (v3) -- (v5) -- (v7) -- (v8) -- (v6) -- (v4) -- (v2) -- (v1);
     \draw [line width=1.5,black](v3) -- (v4);
     \draw [line width=1.5,black](v5) -- (v6);

       \filldraw[white] (-1.5, 0) circle(2pt);
       \filldraw[white] (1.5, 0) circle(2pt);
       \filldraw[white] (-1.5, 1) circle(2pt);
       \filldraw[white] (1.5, 1) circle(2pt);
       \filldraw[white] (-1.5, 2) circle(2pt);
       \filldraw[white] (1.5, 2) circle(2pt);
       \filldraw[white] (-1.5, 3) circle(2pt);
       \filldraw[white] (1.5, 3) circle(2pt);

       \filldraw[white] (v6) circle(3pt);
       \filldraw[white] (v7) circle(3pt);
       \filldraw[white] (v8) circle(3pt);

       \filldraw[black] (v6) circle(2pt);
       \filldraw[black] (v7) circle(2pt);
       \filldraw[black] (v8) circle(2pt);

\end{tikzpicture}
}
\caption{$G_e$ is a cube $Q_3\/$}     \label{cube Q3}
\end{minipage}
\begin{minipage}{.37\textwidth}
\resizebox{4cm}{!}{
\begin{tikzpicture}
\coordinate (center) at (0,0);
   \def\radius{2.2cm}
   \foreach \x in {0, 60,...,360} {
               }


  \draw [line width=1.5, black] plot [smooth,  tension=.9] coordinates {(-1.5, 3) (-2.2, 2.4)  (-2.2, .6) (-1.5, 0) };

  \draw [line width=1.5, black] plot [smooth,  tension=.9] coordinates {(1.5, 3) (2.2, 2.4)  (2.2, .6) (1.5, 0) };

       \coordinate (v1) at (-1.5, 0);\filldraw[black] (v1) circle(4pt);\node[below left] at (v1)   {\LARGE${v_2}$};
       \coordinate (v2) at (1.5,0);\filldraw[black] (v2) circle(4pt);\node[below right] at (v2)    {\LARGE${v_{7}}$};
       \coordinate (v3) at (-1.5, 1);\filldraw[black] (v3) circle(4pt);\node[left] at (v3)   {\LARGE${v_8}$};
       \coordinate (v4) at (1.5, 1);\filldraw[black] (v4) circle(4pt);\node[right] at (v4)   {\LARGE${v_5}$};
       \coordinate (v5) at (-1.5, 2);\filldraw[black] (v5) circle(4pt);\node[left] at (v5)   {\LARGE${v_6}$};
       \coordinate (v6) at (1.5, 2);\filldraw[black] (v6) circle(4pt);\node[right] at (v6)   {\Large${v_{10}}$};
       \coordinate (v7) at (-1.5, 3);\filldraw[black] (v7) circle(4pt);\node[above left] at (v7)   {\LARGE${v_4}$};
       \coordinate (v8) at (1.5, 3);\filldraw[black] (v8) circle(4pt);\node[above right] at (v8)   {\LARGE${v_{3}}$};

       \coordinate (v9) at (-1.5, 3.7);\filldraw[black] (v9) circle(4pt);\node[above] at (v9)   {\LARGE${v_{9}}$};
       \coordinate (v10) at (1.5, 3.7);\filldraw[black] (v10) circle(4pt);\node[above] at (v10)   {\LARGE${v_{1}}$};

      \draw [line width=1.5,black] (v10)-- (v8) -- (v6) -- (v4) -- (v2) -- (v1) -- (v3) -- (v5) -- (v7)--(v9);
     \draw [line width=1.5,black](v3) -- (v4);
     \draw [line width=1.5,black](v5) -- (v6);

     \draw [line width=1.5,black](v9) to[bend right] (v10);
     \draw [line width=1.5,black](v10) to[bend right] (v9);

       \filldraw[white] (-1.5, 0) circle(2pt);
       \filldraw[white] (1.5, 0) circle(2pt);
       \filldraw[white] (-1.5, 1) circle(2pt);
       \filldraw[white] (1.5, 1) circle(2pt);
       \filldraw[white] (-1.5, 2) circle(2pt);
       \filldraw[red] (1.5, 2) circle(2pt);
       \filldraw[white] (-1.5, 3) circle(2pt);
       \filldraw[red] (1.5, 3) circle(2pt);

       \filldraw[red] (-1.5, 3.7) circle(2pt);
       \filldraw[red] (1.5, 3.7) circle(2pt);

       \filldraw[red] (v10) circle(3pt);
       \filldraw[white] (v6) circle(3pt);
       \filldraw[white] (v7) circle(3pt);
       \filldraw[white] (v8) circle(3pt);

       \filldraw[black] (v10) circle(2pt);
       \filldraw[black] (v6) circle(2pt);
       \filldraw[black] (v7) circle(2pt);
       \filldraw[black] (v8) circle(2pt);

\end{tikzpicture}
}
\caption{The cubic $G$}    \label{cubic G on 10 vertices has multiple edges}

\end{minipage}

\end{center}
\end{figure}

 \vspace{5mm}   Hence we assume that  $G_e\/$ is not the cube $Q_3\/$.

 \vspace{5mm}   By induction, $V(G_e)$ can be labeled as $v_1, v_2, \ldots, v_{n-2}$ which satisfies condition (i) of Lemma \ref{4labeling}. We shall use this labeling on $V(G_e)\/$ to obtain a labeling for $G\/$ that satisfies the conditions of Lemma \ref{4labeling}.

 \vspace{3mm} (a) Suppose $v_1 \in \{x_1,  y_1\}\/$.

 \vspace{5mm}   Assume without loss of generality that $x_1 = v_1\/$ and $y_1 = v_a\/$. In $G_e\/$, change the labels of $x_1, y_1\/$ to $v_{n-1}, v_n\/$ respectively.  Extend this new labeling on $V(G_e)\/$ to a required labeling of $G\/$ by labeling $x, y\/$ with $v_1, v_a\/$ respectively.

\vspace{3mm} (b) Suppose $v_1 \not \in \{x_1,  y_1\}\/$.

 \vspace{5mm}   (b1) $v_{n-2} \in \{x_1, y_1\}\/$.

 \vspace{5mm}   Suppose $x_1 = v_{n-2}\/$.

 \vspace{5mm}   If no vertex in $N_{G_e}(v_{n-2})-y_1\/$ is $v_{n-4}\/$, then in $G_e\/$, change the labels for $v_1\/$ to $v_n\/$, and interchange the labels of $v_{n-3}\/$ and $v_{n-2}\/$. Now, extend this new labeling on $G_e\/$ to a required labeling in $G\/$ by assigning $v_1 , v_{n-1}\/$ to $x, y\/$ respectively.

 \vspace{5mm}   If some vertex in $N_{G_e}(v_{n-2})-y_1\/$ is $v_{n-4}\/$, then in $G_e\/$, change the labels for $x_1, v_1, v_{n-3}\/$ to $v_{n-1}, v_n, v_{n-2}\/$ respectively.  Now, extend this new labeling on $G_e\/$ to a required labeling in $G\/$ by assigning $v_1 , v_{n-3}\/$ to $x, y\/$ respectively.


\vspace{3mm} (b2) $v_{n-2} \not \in \{x_1, y_1\}\/$.

 \vspace{5mm}   Suppose the label of $x_1 \/$ is $v_i\/$ where $i\/$ is odd. Then in $G_e\/$, change the label of $v_1\/$ to $v_n\/$. In $G\/$, assign $v_1, v_{n-1}\/$ to $x, y\/$ respectively.

 \vspace{5mm}   Hence assume that $x_1 =v_i, y_1 = v_j\/$ are such that $i\/$ and $j\/$ are both even.

 \vspace{5mm}   Suppose $2\in \{i, j\}$. Assume without loss of generality that $i=2$.

 \vspace{5mm}    If $v_4 \notin N_{G_e}(x_1)-y_1\/$, then in $G_e$, interchange the labels of $v_{2}\/$ and $v_{3}\/$. Also, change the label of $v_1\/$ to $v_n\/$. Now, extend this new labeling on $G_e\/$ to a required labeling in $G\/$ by assigning $v_1 , v_{n-1}\/$ to $x, y\/$ respectively.

 \vspace{5mm}   If $v_4 \in N_{G_e}(x_1)-y_1\/$, then in $G_e\/$, change the label of $x_1\/$ from $v_{2}\/$ to $v_{n-1}\/$. Also, change the label of $v_1, v_3\/$ to $v_n, v_2\/$ respectively. Now, extend this new labeling on $G_e\/$ to a required labeling in $G\/$ by assigning $v_1 , v_{3}\/$ to $x, y\/$ respectively.

 \vspace{5mm}   Hence assume that $2\notin \{i, j\}$ and that  $i<j$.

 \vspace{5mm}   Let $u$ be the vertex in $G_e$ having the label $v_{i+1}$. Note that such a vertex exists because $i < n-2$.

 \vspace{5mm}   To extend the labeling of $G_e\/$ to a required labeling for $G\/$, we first  interchange the labels of $v_i\/$ and $v_{i+1}\/$ (in $G_e\/$). There are two cases to consider.

\vspace{3mm} Now, if $u$ is not adjacent to a vertex with the label $v_{i-1}$ in $G_e$, then change  the label for $v_1\/$  to $v_n\/$. In $G\/$, we  assign  $v_1 , v_{n-1}\/$ to $x, y\/$ respectively.

 \vspace{5mm}   If $u$ is adjacent to a vertex with the label $v_{i-1}$ in $G_e$, then  change the label for $v_1, v_{i-1}\/$ to $v_{n-1}, v_n\/$ respectively.
  In $G\/$, assign $v_1 , v_{i-1}\/$ to $x, y\/$ respectively.

 \vspace{5mm}   (2.2) $G_e$   has a $2\/$-cycle.

 \vspace{5mm}   In this case, $x_1x_2 \in E(G)\/$.   By induction, $G_e\/$ has a labeling $v_1, v_2, \ldots, v_{n-2}\/$ that satisfies condition (ii) of the lemma. To obtain a required labeling for $G\/$, we first change the label of the  vertex $v\/$ (which is adjacent to $x_1\/$) in $G_e\/$ from $v_b\/$ to $v_n\/$, and then extend this new labeling to $G\/$ by assigning $v_1, v_{n-1}, v_b\/$ to $x, y , x_1\/$ respectively.

 \vspace{5mm}   This completes the proof.  $\qed$


\begin{result} \label{4-angulation}
Let $G\/$ be a $3\/$-regular graph with $n =4 + 2t\/$ vertices where $t \geq 2\/$. Then $G\/ $ is potentially $4$-angulable.
\end{result}

 \vspace{5mm}    \noindent
{\bf Proof:}   Suppose $G\/$ is the cube $Q_3\/$ and has labeling as shown in Figure \ref{cube Q3}. Then $K_n-G\/ $ admits a $4$-angulation that has the diagonal edges $v_1v_6\/$ and $v_2v_5\/$.
 Now suppose that $G\/$ is a cubic not $Q_3\/$. By Lemma \ref{labeling}, $V(G)\/$ admits a labeling as described in Lemma \ref{4labeling}.  Then $G\/ $ is potentially $4$-angulable.        $\qed$

\break

\break


\begin{figure}[htp]
\begin{center}
\begin{minipage}{.5\textwidth}
\begin{center}
\resizebox{6.5cm}{!}{
\begin{tikzpicture}[rotate=-36,.style={draw}]
\coordinate (center) at (0,0);
   \def\radius{2cm}
   \foreach \x in {0, 18,...,360} {
             \filldraw[] (\x:2cm) circle(1pt);
               }


      \coordinate (v0) at (-1.62, 1.18);\filldraw[blue] (v0) circle(1pt); 

       \coordinate (v1) at (-1.18, 1.62);\filldraw[blue] (v1) circle(1pt);
       \coordinate (v2) at (-.62, 1.9);\filldraw[blue] (v2) circle(1pt);  

       \coordinate (v3) at (0, 2);\filldraw[blue] (v3) circle(1pt);
       \coordinate (v4) at (.62, 1.9);\filldraw[blue] (v4) circle(1pt);

       \coordinate (v5) at (1.18, 1.62);\filldraw[blue] (v5) circle(1pt);
       \coordinate (v6) at (1.62, 1.18);\filldraw[blue] (v6) circle(1pt);\node[right] at (v6) {${v_{1}}$};

       \coordinate (v7) at (1.89, .62);\filldraw[blue] (v7) circle(1pt);\node[below right] at (v7) {${v_2}$};
      \coordinate (v8) at (2, 0);\filldraw[blue] (v8) circle(1pt);

       \coordinate (v9) at (1.89, -.62);\filldraw[black] (v9) circle(1pt);
       \coordinate (v10) at (1.62, -1.18);\filldraw[black] (v10) circle(1pt);

       \coordinate (v11) at (1.18, -1.62);\filldraw[black] (v11) circle(1pt);

       \coordinate (v12) at (.62, -1.9);\filldraw[black] (v12) circle(1pt);

       \coordinate (v13) at (0, -2);\filldraw[black] (v13) circle(1pt);

       \coordinate (v14) at (-.62, -1.9);\filldraw[black] (v14) circle(1pt);

       \coordinate (v15) at (-1.18, -1.62);\filldraw[black] (v15) circle(1pt);

       \coordinate (v16) at (-1.62, -1.18);\filldraw[black] (v16) circle(1pt);

       \coordinate (v17) at (-1.89, -.62);\filldraw[black] (v17) circle(1pt);

       \coordinate (v18) at (-2, 0);\filldraw[black] (v18) circle(1pt);

       \coordinate (v19) at (-1.89, .62);\filldraw[black] (v19) circle(1pt);

     \draw [line width=1,black](v6)-- (v13);
     \draw [dashed,white](v6)-- (v13);

     \draw [line width=1,black](v16)-- (v12);

     \draw [line width=1,black](v3)-- (v7);

     \draw [line width=1.3,black](v16)-- (v0);
     \draw [line width=1.3,black](v17)-- (v1);
     \draw [line width=1.3,black](v18)-- (v2)-- (v9);

     \draw [line width=1.3,black](v19)-- (v3);
     \draw [dashed, white](v19)-- (v3);

    \draw [line width=1, black](v0)-- (v4);
    \draw [line width=1, black](v1)-- (v5);

    \draw [line width=1.3,black](v10)-- (v0);
    \draw [line width=1, black](v11)-- (v1);
     \draw [line width=1,black](v12)-- (v2);

     \draw [line width=1, black](v7)-- (v0);
     \draw [line width=1, black](v8)-- (v1);
     \draw [dashed,black](v7)-- (v0);
     \draw [dashed,black](v8)-- (v1);
     \draw [dashed,black](v5)-- (v1);
     \draw [dashed,black](v4)-- (v0);

     \draw [line width=1.3,orange](v6) -- (v19);

     \draw [line width=1.3,green](v13)-- (v3);
     \draw [line width=1, black](v15)-- (v19);
     \draw [line width=1, black](v14)-- (v18);

    \filldraw[black] (v0) circle(2.5pt);
    \filldraw[black] (v1) circle(2.5pt);
    \filldraw[black] (v2) circle(2.5pt);
   \filldraw[black] (v3) circle(2.5pt);
    \filldraw[black] (v4) circle(2.5pt);
   \filldraw[black] (v5) circle(2.5pt);
    \filldraw[black] (v6) circle(2.5pt);
    \filldraw[black] (v7) circle(2.5pt);
   \filldraw[black] (v8) circle(2.5pt);
    \filldraw[black] (v9) circle(2.5pt);
   \filldraw[black] (v10) circle(2.5pt);
    \filldraw[black] (v11) circle(2.5pt);
    \filldraw[black] (v12) circle(2.5pt);
    \filldraw[black] (v13) circle(2.5pt);
    \filldraw[black] (v14) circle(2.5pt);
    \filldraw[black] (v15) circle(2.5pt);
    \filldraw[black] (v16) circle(2.5pt);
    \filldraw[black] (v17) circle(2.5pt);
    \filldraw[black] (v18) circle(2.5pt);
    \filldraw[black] (v19) circle(2.5pt);

    \filldraw[red] (v18) circle(1.5pt);
    \filldraw[red] (v19) circle(1.5pt);
    \filldraw[red] (v0) circle(1.5pt);
    \filldraw[red] (v1) circle(1.5pt);
    \filldraw[red] (v2) circle(1.5pt);
    \filldraw[red] (v3) circle(1.5pt);

    \filldraw[white] (v4) circle(1.5pt);
    \filldraw[white] (v5) circle(1.5pt);
    \filldraw[white] (v9) circle(1.5pt);
    \filldraw[white] (v7) circle(1.5pt);
    \filldraw[white] (v8) circle(1.5pt);

    \filldraw[white] (v10) circle(1.5pt);
    \filldraw[white] (v11) circle(1.5pt);
    \filldraw[white] (v12) circle(1.5pt);

    \filldraw[white] (v14) circle(1.5pt);
    \filldraw[white] (v15) circle(1.5pt);
    \filldraw[white] (v16) circle(1.5pt);
    \filldraw[white] (v17) circle(1.5pt);

    \filldraw[white] (v4) circle(1.5pt);
    \filldraw[yellow] (v5) circle(1.5pt);
    \filldraw[orange] (v6) circle(1.5pt);
    \filldraw[yellow] (v7) circle(1.5pt);
    \filldraw[white] (v8) circle(1.5pt);

    \filldraw[green] (v11) circle(1.5pt);
    \filldraw[green] (v12) circle(1.5pt);
    \filldraw[green] (v13) circle(1.5pt);
    \filldraw[white] (v14) circle(1.5pt);
    \filldraw[white] (v15) circle(1.5pt);

\end{tikzpicture}
}

\centering {${\cal J}_{20,5}(^*_2)\/$ has an $\alpha$-pair of vertices  in $F_m={\cal J}_{20,5}(^*_2)-\{v_{q+1}, \ldots, v_{q+g-2}\}\/$ for some $q\in\{1, 2, 4\}\/$}

\end{center}
\end{minipage}

\vspace{15mm}

\begin{minipage}{.33\textwidth}
\resizebox{5.cm}{!}{
\begin{tikzpicture}[rotate=-36,.style={draw}]
\coordinate (center) at (0,0);
   \def\radius{2cm}
   \foreach \x in {0, 18,...,360} {
             \filldraw[] (\x:2cm) circle(1pt);
               }


      \coordinate (v0) at (-1.62, 1.18);\filldraw[blue] (v0) circle(1pt); 

       \coordinate (v1) at (-1.18, 1.62);\filldraw[blue] (v1) circle(1pt);
       \coordinate (v2) at (-.62, 1.9);\filldraw[blue] (v2) circle(1pt);  

       \coordinate (v3) at (0, 2);\filldraw[blue] (v3) circle(1pt);
       \coordinate (v4) at (.62, 1.9);\filldraw[blue] (v4) circle(1pt);

       \coordinate (v5) at (1.18, 1.62);\filldraw[blue] (v5) circle(1pt);
       \coordinate (v6) at (1.62, 1.18);\filldraw[blue] (v6) circle(1pt);\node[right] at (v6) {${v_{1}}$};
       \coordinate (v7) at (1.89, .62);\filldraw[blue] (v7) circle(1pt);\node[below right] at (v7) {${v_2}$};
      \coordinate (v8) at (2, 0);\filldraw[blue] (v8) circle(1pt);\node[below right] at (v8) {${v_3}$};

       \coordinate (v9) at (1.89, -.62);\filldraw[black] (v9) circle(1pt);\node[below right] at (v9) {${v_{4}}$};
       \coordinate (v10) at (1.62, -1.18);\filldraw[black] (v10) circle(1pt);\node[below right] at (v10) {${v_{5}}$};

       \coordinate (v11) at (1.18, -1.62);\filldraw[black] (v11) circle(1pt);\node[below ] at (v11) {${v_{6}}$};

       \coordinate (v12) at (.62, -1.9);\filldraw[black] (v12) circle(1pt);\node[below ] at (v12) {${v_{7}}$};

       \coordinate (v13) at (0, -2);\filldraw[black] (v13) circle(1pt);

       \coordinate (v14) at (-.62, -1.9);\filldraw[black] (v14) circle(1pt);

       \coordinate (v15) at (-1.18, -1.62);\filldraw[black] (v15) circle(1pt);

       \coordinate (v16) at (-1.62, -1.18);\filldraw[black] (v16) circle(1pt);

       \coordinate (v17) at (-1.89, -.62);\filldraw[black] (v17) circle(1pt);

       \coordinate (v18) at (-2, 0);\filldraw[black] (v18) circle(1pt);

       \coordinate (v19) at (-1.89, .62);\filldraw[black] (v19) circle(1pt);

     \draw [line width=1,black](v6)-- (v13);
     \draw [dashed,white](v6)-- (v13);

     \draw [line width=1,black](v16)-- (v12);

     \draw [line width=1.3,black](v16)-- (v0);
     \draw [line width=1.3,black](v17)-- (v1);
     \draw [line width=1.3,black](v18)-- (v2);

     \draw [line width=1.3,black](v19)-- (v3);
     \draw [dashed, white](v19)-- (v3);

    \draw [line width=1, black](v0)-- (v4);
    \draw [line width=1, black](v1)-- (v5);

    \draw [line width=1.3,black](v10)-- (v0);
    \draw [line width=1, black](v11)-- (v1);
     \draw [line width=1,black](v12)-- (v2);

     \draw [dashed,black](v5)-- (v1);
     \draw [dashed,black](v4)-- (v0);

     \draw [line width=1.3,orange](v6) -- (v19);

     \draw [line width=1.3,green](v13)-- (v3);
     \draw [line width=1, black](v15)-- (v19);
     \draw [line width=1, black](v14)-- (v18);

    \filldraw[black] (v0) circle(2.5pt);
    \filldraw[black] (v1) circle(2.5pt);
    \filldraw[black] (v2) circle(2.5pt);
   \filldraw[black] (v3) circle(2.5pt);
    \filldraw[black] (v4) circle(2.5pt);
   \filldraw[black] (v5) circle(2.5pt);
    \filldraw[black] (v6) circle(2.5pt);
    \filldraw[black] (v7) circle(2.5pt);
   \filldraw[black] (v8) circle(2.5pt);
    \filldraw[black] (v9) circle(2.5pt);
   \filldraw[black] (v10) circle(2.5pt);
    \filldraw[black] (v11) circle(2.5pt);
    \filldraw[black] (v12) circle(2.5pt);
    \filldraw[black] (v13) circle(2.5pt);
    \filldraw[black] (v14) circle(2.5pt);
    \filldraw[black] (v15) circle(2.5pt);
    \filldraw[black] (v16) circle(2.5pt);
    \filldraw[black] (v17) circle(2.5pt);
    \filldraw[black] (v18) circle(2.5pt);
    \filldraw[black] (v19) circle(2.5pt);

    \filldraw[red] (v18) circle(1.5pt);
    \filldraw[red] (v19) circle(1.5pt);
    \filldraw[red] (v0) circle(1.5pt);
    \filldraw[red] (v1) circle(1.5pt);
    \filldraw[red] (v2) circle(1.5pt);
    \filldraw[red] (v3) circle(1.5pt);

    \filldraw[white] (v4) circle(1.5pt);
    \filldraw[white] (v5) circle(1.5pt);
    \filldraw[white] (v9) circle(1.5pt);
    \filldraw[white] (v7) circle(1.5pt);
    \filldraw[white] (v8) circle(1.5pt);

    \filldraw[white] (v10) circle(1.5pt);
    \filldraw[white] (v11) circle(1.5pt);
    \filldraw[white] (v12) circle(1.5pt);

    \filldraw[white] (v14) circle(1.5pt);
    \filldraw[white] (v15) circle(1.5pt);
    \filldraw[white] (v16) circle(1.5pt);
    \filldraw[white] (v17) circle(1.5pt);

    \filldraw[white] (v4) circle(1.5pt);
    \filldraw[yellow] (v5) circle(1.5pt);
    \filldraw[orange] (v6) circle(1.5pt);
    \filldraw[yellow] (v7) circle(1.5pt);
    \filldraw[white] (v8) circle(1.5pt);

    \filldraw[green] (v11) circle(1.5pt);
    \filldraw[green] (v12) circle(1.5pt);
    \filldraw[green] (v13) circle(1.5pt);
    \filldraw[white] (v14) circle(1.5pt);
    \filldraw[white] (v15) circle(1.5pt);

\end{tikzpicture}
}

\centering {(a) $F_m={\cal J}_{20,5}(^*_2)-\{v_2, v_3, v_4\}= {\cal J}_{17,5}(^*_2)\/$}

\end{minipage}
\begin{minipage}{.32\textwidth}
\resizebox{5.cm}{!}{
\begin{tikzpicture}[rotate=-36,.style={draw}]
\coordinate (center) at (0,0);
   \def\radius{2cm}
   \foreach \x in {0, 18,...,360} {
             \filldraw[] (\x:2cm) circle(1pt);
               }


      \coordinate (v0) at (-1.62, 1.18);\filldraw[blue] (v0) circle(1pt); 

       \coordinate (v1) at (-1.18, 1.62);\filldraw[blue] (v1) circle(1pt);
       \coordinate (v2) at (-.62, 1.9);\filldraw[blue] (v2) circle(1pt);  

       \coordinate (v3) at (0, 2);\filldraw[blue] (v3) circle(1pt);
       \coordinate (v4) at (.62, 1.9);\filldraw[blue] (v4) circle(1pt);

       \coordinate (v5) at (1.18, 1.62);\filldraw[blue] (v5) circle(1pt);
       \coordinate (v6) at (1.62, 1.18);\filldraw[blue] (v6) circle(1pt);\node[right] at (v6) {${v_{1}}$};
       \coordinate (v7) at (1.89, .62);\filldraw[blue] (v7) circle(1pt);\node[below right] at (v7) {${v_2}$};
      \coordinate (v8) at (2, 0);\filldraw[blue] (v8) circle(1pt);\node[below right] at (v8) {${v_3}$};

       \coordinate (v9) at (1.89, -.62);\filldraw[black] (v9) circle(1pt);\node[below right] at (v9) {${v_{4}}$};
       \coordinate (v10) at (1.62, -1.18);\filldraw[black] (v10) circle(1pt);\node[below right] at (v10) {${v_{5}}$};

       \coordinate (v11) at (1.18, -1.62);\filldraw[black] (v11) circle(1pt);\node[below ] at (v11) {${v_{6}}$};

       \coordinate (v12) at (.62, -1.9);\filldraw[black] (v12) circle(1pt);\node[below ] at (v12) {${v_{7}}$};

       \coordinate (v13) at (0, -2);\filldraw[black] (v13) circle(1pt);

       \coordinate (v14) at (-.62, -1.9);\filldraw[black] (v14) circle(1pt);

       \coordinate (v15) at (-1.18, -1.62);\filldraw[black] (v15) circle(1pt);

       \coordinate (v16) at (-1.62, -1.18);\filldraw[black] (v16) circle(1pt);

       \coordinate (v17) at (-1.89, -.62);\filldraw[black] (v17) circle(1pt);

       \coordinate (v18) at (-2, 0);\filldraw[black] (v18) circle(1pt);

       \coordinate (v19) at (-1.89, .62);\filldraw[black] (v19) circle(1pt);

     \draw [line width=1,black](v6)-- (v13);
     \draw [dashed,white](v6)-- (v13);

     \draw [line width=1,black](v16)-- (v12);

     \draw [line width=1,black](v3)-- (v7);

     \draw [line width=1.3,black](v16)-- (v0);
     \draw [line width=1.3,black](v17)-- (v1);
     \draw [line width=1.3,black](v18)-- (v2);

     \draw [line width=1.3,black](v19)-- (v3);
     \draw [dashed, white](v19)-- (v3);

    \draw [line width=1, black](v0)-- (v4);
    \draw [line width=1, black](v1)-- (v5);

    \draw [line width=1, black](v11)-- (v1);
     \draw [line width=1,black](v12)-- (v2);

     \draw [line width=1, black](v7)-- (v0);
     \draw [dashed,black](v7)-- (v0);
     \draw [dashed,black](v5)-- (v1);
     \draw [dashed,black](v4)-- (v0);

     \draw [line width=1.3,orange](v6) -- (v19);

     \draw [line width=1.3,green](v13)-- (v3);
     \draw [line width=1, black](v15)-- (v19);
     \draw [line width=1, black](v14)-- (v18);

    \filldraw[black] (v0) circle(2.5pt);
    \filldraw[black] (v1) circle(2.5pt);
    \filldraw[black] (v2) circle(2.5pt);
   \filldraw[black] (v3) circle(2.5pt);
    \filldraw[black] (v4) circle(2.5pt);
   \filldraw[black] (v5) circle(2.5pt);
    \filldraw[black] (v6) circle(2.5pt);
    \filldraw[black] (v7) circle(2.5pt);
   \filldraw[black] (v8) circle(2.5pt);
    \filldraw[black] (v9) circle(2.5pt);
   \filldraw[black] (v10) circle(2.5pt);
    \filldraw[black] (v11) circle(2.5pt);
    \filldraw[black] (v12) circle(2.5pt);
    \filldraw[black] (v13) circle(2.5pt);
    \filldraw[black] (v14) circle(2.5pt);
    \filldraw[black] (v15) circle(2.5pt);
    \filldraw[black] (v16) circle(2.5pt);
    \filldraw[black] (v17) circle(2.5pt);
    \filldraw[black] (v18) circle(2.5pt);
    \filldraw[black] (v19) circle(2.5pt);

    \filldraw[red] (v18) circle(1.5pt);
    \filldraw[red] (v19) circle(1.5pt);
    \filldraw[red] (v0) circle(1.5pt);
    \filldraw[red] (v1) circle(1.5pt);
    \filldraw[red] (v2) circle(1.5pt);
    \filldraw[red] (v3) circle(1.5pt);

    \filldraw[white] (v4) circle(1.5pt);
    \filldraw[white] (v5) circle(1.5pt);
    \filldraw[white] (v9) circle(1.5pt);
    \filldraw[white] (v7) circle(1.5pt);
    \filldraw[white] (v8) circle(1.5pt);

    \filldraw[white] (v10) circle(1.5pt);
    \filldraw[white] (v11) circle(1.5pt);
    \filldraw[white] (v12) circle(1.5pt);

    \filldraw[white] (v14) circle(1.5pt);
    \filldraw[white] (v15) circle(1.5pt);
    \filldraw[white] (v16) circle(1.5pt);
    \filldraw[white] (v17) circle(1.5pt);

    \filldraw[white] (v4) circle(1.5pt);
    \filldraw[yellow] (v5) circle(1.5pt);
    \filldraw[orange] (v6) circle(1.5pt);
    \filldraw[yellow] (v7) circle(1.5pt);
    \filldraw[white] (v8) circle(1.5pt);

    \filldraw[green] (v11) circle(1.5pt);
    \filldraw[green] (v12) circle(1.5pt);
    \filldraw[green] (v13) circle(1.5pt);
    \filldraw[white] (v14) circle(1.5pt);
    \filldraw[white] (v15) circle(1.5pt);

\end{tikzpicture}
}

\centering {(b)  $F_m={\cal J}_{20,5}(^*_2)-\{v_3, v_4, v_5\}={\cal J}_{17,5}(^*_2)\/$}

\end{minipage}
\begin{minipage}{.33\textwidth}
\resizebox{5.cm}{!}{
\begin{tikzpicture}[rotate=-36,.style={draw}]
\coordinate (center) at (0,0);
   \def\radius{2cm}
   \foreach \x in {0, 18,...,360} {
             \filldraw[] (\x:2cm) circle(1pt);
               }


      \coordinate (v0) at (-1.62, 1.18);\filldraw[blue] (v0) circle(1pt); 

       \coordinate (v1) at (-1.18, 1.62);\filldraw[blue] (v1) circle(1pt);
       \coordinate (v2) at (-.62, 1.9);\filldraw[blue] (v2) circle(1pt);  

       \coordinate (v3) at (0, 2);\filldraw[blue] (v3) circle(1pt);
       \coordinate (v4) at (.62, 1.9);\filldraw[blue] (v4) circle(1pt);

       \coordinate (v5) at (1.18, 1.62);\filldraw[blue] (v5) circle(1pt);
       \coordinate (v6) at (1.62, 1.18);\filldraw[blue] (v6) circle(1pt);\node[right] at (v6) {${v_{1}}$};
       \coordinate (v7) at (1.89, .62);\filldraw[blue] (v7) circle(1pt);\node[below right] at (v7) {${v_2}$};
      \coordinate (v8) at (2, 0);\filldraw[blue] (v8) circle(1pt);\node[below right] at (v8) {${v_3}$};

       \coordinate (v9) at (1.89, -.62);\filldraw[black] (v9) circle(1pt);\node[below right] at (v9) {${v_{4}}$};
       \coordinate (v10) at (1.62, -1.18);\filldraw[black] (v10) circle(1pt);\node[below right] at (v10) {${v_{5}}$};

       \coordinate (v11) at (1.18, -1.62);\filldraw[black] (v11) circle(1pt);\node[below ] at (v11) {${v_{6}}$};

       \coordinate (v12) at (.62, -1.9);\filldraw[black] (v12) circle(1pt);\node[below ] at (v12) {${v_{7}}$};

       \coordinate (v13) at (0, -2);\filldraw[black] (v13) circle(1pt);

       \coordinate (v14) at (-.62, -1.9);\filldraw[black] (v14) circle(1pt);

       \coordinate (v15) at (-1.18, -1.62);\filldraw[black] (v15) circle(1pt);

       \coordinate (v16) at (-1.62, -1.18);\filldraw[black] (v16) circle(1pt);

       \coordinate (v17) at (-1.89, -.62);\filldraw[black] (v17) circle(1pt);

       \coordinate (v18) at (-2, 0);\filldraw[black] (v18) circle(1pt);

       \coordinate (v19) at (-1.89, .62);\filldraw[black] (v19) circle(1pt);

     \draw [line width=1.3,black](v6)-- (v13);
     \draw [dashed,white](v6)-- (v13);

     \draw [dashed,white](v6)-- (v2);
1
     \draw [line width=1,black](v3)-- (v7);

     \draw [line width=1.3,black](v16)-- (v0);
     \draw [line width=1.3,black](v17)-- (v1);
     \draw [line width=1.3,black](v18)-- (v2)-- (v9);

     \draw [line width=1.3,black](v19)-- (v3);
     \draw [dashed, white](v19)-- (v3);

    \draw [line width=1, black](v0)-- (v4);
    \draw [line width=1, black](v1)-- (v5);

     \draw [line width=1, black](v7)-- (v0);
     \draw [line width=1, black](v8)-- (v1);
     \draw [dashed,black](v7)-- (v0);
     \draw [dashed,black](v8)-- (v1);
     \draw [dashed,black](v5)-- (v1);
     \draw [dashed,black](v4)-- (v0);

     \draw [line width=1.3,orange](v6) -- (v19);

     \draw [line width=1.3,green](v13)-- (v3);
     \draw [line width=1, black](v15)-- (v19);
     \draw [line width=1, black](v14)-- (v18);

    \filldraw[black] (v0) circle(2.5pt);
    \filldraw[black] (v1) circle(2.5pt);
    \filldraw[black] (v2) circle(2.5pt);
   \filldraw[black] (v3) circle(2.5pt);
    \filldraw[black] (v4) circle(2.5pt);
   \filldraw[black] (v5) circle(2.5pt);
    \filldraw[black] (v6) circle(2.5pt);
    \filldraw[black] (v7) circle(2.5pt);
   \filldraw[black] (v8) circle(2.5pt);
    \filldraw[black] (v9) circle(2.5pt);
   \filldraw[black] (v10) circle(2.5pt);
    \filldraw[black] (v11) circle(2.5pt);
    \filldraw[black] (v12) circle(2.5pt);
    \filldraw[black] (v13) circle(2.5pt);
    \filldraw[black] (v14) circle(2.5pt);
    \filldraw[black] (v15) circle(2.5pt);
    \filldraw[black] (v16) circle(2.5pt);
    \filldraw[black] (v17) circle(2.5pt);
    \filldraw[black] (v18) circle(2.5pt);
    \filldraw[black] (v19) circle(2.5pt);

    \filldraw[red] (v18) circle(1.5pt);
    \filldraw[red] (v19) circle(1.5pt);
    \filldraw[red] (v0) circle(1.5pt);
    \filldraw[red] (v1) circle(1.5pt);
    \filldraw[red] (v2) circle(1.5pt);
    \filldraw[red] (v3) circle(1.5pt);

    \filldraw[white] (v4) circle(1.5pt);
    \filldraw[white] (v5) circle(1.5pt);
    \filldraw[white] (v9) circle(1.5pt);
    \filldraw[white] (v7) circle(1.5pt);
    \filldraw[white] (v8) circle(1.5pt);

    \filldraw[white] (v10) circle(1.5pt);
    \filldraw[white] (v11) circle(1.5pt);
    \filldraw[white] (v12) circle(1.5pt);

    \filldraw[white] (v14) circle(1.5pt);
    \filldraw[white] (v15) circle(1.5pt);
    \filldraw[white] (v16) circle(1.5pt);
    \filldraw[white] (v17) circle(1.5pt);

    \filldraw[white] (v4) circle(1.5pt);
    \filldraw[yellow] (v5) circle(1.5pt);
    \filldraw[orange] (v6) circle(1.5pt);
    \filldraw[yellow] (v7) circle(1.5pt);
    \filldraw[white] (v8) circle(1.5pt);

    \filldraw[green] (v11) circle(1.5pt);
    \filldraw[green] (v12) circle(1.5pt);
    \filldraw[green] (v13) circle(1.5pt);
    \filldraw[white] (v14) circle(1.5pt);
    \filldraw[white] (v15) circle(1.5pt);

\end{tikzpicture}
}

\centering {(c) $F_m={\cal J}_{20,5}(^*_2)-\{v_5, v_6, v_7\}= {\cal J}_{17,5}(^*_2)\/$}

\end{minipage}
\end{center}
\caption{Deleting any set of $g-2\/$ vertices, containing neither a vertex of an $\alpha$-pair nor neighbor of a vertex of an $\alpha$-pair, form ${\cal J}_{20,5}(^*_2)\/$ results in ${\cal J}_{17,5}(^*_2)\/$}     \label{alpha-pair of vertices  in Fm gives type-2}
\end{figure}

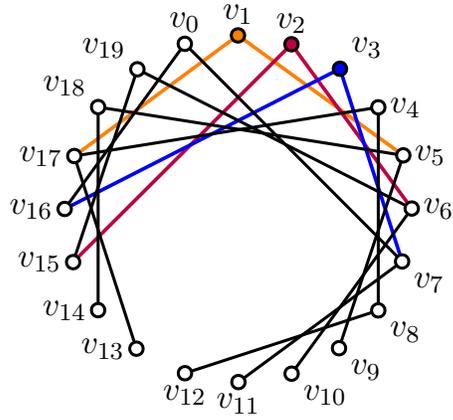
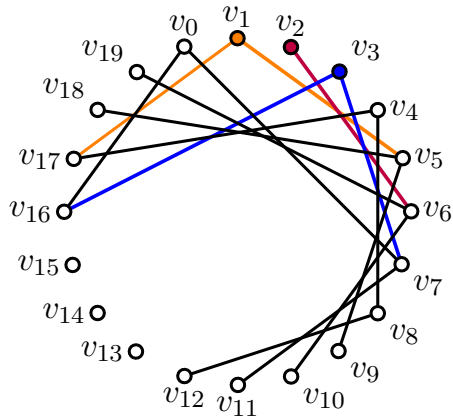
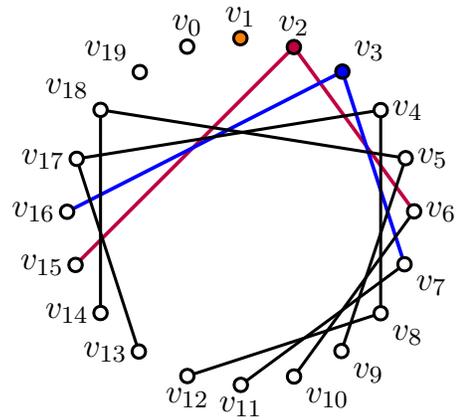
\begin{figure}[htp]
\begin{center}
\begin{minipage}{.6\textwidth}
\begin{center}
\resizebox{6.5cm}{!}{
\begin{tikzpicture}[rotate=-36,.style={draw}]
\coordinate (center) at (0,0);
   \def\radius{2cm}
   \foreach \x in {0, 18,...,360} {
             \filldraw[] (\x:2cm) circle(1pt);
               }


      \coordinate (v0) at (-1.62, 1.18);\filldraw[blue] (v0) circle(1pt);\node[above ] at (v0) {${v_0}$};

       \coordinate (v1) at (-1.18, 1.62);\filldraw[blue] (v1) circle(1pt); \node[above] at (v1) {${v_1}$};
       \coordinate (v2) at (-.62, 1.9);\filldraw[blue] (v2) circle(1pt);  \node[above] at (v2) {${v_2}$};

       \coordinate (v3) at (0, 2);\filldraw[blue] (v3) circle(1pt);\node[above right] at (v3) {${v_{3}}$};
       \coordinate (v4) at (.62, 1.9);\filldraw[blue] (v4) circle(1pt);\node[ right] at (v4) {${v_4}$};

       \coordinate (v5) at (1.18, 1.62);\filldraw[blue] (v5) circle(1pt);\node[ right] at (v5) {${v_5}$};
       \coordinate (v6) at (1.62, 1.18);\filldraw[blue] (v6) circle(1pt);\node[right] at (v6) {${v_{6}}$};

       \coordinate (v7) at (1.89, .62);\filldraw[blue] (v7) circle(1pt);\node[below right] at (v7) {${v_7}$};
      \coordinate (v8) at (2, 0);\filldraw[blue] (v8) circle(1pt);\node[below right] at (v8) {${v_8}$};

       \coordinate (v9) at (1.89, -.62);\filldraw[black] (v9) circle(1pt);\node[below right] at (v9) {${v_{9}}$};
       \coordinate (v10) at (1.62, -1.18);\filldraw[black] (v10) circle(1pt);\node[below right] at (v10) {${v_{10}}$};

       \coordinate (v11) at (1.18, -1.62);\filldraw[black] (v11) circle(1pt);\node[below ] at (v11) {${v_{11}}$};

       \coordinate (v12) at (.62, -1.9);\filldraw[black] (v12) circle(1pt);\node[below ] at (v12) {${v_{12}}$};

       \coordinate (v13) at (0, -2);\filldraw[black] (v13) circle(1pt);\node[left] at (v13) {${v_{13}}$};

       \coordinate (v14) at (-.62, -1.9);\filldraw[black] (v14) circle(1pt);\node[ left] at (v14) {${v_{14}}$};

       \coordinate (v15) at (-1.18, -1.62);\filldraw[black] (v15) circle(1pt);\node[ left] at (v15) {${v_{15}}$};

       \coordinate (v16) at (-1.62, -1.18);\filldraw[black] (v16) circle(1pt);\node[ left] at (v16) {${v_{16}}$};

       \coordinate (v17) at (-1.89, -.62);\filldraw[black] (v17) circle(1pt);\node[left] at (v17) {${v_{17}}$};

       \coordinate (v18) at (-2, 0);\filldraw[black] (v18) circle(1pt);\node[above left] at (v18) {${v_{18}}$};

       \coordinate (v19) at (-1.89, .62);\filldraw[black] (v19) circle(1pt);\node[above left] at (v19) {${v_{19}}$};

     \draw [line width=1.2,orange](v1)-- (v17);
     \draw [line width=1.2,orange](v1)-- (v5);
     \draw [line width=1.2,purple](v15)-- (v2)-- (v6);
     \draw [line width=1.2,blue](v16)-- (v3)-- (v7);

     \draw [line width=1,black](v13)-- (v17)-- (v4)-- (v8)-- (v12);
     \draw [line width=1,black](v14)-- (v18)-- (v5)-- (v9);
     \draw [line width=1,black](v15)-- (v19)-- (v6)-- (v10);
     \draw [line width=1,black](v16)-- (v0)-- (v7)-- (v11);

    \filldraw[black] (v0) circle(2.5pt);
    \filldraw[black] (v1) circle(2.5pt);
    \filldraw[black] (v2) circle(2.5pt);
   \filldraw[black] (v3) circle(2.5pt);
    \filldraw[black] (v4) circle(2.5pt);
   \filldraw[black] (v5) circle(2.5pt);
    \filldraw[black] (v6) circle(2.5pt);
    \filldraw[black] (v7) circle(2.5pt);
   \filldraw[black] (v8) circle(2.5pt);
    \filldraw[black] (v9) circle(2.5pt);
   \filldraw[black] (v10) circle(2.5pt);
    \filldraw[black] (v11) circle(2.5pt);
    \filldraw[black] (v12) circle(2.5pt);
    \filldraw[black] (v13) circle(2.5pt);
    \filldraw[black] (v14) circle(2.5pt);
    \filldraw[black] (v15) circle(2.5pt);
    \filldraw[black] (v16) circle(2.5pt);
    \filldraw[black] (v17) circle(2.5pt);
    \filldraw[black] (v18) circle(2.5pt);
    \filldraw[black] (v19) circle(2.5pt);

    \filldraw[orange] (v1) circle(1.5pt);
    \filldraw[purple] (v2) circle(1.5pt);
    \filldraw[blue] (v3) circle(1.5pt);

    \filldraw[white] (v4) circle(1.5pt);
    \filldraw[white] (v5) circle(1.5pt);
    \filldraw[white] (v6) circle(1.5pt);
    \filldraw[white] (v7) circle(1.5pt);
    \filldraw[white] (v8) circle(1.5pt);
    \filldraw[white] (v9) circle(1.5pt);

    \filldraw[white] (v10) circle(1.5pt);
    \filldraw[white] (v11) circle(1.5pt);
    \filldraw[white] (v12) circle(1.5pt);
    \filldraw[white] (v13) circle(1.5pt);

    \filldraw[white] (v14) circle(1.5pt);
    \filldraw[white] (v15) circle(1.5pt);
    \filldraw[white] (v16) circle(1.5pt);
    \filldraw[white] (v17) circle(1.5pt);
    \filldraw[white] (v18) circle(1.5pt);
    \filldraw[white] (v19) circle(1.5pt);
    \filldraw[white] (v0) circle(1.5pt);

\end{tikzpicture}
}
\end{center}

\centering{(a) ${\cal J}_{20,5}(^*_3)\/$}

\end{minipage}

\vspace{10mm}

\begin{minipage}{.4\textwidth}
\resizebox{6.5cm}{!}{
\begin{tikzpicture}[rotate=-36,.style={draw}]
\coordinate (center) at (0,0);
   \def\radius{2cm}
   \foreach \x in {0, 18,...,360} {
             \filldraw[] (\x:2cm) circle(1pt);
               }


      \coordinate (v0) at (-1.62, 1.18);\filldraw[blue] (v0) circle(1pt);\node[above ] at (v0) {${v_0}$};

       \coordinate (v1) at (-1.18, 1.62);\filldraw[blue] (v1) circle(1pt); \node[above] at (v1) {${v_1}$};
       \coordinate (v2) at (-.62, 1.9);\filldraw[blue] (v2) circle(1pt);  \node[above] at (v2) {${v_2}$};

       \coordinate (v3) at (0, 2);\filldraw[blue] (v3) circle(1pt);\node[above right] at (v3) {${v_{3}}$};
       \coordinate (v4) at (.62, 1.9);\filldraw[blue] (v4) circle(1pt);\node[ right] at (v4) {${v_4}$};

       \coordinate (v5) at (1.18, 1.62);\filldraw[blue] (v5) circle(1pt);\node[ right] at (v5) {${v_5}$};
       \coordinate (v6) at (1.62, 1.18);\filldraw[blue] (v6) circle(1pt);\node[right] at (v6) {${v_{6}}$};

       \coordinate (v7) at (1.89, .62);\filldraw[blue] (v7) circle(1pt);\node[below right] at (v7) {${v_7}$};
      \coordinate (v8) at (2, 0);\filldraw[blue] (v8) circle(1pt);\node[below right] at (v8) {${v_8}$};

       \coordinate (v9) at (1.89, -.62);\filldraw[black] (v9) circle(1pt);\node[below right] at (v9) {${v_{9}}$};
       \coordinate (v10) at (1.62, -1.18);\filldraw[black] (v10) circle(1pt);\node[below right] at (v10) {${v_{10}}$};

       \coordinate (v11) at (1.18, -1.62);\filldraw[black] (v11) circle(1pt);\node[below ] at (v11) {${v_{11}}$};

       \coordinate (v12) at (.62, -1.9);\filldraw[black] (v12) circle(1pt);\node[below ] at (v12) {${v_{12}}$};

       \coordinate (v13) at (0, -2);\filldraw[black] (v13) circle(1pt);\node[left] at (v13) {${v_{13}}$};

       \coordinate (v14) at (-.62, -1.9);\filldraw[black] (v14) circle(1pt);\node[ left] at (v14) {${v_{14}}$};

       \coordinate (v15) at (-1.18, -1.62);\filldraw[black] (v15) circle(1pt);\node[ left] at (v15) {${v_{15}}$};

       \coordinate (v16) at (-1.62, -1.18);\filldraw[black] (v16) circle(1pt);\node[ left] at (v16) {${v_{16}}$};

       \coordinate (v17) at (-1.89, -.62);\filldraw[black] (v17) circle(1pt);\node[left] at (v17) {${v_{17}}$};

       \coordinate (v18) at (-2, 0);\filldraw[black] (v18) circle(1pt);\node[above  left] at (v18) {${v_{18}}$};

       \coordinate (v19) at (-1.89, .62);\filldraw[black] (v19) circle(1pt);\node[above left] at (v19) {${v_{19}}$};

     \draw [line width=1.,orange](v1)-- (v17);
     \draw [line width=1.2,orange](v1)-- (v5);
     \draw [line width=1.2,purple] (v2)-- (v6);
     \draw [line width=1.2,blue](v16)-- (v3)-- (v7);

     \draw [line width=1,black](v17)-- (v4)-- (v8)-- (v12);
     \draw [line width=1,black](v18)-- (v5)-- (v9);
     \draw [line width=1,black] (v19)-- (v6)-- (v10);
     \draw [line width=1,black](v16)--(v0)-- (v7)-- (v11);

    \filldraw[black] (v0) circle(2.5pt);
    \filldraw[black] (v1) circle(2.5pt);
    \filldraw[black] (v2) circle(2.5pt);
   \filldraw[black] (v3) circle(2.5pt);
    \filldraw[black] (v4) circle(2.5pt);
   \filldraw[black] (v5) circle(2.5pt);
    \filldraw[black] (v6) circle(2.5pt);
    \filldraw[black] (v7) circle(2.5pt);
   \filldraw[black] (v8) circle(2.5pt);
    \filldraw[black] (v9) circle(2.5pt);
   \filldraw[black] (v10) circle(2.5pt);
    \filldraw[black] (v11) circle(2.5pt);
    \filldraw[black] (v12) circle(2.5pt);
    \filldraw[black] (v13) circle(2.5pt);
    \filldraw[black] (v14) circle(2.5pt);
    \filldraw[black] (v15) circle(2.5pt);
    \filldraw[black] (v16) circle(2.5pt);
    \filldraw[black] (v17) circle(2.5pt);
    \filldraw[black] (v18) circle(2.5pt);
    \filldraw[black] (v19) circle(2.5pt);

    \filldraw[orange] (v1) circle(1.5pt);
    \filldraw[purple] (v2) circle(1.5pt);
    \filldraw[blue] (v3) circle(1.5pt);

    \filldraw[white] (v4) circle(1.5pt);
    \filldraw[white] (v5) circle(1.5pt);
    \filldraw[white] (v6) circle(1.5pt);
    \filldraw[white] (v7) circle(1.5pt);
    \filldraw[white] (v8) circle(1.5pt);
    \filldraw[white] (v9) circle(1.5pt);

    \filldraw[white] (v10) circle(1.5pt);
    \filldraw[white] (v11) circle(1.5pt);
    \filldraw[white] (v12) circle(1.5pt);
    \filldraw[white] (v13) circle(1.5pt);

    \filldraw[white] (v14) circle(1.5pt);
    \filldraw[white] (v15) circle(1.5pt);
    \filldraw[white] (v16) circle(1.5pt);
    \filldraw[white] (v17) circle(1.5pt);
    \filldraw[white] (v18) circle(1.5pt);
    \filldraw[white] (v19) circle(1.5pt);
    \filldraw[white] (v0) circle(1.5pt);

\end{tikzpicture}
}

\centering{(b) ${\cal J}_{17,5}(^*_2)\/$}

\end{minipage}
\hspace{8mm}
\begin{minipage}{.4\textwidth}
\resizebox{6.5cm}{!}{
\begin{tikzpicture}[rotate=-36,.style={draw}]
\coordinate (center) at (0,0);
   \def\radius{2cm}
   \foreach \x in {0, 18,...,360} {
             \filldraw[] (\x:2cm) circle(1pt);
               }


      \coordinate (v0) at (-1.62, 1.18);\filldraw[blue] (v0) circle(1pt);\node[above ] at (v0) {${v_0}$};

       \coordinate (v1) at (-1.18, 1.62);\filldraw[blue] (v1) circle(1pt); \node[above] at (v1) {${v_1}$};
       \coordinate (v2) at (-.62, 1.9);\filldraw[blue] (v2) circle(1pt);  \node[above] at (v2) {${v_2}$};

       \coordinate (v3) at (0, 2);\filldraw[blue] (v3) circle(1pt);\node[above right] at (v3) {${v_{3}}$};
       \coordinate (v4) at (.62, 1.9);\filldraw[blue] (v4) circle(1pt);\node[ right] at (v4) {${v_4}$};

       \coordinate (v5) at (1.18, 1.62);\filldraw[blue] (v5) circle(1pt);\node[ right] at (v5) {${v_5}$};
       \coordinate (v6) at (1.62, 1.18);\filldraw[blue] (v6) circle(1pt);\node[right] at (v6) {${v_{6}}$};

       \coordinate (v7) at (1.89, .62);\filldraw[blue] (v7) circle(1pt);\node[below right] at (v7) {${v_7}$};
      \coordinate (v8) at (2, 0);\filldraw[blue] (v8) circle(1pt);\node[below right] at (v8) {${v_8}$};

       \coordinate (v9) at (1.89, -.62);\filldraw[black] (v9) circle(1pt);\node[below right] at (v9) {${v_{9}}$};
       \coordinate (v10) at (1.62, -1.18);\filldraw[black] (v10) circle(1pt);\node[below right] at (v10) {${v_{10}}$};

       \coordinate (v11) at (1.18, -1.62);\filldraw[black] (v11) circle(1pt);\node[below ] at (v11) {${v_{11}}$};

       \coordinate (v12) at (.62, -1.9);\filldraw[black] (v12) circle(1pt);\node[below ] at (v12) {${v_{12}}$};

       \coordinate (v13) at (0, -2);\filldraw[black] (v13) circle(1pt);\node[left] at (v13) {${v_{13}}$};

       \coordinate (v14) at (-.62, -1.9);\filldraw[black] (v14) circle(1pt);\node[ left] at (v14) {${v_{14}}$};

       \coordinate (v15) at (-1.18, -1.62);\filldraw[black] (v15) circle(1pt);\node[ left] at (v15) {${v_{15}}$};

       \coordinate (v16) at (-1.62, -1.18);\filldraw[black] (v16) circle(1pt);\node[ left] at (v16) {${v_{16}}$};

       \coordinate (v17) at (-1.89, -.62);\filldraw[black] (v17) circle(1pt);\node[left] at (v17) {${v_{17}}$};

       \coordinate (v18) at (-2, 0);\filldraw[black] (v18) circle(1pt);\node[above left] at (v18) {${v_{18}}$};

       \coordinate (v19) at (-1.89, .62);\filldraw[black] (v19) circle(1pt);\node[above left] at (v19) {${v_{19}}$};

     \draw [line width=1.2,purple](v15)-- (v2)-- (v6);
     \draw [line width=1.2,blue](v16)-- (v3)-- (v7);

     \draw [line width=1,black](v13)-- (v17)-- (v4)-- (v8)-- (v12);
     \draw [line width=1,black](v14)-- (v18)-- (v5)-- (v9);
     \draw [line width=1,black] (v6)-- (v10);
     \draw [line width=1,black](v7)-- (v11);

    \filldraw[black] (v0) circle(2.5pt);
    \filldraw[black] (v1) circle(2.5pt);
    \filldraw[black] (v2) circle(2.5pt);
   \filldraw[black] (v3) circle(2.5pt);
    \filldraw[black] (v4) circle(2.5pt);
   \filldraw[black] (v5) circle(2.5pt);
    \filldraw[black] (v6) circle(2.5pt);
    \filldraw[black] (v7) circle(2.5pt);
   \filldraw[black] (v8) circle(2.5pt);
    \filldraw[black] (v9) circle(2.5pt);
   \filldraw[black] (v10) circle(2.5pt);
    \filldraw[black] (v11) circle(2.5pt);
    \filldraw[black] (v12) circle(2.5pt);
    \filldraw[black] (v13) circle(2.5pt);
    \filldraw[black] (v14) circle(2.5pt);
    \filldraw[black] (v15) circle(2.5pt);
    \filldraw[black] (v16) circle(2.5pt);
    \filldraw[black] (v17) circle(2.5pt);
    \filldraw[black] (v18) circle(2.5pt);
    \filldraw[black] (v19) circle(2.5pt);

    \filldraw[orange] (v1) circle(1.5pt);
    \filldraw[purple] (v2) circle(1.5pt);
    \filldraw[blue] (v3) circle(1.5pt);

    \filldraw[white] (v4) circle(1.5pt);
    \filldraw[white] (v5) circle(1.5pt);
    \filldraw[white] (v6) circle(1.5pt);
    \filldraw[white] (v7) circle(1.5pt);
    \filldraw[white] (v8) circle(1.5pt);
    \filldraw[white] (v9) circle(1.5pt);

    \filldraw[white] (v10) circle(1.5pt);
    \filldraw[white] (v11) circle(1.5pt);
    \filldraw[white] (v12) circle(1.5pt);
    \filldraw[white] (v13) circle(1.5pt);

    \filldraw[white] (v14) circle(1.5pt);
    \filldraw[white] (v15) circle(1.5pt);
    \filldraw[white] (v16) circle(1.5pt);
    \filldraw[white] (v17) circle(1.5pt);
    \filldraw[white] (v18) circle(1.5pt);
    \filldraw[white] (v19) circle(1.5pt);
    \filldraw[white] (v0) circle(1.5pt);

\end{tikzpicture}
}

\centering{(c) ${\cal F}_{17,5}(^*)\/$}

\end{minipage}
\end{center}
\caption{${\cal J}_{20,5}(^*_3)\/$}   \label{example type 3}
\end{figure}

\vspace{5mm} An example of a convex geometric graph ${\cal J}_{n,g}(^*_3)\/$  is depicted in Figure \ref{example type 3}.
 In Figure \ref{example type 3} (a), $\{v_{15}, v_{19}\}$ and $\{v_{16}, v_{0}\}$ are two $\alpha\/$-pairs in  ${\cal J}_{20,5}(^*_3)\/$,
 in (b) $\{v_{16}, v_{0}\}$ is the $\alpha\/$-pair in the convex geometric graph ${\cal J}_{17,5}(^*_2)\/$ which obtained from  ${\cal J}_{20,5}(^*_3)\/$ by deleting $\{v_{13}, v_{14}, v_{15}\}\/$,
 in (c) ${\cal F}_{17,5}(^*)\/$ is obtained from  ${\cal J}_{20,5}(^*_3)\/$ by deleting $\{v_{0}, v_{1}, v_{19}\}\/$.

\begin{figure}[htp]
\begin{center}

\begin{minipage}{.3\textwidth}            

\resizebox{4.5cm}{!}{
\begin{tikzpicture}
\coordinate (center) at (0,0);
   \def\radius{2.2cm}
   \foreach \x in {0, 60,...,360} {
               }

\draw [line width=1.5, white] plot [smooth,  tension=.9] coordinates {(0,4) (-2,3) (-2,0) (0,-1)};

       \coordinate (v1) at (-1.2, 0);\filldraw[black] (v1) circle(4pt);    \node[left] at (v1)    {\LARGE${v_2}$};
       \coordinate (v2) at (1.2,0);\filldraw[black] (v2) circle(4pt);      \node[right] at (v2)   {\LARGE${v_{6}}$};
       \coordinate (v3) at (0.2, 1);\filldraw[black] (v3) circle(4pt);    \node[above] at (0.3, 1)    {\LARGE${v_8}$};
       \coordinate (v4) at (-0.2, 2);\filldraw[black] (v4) circle(4pt);     \node[below] at (-0.3, 2)   {\LARGE${v_5}$};
       \coordinate (v5) at (-1.2, 3);\filldraw[black] (v5) circle(4pt);    \node[left] at (v5)    {\LARGE${v_7}$};
       \coordinate (v6) at (1.2, 3);\filldraw[black] (v6) circle(4pt);     \node[below right] at (v6)   {\LARGE${v_3}$};

      \coordinate (v9) at (0.2, 4);  \filldraw[black] (v9) circle(4pt);     \node[left] at (0.1, 4)     {\LARGE${v_1}$};
      \coordinate (v10) at (-0.2,-1);    \filldraw[black] (v10) circle(4pt);     \node[right] at (v10) {\LARGE${v_{4}}$};

      \coordinate (v11) at (0.8, 4.);    \filldraw[black] (v11) circle(3pt);     \node[right] at (0.9, 4.) {\LARGE${v_{x}}$};
      \coordinate (v12) at (1.8, 3.);    \filldraw[black] (v12) circle(3pt);     \node[right] at (v12) {\LARGE${v_{y}}$};

     \draw [line width=1.5,black](v10) -- (v1) -- (v5) -- (v9);
      \draw [line width=1.5,black] (v6) -- (v2) -- (v10);
   \draw [line width=1.5,black] (v6) -- (v4)-- (v5);
     \draw [line width=1.5,black](v1) -- (v3) -- (v2);
     \draw [line width=1.5,black](v10) -- (v3) ;
     \draw [line width=1.5,black](v4) -- (v9);

     \draw [dashed,black](v6) -- (v9);

     \draw [line width=1.5,black](v9) -- (v11) to [bend right] (v12) to [bend right] (v11);
     \draw [line width=1.5,black](v12) -- (v6);

       \filldraw[white] (-1.2, 0) circle(2pt);
       \filldraw[white] (1.2,0) circle(2pt);
       \filldraw[white] (0.2, 1) circle(2pt);
       \filldraw[white] (-0.2, 2) circle(2pt);
       \filldraw[white] (-1.2, 3) circle(2pt);
       \filldraw[white] (1.2, 3) circle(2pt);

       \filldraw[white] (0.2, 4) circle(2pt);
       \filldraw[white] (-0.2,-1) circle(2pt);

\end{tikzpicture}
}

\centering   (a) $v_1 \in \{x_1,  y_1\}\/$

\end{minipage}
\begin{minipage}{.3\textwidth}            

\resizebox{4.5cm}{!}{
\begin{tikzpicture}
\coordinate (center) at (0,0);
   \def\radius{2.2cm}
   \foreach \x in {0, 60,...,360} {
               }

\draw [line width=1.5, white] plot [smooth,  tension=.9] coordinates {(0,4) (-2,3) (-2,0) (0,-1)};

       \coordinate (v1) at (-1.2, 0);\filldraw[black] (v1) circle(4pt);    \node[left] at (v1)    {\LARGE${v_2}$};
       \coordinate (v2) at (1.2,0);\filldraw[black] (v2) circle(4pt);      \node[below right] at (v2)   {\LARGE${v_{6}}$};
       \coordinate (v3) at (0.2, 1);\filldraw[black] (v3) circle(4pt);    \node[above] at (0.3, 1)    {\LARGE${v_8}$};
       \coordinate (v4) at (-0.2, 2);\filldraw[black] (v4) circle(4pt);     \node[below] at (-0.3, 2)   {\LARGE${v_5}$};
       \coordinate (v5) at (-1.2, 3);\filldraw[black] (v5) circle(4pt);    \node[left] at (v5)    {\LARGE${v_7}$};
       \coordinate (v6) at (1.2, 3);\filldraw[black] (v6) circle(4pt);     \node[right] at (v6)   {\LARGE${v_3}$};

      \coordinate (v9) at (0.2, 4);  \filldraw[black] (v9) circle(4pt);     \node[right] at (v9)     {\LARGE${v_1}$};
      \coordinate (v10) at (-0.2,-1);    \filldraw[black] (v10) circle(4pt);     \node[right] at (v10) {\LARGE${v_{4}}$};

      \coordinate (v11) at (0.8, 1.);    \filldraw[black] (v11) circle(3pt);     \node[above] at (v11) {\LARGE${v_{x}}$};
      \coordinate (v12) at (1.8, 0.);    \filldraw[black] (v12) circle(3pt);     \node[right] at (v12) {\LARGE${v_{y}}$};

     \draw [line width=1.5,black](v10) -- (v1) -- (v5) -- (v9) -- (v6) -- (v2) -- (v10);
   \draw [line width=1.5,black] (v6) -- (v4)-- (v5);
     \draw [line width=1.5,black](v1) -- (v3);
     \draw [dashed,black](v2) -- (v3);
     \draw [line width=1.5,black](v10) -- (v3) ;
     \draw [line width=1.5,black](v4) -- (v9);

     \draw [line width=1.5,black](v3) -- (v11) to [bend right] (v12) to [bend right] (v11);
     \draw [line width=1.5,black](v12) -- (v2);

       \filldraw[white] (-1.2, 0) circle(2pt);
       \filldraw[white] (1.2,0) circle(2pt);
       \filldraw[white] (0.2, 1) circle(2pt);
       \filldraw[white] (-0.2, 2) circle(2pt);
       \filldraw[white] (-1.2, 3) circle(2pt);
       \filldraw[white] (1.2, 3) circle(2pt);

       \filldraw[white] (0.2, 4) circle(2pt);
       \filldraw[white] (-0.2,-1) circle(2pt);

\end{tikzpicture}
}

\centering  (b1) $v_{n-2} \in \{x_1,  y_1\}\/$,       $v_{n-4}\notin N_{G_e}(v_{n-2})-y_1\/$

\end{minipage}
\begin{minipage}{.3\textwidth}            

\resizebox{4.33cm}{!}{
\begin{tikzpicture}
\coordinate (center) at (0,0);
   \def\radius{2.2cm}
   \foreach \x in {0, 60,...,360} {
               }

\draw [line width=1.5, white] plot [smooth,  tension=.9] coordinates {(0,4) (-2,3) (-2,0) (0,-1)};

       \coordinate (v1) at (-1.2, 0);\filldraw[black] (v1) circle(4pt);    \node[left] at (v1)    {\LARGE${v_2}$};
       \coordinate (v2) at (1.2,0);\filldraw[black] (v2) circle(4pt);      \node[right] at (v2)   {\LARGE${v_{6}}$};
       \coordinate (v3) at (0.2, 1);\filldraw[black] (v3) circle(4pt);    \node[above] at (0.3, 1)    {\LARGE${v_8}$};
       \coordinate (v4) at (-0.2, 2);\filldraw[black] (v4) circle(4pt);     \node[below] at (-0.3, 2)   {\LARGE${v_5}$};
       \coordinate (v5) at (-1.2, 3);\filldraw[black] (v5) circle(4pt);    \node[left] at (v5)    {\LARGE${v_7}$};
       \coordinate (v6) at (1.2, 3);\filldraw[black] (v6) circle(4pt);     \node[right] at (v6)   {\LARGE${v_3}$};

      \coordinate (v9) at (0.2, 4);  \filldraw[black] (v9) circle(4pt);     \node[right] at (v9)     {\LARGE${v_1}$};
      \coordinate (v10) at (-0.2,-1);    \filldraw[black] (v10) circle(4pt);     \node[left] at (v10) {\LARGE${v_{4}}$};

      \coordinate (v11) at (0.8, 1.);    \filldraw[black] (v11) circle(3pt);     \node[above] at (v11) {\LARGE${v_{x}}$};
      \coordinate (v12) at (.4, -1);    \filldraw[black] (v12) circle(3pt);     \node[right] at (v12) {\LARGE${v_{y}}$};

     \draw [line width=1.5,black](v10) -- (v1) -- (v5) -- (v9) -- (v6) -- (v2) -- (v10);
   \draw [line width=1.5,black] (v6) -- (v4)-- (v5);
     \draw [line width=1.5,black](v1) -- (v3) -- (v2);
     \draw [dashed,black](v10) -- (v3) ;
     \draw [line width=1.5,black](v4) -- (v9);

     \draw [line width=1.5,black](v3) -- (v11) to [bend right] (v12) to [bend right] (v11);
     \draw [line width=1.5,black](v12) -- (v10);

       \filldraw[white] (-1.2, 0) circle(2pt);
       \filldraw[white] (1.2,0) circle(2pt);
       \filldraw[white] (0.2, 1) circle(2pt);
       \filldraw[white] (-0.2, 2) circle(2pt);
       \filldraw[white] (-1.2, 3) circle(2pt);
       \filldraw[white] (1.2, 3) circle(2pt);

       \filldraw[white] (0.2, 4) circle(2pt);
       \filldraw[white] (-0.2,-1) circle(2pt);

\end{tikzpicture}
}

\centering   (b1)  $v_{n-2} \in \{x_1,  y_1\}\/$,      $v_{n-4}\in N_{G_e}(v_{n-2})-y_1\/$

\end{minipage}

\begin{minipage}{.3\textwidth}              
\resizebox{4.5cm}{!}{
\begin{tikzpicture}
\coordinate (center) at (0,0);
   \def\radius{2.2cm}
   \foreach \x in {0, 60,...,360} {
               }


       \coordinate (v1) at (-1, 0);\filldraw[black] (v1) circle(4pt);    \node[left] at (v1)    {\LARGE${v_4}$};
       \coordinate (v2) at (1,0);\filldraw[black] (v2) circle(4pt);      \node[right] at (v2)   {\LARGE${v_{6}}$};
       \coordinate (v3) at (-1, 1.5);\filldraw[black] (v3) circle(4pt);    \node[left] at (v3)    {\LARGE${v_8}$};
       \coordinate (v4) at (1, 1.5);\filldraw[black] (v4) circle(4pt);     \node[right] at (v4)   {\LARGE${v_3}$};
       \coordinate (v5) at (-1, 3);\filldraw[black] (v5) circle(4pt);    \node[left] at (v5)    {\LARGE${v_5}$};
       \coordinate (v6) at (1, 3);\filldraw[black] (v6) circle(4pt);     \node[right] at (v6)   {\LARGE${v_1}$};

      \coordinate (v9) at (0, 4);  \filldraw[black] (v9) circle(4pt);     \node[right] at (v9)     {\LARGE${v_7}$};
      \coordinate (v10) at (0,-1);    \filldraw[black] (v10) circle(4pt);     \node[right] at (v10) {\LARGE${v_{2}}$};

      \coordinate (v11) at (-0.5, .75);    \filldraw[black] (v11) circle(3pt);     \node[above] at (v11) {\LARGE${v_{y}}$};
      \coordinate (v12) at (0.5, .75);    \filldraw[black] (v12) circle(3pt);     \node[above] at (v12) {\LARGE${v_{x}}$};

     \draw [line width=1.5,black] (v1) -- (v3) -- (v5) -- (v9) -- (v6) -- (v4) -- (v2) -- (v10);
     \draw [line width=1.5,black](v3) -- (v4);
     \draw [line width=1.5,black](v5) -- (v6);
     \draw [line width=1.5,black](v1) -- (v2);

\draw [line width=1.5, black] plot [smooth,  tension=.9] coordinates {(0,4) (-2,3) (-2,0) (0,-1)};

     \draw [line width=1.5,black](v1) -- (v11) to [bend right] (v12) to [bend right] (v11);
     \draw [line width=1.5,black](v12) -- (v10);

    \draw [dashed,black](v1) -- (v10);

       \filldraw[white] (-1, 0) circle(2pt);
       \filldraw[white] (1, 0) circle(2pt);
       \filldraw[white] (-1, 1.5) circle(2pt);
       \filldraw[white] (1, 1.5) circle(2pt);
       \filldraw[white] (-1, 3) circle(2pt);
       \filldraw[white] (1, 3) circle(2pt);
       \filldraw[white] (1, 4) circle(2pt);

       \filldraw[white] (0, 4) circle(2pt);
       \filldraw[white] (0,-1) circle(2pt);

\end{tikzpicture}
}

\centering  $2\in\{i, j\}\/$ and $v_4 \notin N_{G_e}(x_1)-y_1\/$

\end{minipage}\begin{minipage}{.3\textwidth}              
\resizebox{4.5cm}{!}{
\begin{tikzpicture}
\coordinate (center) at (0,0);
   \def\radius{2.2cm}
   \foreach \x in {0, 60,...,360} {
               }


       \coordinate (v1) at (-1, 0);\filldraw[black] (v1) circle(4pt);    \node[left] at (v1)    {\LARGE${v_4}$};
       \coordinate (v2) at (1,0);\filldraw[black] (v2) circle(4pt);      \node[right] at (v2)   {\LARGE${v_{6}}$};
       \coordinate (v3) at (-1, 1.5);\filldraw[black] (v3) circle(4pt);    \node[left] at (v3)    {\LARGE${v_8}$};
       \coordinate (v4) at (1, 1.5);\filldraw[black] (v4) circle(4pt);     \node[right] at (v4)   {\LARGE${v_3}$};
       \coordinate (v5) at (-1, 3);\filldraw[black] (v5) circle(4pt);    \node[left] at (v5)    {\LARGE${v_5}$};
       \coordinate (v6) at (1, 3);\filldraw[black] (v6) circle(4pt);     \node[right] at (v6)   {\LARGE${v_1}$};

      \coordinate (v9) at (0, 4);  \filldraw[black] (v9) circle(4pt);     \node[right] at (v9)     {\LARGE${v_7}$};
      \coordinate (v10) at (0,-1);    \filldraw[black] (v10) circle(4pt);     \node[right] at (v10) {\LARGE${v_{2}}$};

      \coordinate (v11) at (-0.5, .75);    \filldraw[black] (v11) circle(3pt);     \node[above] at (v11) {\LARGE${v_{x}}$};
      \coordinate (v12) at (0.5, .75);    \filldraw[black] (v12) circle(3pt);     \node[above] at (v12) {\LARGE${v_{y}}$};

     \draw [line width=1.5,black](v10) -- (v1) -- (v3) -- (v5) -- (v9) -- (v6) -- (v4) -- (v2);
     \draw [line width=1.5,black](v3) -- (v4);
     \draw [line width=1.5,black](v5) -- (v6);
     \draw [line width=1.5,black](v1) -- (v2);

\draw [line width=1.5, black] plot [smooth,  tension=.9] coordinates {(0,4) (-2,3) (-2,0) (0,-1)};

     \draw [line width=1.5,black](v10) -- (v11) to [bend right] (v12) to [bend right] (v11);
     \draw [line width=1.5,black](v12) -- (v2);

     \draw [dashed,black](v10) -- (v2);

       \filldraw[white] (-1, 0) circle(2pt);
       \filldraw[white] (1, 0) circle(2pt);
       \filldraw[white] (-1, 1.5) circle(2pt);
       \filldraw[white] (1, 1.5) circle(2pt);
       \filldraw[white] (-1, 3) circle(2pt);
       \filldraw[white] (1, 3) circle(2pt);
       \filldraw[white] (1, 4) circle(2pt);

       \filldraw[white] (0, 4) circle(2pt);
       \filldraw[white] (0,-1) circle(2pt);

\end{tikzpicture}
}

\centering $2\in\{i, j\}\/$ and $v_4 \in N_{G_e}(x_1)-y_1\/$

\end{minipage}\begin{minipage}{.3\textwidth}              
\resizebox{4.5cm}{!}{
\begin{tikzpicture}
\coordinate (center) at (0,0);
   \def\radius{2.2cm}
   \foreach \x in {0, 60,...,360} {
               }


       \coordinate (v1) at (-1, 0);\filldraw[black] (v1) circle(4pt);    \node[left] at (v1)    {\LARGE${v_4}$};
       \coordinate (v2) at (1,0);\filldraw[black] (v2) circle(4pt);      \node[right] at (v2)   {\LARGE${v_{6}}$};
       \coordinate (v3) at (-1, 1.5);\filldraw[black] (v3) circle(4pt);    \node[left] at (v3)    {\LARGE${v_8}$};
       \coordinate (v4) at (1, 1.5);\filldraw[black] (v4) circle(4pt);     \node[right] at (v4)   {\LARGE${v_3}$};
       \coordinate (v5) at (-1, 3);\filldraw[black] (v5) circle(4pt);    \node[left] at (v5)    {\LARGE${v_5}$};
       \coordinate (v6) at (1, 3);\filldraw[black] (v6) circle(4pt);     \node[right] at (v6)   {\LARGE${v_1}$};

      \coordinate (v9) at (0, 4);  \filldraw[black] (v9) circle(4pt);     \node[right] at (v9)     {\LARGE${v_7}$};
      \coordinate (v10) at (0,-1);    \filldraw[black] (v10) circle(4pt);     \node[right] at (v10) {\LARGE${v_{2}}$};

     \coordinate (v11) at (-0.5, .75);    \filldraw[black] (v11) circle(3pt);     \node[above] at (v11) {\LARGE${v_{x}}$};
      \coordinate (v12) at (0.5, .75);    \filldraw[black] (v12) circle(3pt);     \node[above] at (v12) {\LARGE${v_{y}}$};

     \draw [line width=1.5,black](v10) -- (v1) -- (v3) -- (v5) -- (v9) -- (v6) -- (v4) -- (v2) -- (v10);
     \draw [line width=1.5,black](v3) -- (v4);
     \draw [line width=1.5,black](v5) -- (v6);
     \draw [dashed,black](v1) -- (v2);

\draw [line width=1.5, black] plot [smooth,  tension=.9] coordinates {(0,4) (-2,3) (-2,0) (0,-1)};

     \draw [line width=1.5,black](v1) -- (v11) to [bend right] (v12) to [bend right] (v11);
     \draw [line width=1.5,black](v12) -- (v2);

       \filldraw[white] (-1, 0) circle(2pt);
       \filldraw[white] (1, 0) circle(2pt);
       \filldraw[white] (-1, 1.5) circle(2pt);
       \filldraw[white] (1, 1.5) circle(2pt);
       \filldraw[white] (-1, 3) circle(2pt);
       \filldraw[white] (1, 3) circle(2pt);
       \filldraw[white] (1, 4) circle(2pt);

       \filldraw[white] (0, 4) circle(2pt);
       \filldraw[white] (0,-1) circle(2pt);

\end{tikzpicture}
}

\centering  $2\notin\{i, j\}\/$

\end{minipage}

\vspace{3mm}
\centering   (b2)  $v_{n-2} \notin \{x_1,  y_1\}\/$,

\end{center}
\caption{$G\/$ has a $2\/$-cycle and $G_e\/$ has no $2\/$-cycle}    \label{Case(2) of labeling (2.1)}
\end{figure}

\begin{figure}[htp]
\begin{center}                     
\begin{minipage}{.3\textwidth}              
\begin{center}
\resizebox{4cm}{!}{
\begin{tikzpicture}
\coordinate (center) at (0,0);
   \def\radius{2.2cm}
   \foreach \x in {0, 60,...,360} {
               }


       \coordinate (v1) at (-1, 0);\filldraw[black] (v1) circle(4pt);\node[left] at (v1)   {\LARGE${v_4}$};
       \coordinate (v2) at (1,0);\filldraw[black] (v2) circle(4pt);\node[right] at (v2)    {\LARGE${v_{2}}$};
       \coordinate (v3) at (-1, 1);\filldraw[black] (v3) circle(4pt);\node[left] at (v3)   {\LARGE${v_7}$};
       \coordinate (v4) at (1, 1);\filldraw[black] (v4) circle(4pt);\node[right] at (v4)   {\LARGE${v_6}$};
       \coordinate (v5) at (-1, 2);\filldraw[black] (v5) circle(4pt);\node[left] at (v5)   {\LARGE${v_5}$};
       \coordinate (v6) at (1, 2);\filldraw[black] (v6) circle(4pt);\node[right] at (v6)   {\LARGE${v_8}$};
       \coordinate (v7) at (-1, 3);\filldraw[black] (v7) circle(4pt);\node[left] at (v7)   {\LARGE${v_1}$};
       \coordinate (v8) at (1, 3);\filldraw[black] (v8) circle(4pt);\node[right] at (v8)   {\LARGE${v_{3}}$};

      \coordinate (v9) at (0, 3.33);\filldraw[black] (v9) circle(3pt);\node[above] at (v9)   {\LARGE${v_{x}}$};
       \coordinate (v10) at (1, 1.5);\filldraw[black] (v10) circle(3pt);\node[right] at (v10)   {\LARGE${v_{y}}$};

     \draw [line width=1.5,black](v8) -- (v6) -- (v4) -- (v2) --(v1) -- (v3) -- (v5) -- (v7);
     \draw [line width=1.5,black](v2) -- (v3);
     \draw [line width=1.5,black](v1) -- (v4);
     \draw [line width=1.5,black](v6) -- (v5);

     \draw [line width=1.5,black](v7) to[bend right] (v8);
     \draw [line width=1.5,black](v8) to[bend right] (v7);

     \draw [line width=1.5,black](v9) to[bend right] (v10);

       \filldraw[white] (-1, 0) circle(2pt);
       \filldraw[white] (1, 0) circle(2pt);
       \filldraw[white] (-1, 1) circle(2pt);
       \filldraw[white] (1, 1) circle(2pt);
       \filldraw[white] (-1, 2) circle(2pt);
       \filldraw[white] (1, 2) circle(2pt);
       \filldraw[white] (-1, 3) circle(2pt);
       \filldraw[white] (1, 3) circle(2pt);

       \filldraw[white] (0, 5) circle(2pt);
       \filldraw[white] (0, -1) circle(2pt);

\end{tikzpicture}
}
\end{center}
\centering (a)      

\end{minipage}
\begin{minipage}{.3\textwidth}              
\begin{center}
\resizebox{4.cm}{!}{
\begin{tikzpicture}
\coordinate (center) at (0,0);
   \def\radius{2.2cm}
   \foreach \x in {0, 60,...,360} {
               }


       \coordinate (v1) at (-1, 0);\filldraw[black] (v1) circle(4pt);\node[left] at (v1)   {\LARGE${v_4}$};
       \coordinate (v2) at (1,0);\filldraw[black] (v2) circle(4pt);\node[right] at (v2)    {\LARGE${v_{2}}$};
       \coordinate (v3) at (-1, 1);\filldraw[black] (v3) circle(4pt);\node[left] at (v3)   {\LARGE${v_7}$};
       \coordinate (v4) at (1, 1);\filldraw[black] (v4) circle(4pt);\node[right] at (v4)   {\LARGE${v_8}$};
       \coordinate (v5) at (-1, 2);\filldraw[black] (v5) circle(4pt);\node[left] at (v5)   {\LARGE${v_5}$};
       \coordinate (v6) at (1, 2);\filldraw[black] (v6) circle(4pt);\node[right] at (v6)   {\LARGE${v_6}$};
       \coordinate (v7) at (-1, 3);\filldraw[black] (v7) circle(4pt);\node[left] at (v7)   {\LARGE${v_1}$};
       \coordinate (v8) at (1, 3);\filldraw[black] (v8) circle(4pt);\node[right] at (v8)   {\LARGE${v_{3}}$};

      \coordinate (v9) at (0, 3.33);\filldraw[black] (v9) circle(3pt);\node[above] at (v9)   {\LARGE${v_{x}}$};
       \coordinate (v10) at (-.5, 1.73);\filldraw[black] (v10) circle(3pt);\node[below] at (v10)   {\LARGE${v_{y}}$};

     \draw [line width=1.5,black](v8) -- (v6) -- (v4) -- (v2) -- (v1) -- (v3) -- (v5) -- (v7);
     \draw [line width=1.5,black](v1) -- (v6);
     \draw [line width=1.5,black](v2) -- (v3);
     \draw [line width=1.5,black](v4) -- (v5);

     \draw [line width=1.5,black](v7) to[bend right] (v8);
     \draw [line width=1.5,black](v8) to[bend right] (v7);

     \draw [line width=1.5,black](v9) to[bend left] (v10);

       \filldraw[white] (-1, 0) circle(2pt);
       \filldraw[white] (1, 0) circle(2pt);
       \filldraw[white] (-1, 1) circle(2pt);
       \filldraw[white] (1, 1) circle(2pt);
       \filldraw[white] (-1, 2) circle(2pt);
       \filldraw[white] (1, 2) circle(2pt);
       \filldraw[white] (-1, 3) circle(2pt);
       \filldraw[white] (1, 3) circle(2pt);

       \filldraw[white] (0, 5) circle(2pt);
       \filldraw[white] (0, -1) circle(2pt);

\end{tikzpicture}
}
\end{center}
\centering  (a)      

\end{minipage}
\begin{minipage}{.3\textwidth}              
\begin{center}
\resizebox{4cm}{!}{
\begin{tikzpicture}
\coordinate (center) at (0,0);
   \def\radius{2.2cm}
   \foreach \x in {0, 60,...,360} {
               }


       \coordinate (v1) at (-1, 0);\filldraw[black] (v1) circle(4pt);\node[left] at (v1)   {\LARGE${v_4}$};
       \coordinate (v2) at (1,0);\filldraw[black] (v2) circle(4pt);\node[right] at (v2)    {\LARGE${v_{2}}$};
       \coordinate (v3) at (-1, 1);\filldraw[black] (v3) circle(4pt);\node[left] at (v3)   {\LARGE${v_7}$};
       \coordinate (v4) at (1, 1);\filldraw[black] (v4) circle(4pt);\node[right] at (v4)   {\LARGE${v_6}$};
       \coordinate (v5) at (-1, 2);\filldraw[black] (v5) circle(4pt);\node[left] at (v5)   {\LARGE${v_5}$};
       \coordinate (v6) at (1, 2);\filldraw[black] (v6) circle(4pt);\node[right] at (v6)   {\LARGE${v_8}$};
       \coordinate (v7) at (-1, 3);\filldraw[black] (v7) circle(4pt);\node[left] at (v7)   {\LARGE${v_1}$};
       \coordinate (v8) at (1, 3);\filldraw[black] (v8) circle(4pt);\node[right] at (v8)   {\LARGE${v_{3}}$};

      \coordinate (v9) at (0, 3.33);\filldraw[black] (v9) circle(3pt);\node[above] at (v9)   {\LARGE${v_{x}}$};
       \coordinate (v10) at (0, 2);\filldraw[black] (v10) circle(3pt);\node[below] at (v10)   {\LARGE${v_{y}}$};

     \draw [line width=1.5,black](v8) -- (v6) -- (v4) -- (v2) --(v1) -- (v3) -- (v5) -- (v7);
     \draw [line width=1.5,black](v2) -- (v3);
     \draw [line width=1.5,black](v1) -- (v4);
     \draw [line width=1.5,black](v6) -- (v5);

     \draw [line width=1.5,black](v7) to[bend right] (v8);
     \draw [line width=1.5,black](v8) to[bend right] (v7);

     \draw [line width=1.5,black](v9) to[bend left] (v10);

       \filldraw[white] (-1, 0) circle(2pt);
       \filldraw[white] (1, 0) circle(2pt);
       \filldraw[white] (-1, 1) circle(2pt);
       \filldraw[white] (1, 1) circle(2pt);
       \filldraw[white] (-1, 2) circle(2pt);
       \filldraw[white] (1, 2) circle(2pt);
       \filldraw[white] (-1, 3) circle(2pt);
       \filldraw[white] (1, 3) circle(2pt);

       \filldraw[white] (0, 5) circle(2pt);
       \filldraw[white] (0, -1) circle(2pt);

\end{tikzpicture}
}
\end{center}
\centering (a)   

\end{minipage}

\begin{minipage}{.25\textwidth}              
\resizebox{4cm}{!}{
\begin{tikzpicture}
\coordinate (center) at (0,0);
   \def\radius{2.2cm}
   \foreach \x in {0, 60,...,360} {
               }


       \coordinate (v1) at (-1, 0);\filldraw[black] (v1) circle(4pt);\node[left] at (v1)   {\LARGE${v_4}$};
       \coordinate (v2) at (1,0);\filldraw[black] (v2) circle(4pt);\node[right] at (v2)    {\LARGE${v_{7}}$};
       \coordinate (v3) at (-1, 1);\filldraw[black] (v3) circle(4pt);\node[left] at (v3)   {\LARGE${v_8}$};
       \coordinate (v4) at (1, 1);\filldraw[black] (v4) circle(4pt);\node[right] at (v4)   {\LARGE${v_2}$};
       \coordinate (v5) at (-1, 2);\filldraw[black] (v5) circle(4pt);\node[left] at (v5)   {\LARGE${v_5}$};
       \coordinate (v6) at (1, 2);\filldraw[black] (v6) circle(4pt);\node[right] at (v6)   {\LARGE${v_6}$};
       \coordinate (v7) at (-1, 3);\filldraw[black] (v7) circle(4pt);\node[left] at (v7)   {\LARGE${v_1}$};
       \coordinate (v8) at (1, 3);\filldraw[black] (v8) circle(4pt);\node[right] at (v8)   {\LARGE${v_{3}}$};

      \coordinate (v9) at (0, 3.33);\filldraw[black] (v9) circle(3pt);\node[above] at (v9)   {\LARGE${v_{x}}$};
       \coordinate (v10) at (0, 2.73);\filldraw[black] (v10) circle(3pt);\node[below] at (v10)   {\LARGE${v_{y}}$};

     \draw [line width=1.5,black](v8) -- (v6) -- (v4) -- (v2) -- (v1) -- (v3) -- (v5) -- (v7);
     \draw [line width=1.5,black](v1) -- (v6);
     \draw [line width=1.5,black](v4) -- (v3);
     \draw [line width=1.5,black](v2) -- (v5);

     \draw [line width=1.5,black](v7) to[bend right] (v8);
     \draw [line width=1.5,black](v8) to[bend right] (v7);

     \draw [line width=1.5,black](v9) to[bend left] (v10);

       \filldraw[white] (-1, 0) circle(2pt);
       \filldraw[white] (1, 0) circle(2pt);
       \filldraw[white] (-1, 1) circle(2pt);
       \filldraw[white] (1, 1) circle(2pt);
       \filldraw[white] (-1, 2) circle(2pt);
       \filldraw[white] (1, 2) circle(2pt);
       \filldraw[white] (-1, 3) circle(2pt);
       \filldraw[white] (1, 3) circle(2pt);

       \filldraw[white] (0, 5) circle(2pt);
       \filldraw[white] (0, -1) circle(2pt);

\end{tikzpicture}
}

\centering (b1)

\end{minipage}
\begin{minipage}{.25\textwidth}              

\resizebox{4cm}{!}{
\begin{tikzpicture}
\coordinate (center) at (0,0);
   \def\radius{2.2cm}
   \foreach \x in {0, 60,...,360} {
               }


       \coordinate (v1) at (-1, 0);\filldraw[black] (v1) circle(4pt);\node[left] at (v1)   {\LARGE${v_8}$};
       \coordinate (v2) at (1,0);\filldraw[black] (v2) circle(4pt);\node[right] at (v2)    {\LARGE${v_{6}}$};
       \coordinate (v3) at (-1, 1);\filldraw[black] (v3) circle(4pt);\node[left] at (v3)   {\LARGE${v_2}$};
       \coordinate (v4) at (1, 1);\filldraw[black] (v4) circle(4pt);\node[right] at (v4)   {\LARGE${v_4}$};
       \coordinate (v5) at (-1, 2);\filldraw[black] (v5) circle(4pt);\node[left] at (v5)   {\LARGE${v_5}$};
       \coordinate (v6) at (1, 2);\filldraw[black] (v6) circle(4pt);\node[right] at (v6)   {\LARGE${v_7}$};
       \coordinate (v7) at (-1, 3);\filldraw[black] (v7) circle(4pt);\node[left] at (v7)   {\LARGE${v_1}$};
       \coordinate (v8) at (1, 3);\filldraw[black] (v8) circle(4pt);\node[right] at (v8)   {\LARGE${v_{3}}$};

       \coordinate (v9) at (0, 3.33);\filldraw[black] (v9) circle(3pt);\node[above] at (v9)   {\LARGE${v_{x}}$};
       \coordinate (v10) at (-1, 2.5);\filldraw[black] (v10) circle(3pt);\node[left] at (v10)   {\LARGE${v_{y}}$};

     \draw [line width=1.5,black](v8) -- (v6) -- (v4) -- (v2) --(v1) -- (v3) -- (v5) -- (v7);
     \draw [line width=1.5,black](v2) -- (v3);
     \draw [line width=1.5,black](v1) -- (v4);
     \draw [line width=1.5,black](v6) -- (v5);

     \draw [line width=1.5,black](v7) to[bend right] (v8);
     \draw [line width=1.5,black](v8) to[bend right] (v7);

     \draw [line width=1.5,black](v9) to[bend left] (v10);

       \filldraw[white] (-1, 0) circle(2pt);
       \filldraw[white] (1, 0) circle(2pt);
       \filldraw[white] (-1, 1) circle(2pt);
       \filldraw[white] (1, 1) circle(2pt);
       \filldraw[white] (-1, 2) circle(2pt);
       \filldraw[white] (1, 2) circle(2pt);
       \filldraw[white] (-1, 3) circle(2pt);
       \filldraw[white] (1, 3) circle(2pt);

       \filldraw[white] (0, 5) circle(2pt);
       \filldraw[white] (0, -1) circle(2pt);

\end{tikzpicture}
}

\centering (b2.1)

\end{minipage}
\begin{minipage}{.25\textwidth}              

\resizebox{4cm}{!}{
\begin{tikzpicture}
\coordinate (center) at (0,0);
   \def\radius{2.2cm}
   \foreach \x in {0, 60,...,360} {
               }


       \coordinate (v1) at (-1, 0);\filldraw[black] (v1) circle(4pt);\node[left] at (v1)   {\LARGE${v_8}$};
       \coordinate (v2) at (1,0);\filldraw[black] (v2) circle(4pt);\node[right] at (v2)    {\LARGE${v_{6}}$};
       \coordinate (v3) at (-1, 1);\filldraw[black] (v3) circle(4pt);\node[left] at (v3)   {\LARGE${v_2}$};
       \coordinate (v4) at (1, 1);\filldraw[black] (v4) circle(4pt);\node[right] at (v4)   {\LARGE${v_4}$};
       \coordinate (v5) at (-1, 2);\filldraw[black] (v5) circle(4pt);\node[left] at (v5)   {\LARGE${v_5}$};
       \coordinate (v6) at (1, 2);\filldraw[black] (v6) circle(4pt);\node[right] at (v6)   {\LARGE${v_7}$};
       \coordinate (v7) at (-1, 3);\filldraw[black] (v7) circle(4pt);\node[left] at (v7)   {\LARGE${v_1}$};
       \coordinate (v8) at (1, 3);\filldraw[black] (v8) circle(4pt);\node[right] at (v8)   {\LARGE${v_{3}}$};

       \coordinate (v9) at (0, 3.33);\filldraw[black] (v9) circle(3pt);\node[above] at (v9)   {\LARGE${v_{x}}$};
       \coordinate (v10) at (0, 2.47);\filldraw[black] (v10) circle(3pt);\node[below right] at (0.05, 2.7)   {\LARGE${v_{y}}$};

     \draw [line width=1.5,black](v6) -- (v4) -- (v2) --(v1) -- (v3) -- (v6);
     \draw [line width=1.5,black](v5) -- (v6);
     \draw [line width=1.5,black](v5) -- (v7);
     \draw [line width=1.5,black](v5) -- (v8);

     \draw [line width=1.5,black](v1) -- (v4);
     \draw [line width=1.5,black](v2) -- (v3);

     \draw [line width=1.5,black](v7) to[bend right] (v8);
     \draw [line width=1.5,black](v8) to[bend right] (v7);

     \draw [line width=1.5,black](v9) to[bend right] (v10);

       \filldraw[white] (-1, 0) circle(2pt);
       \filldraw[white] (1, 0) circle(2pt);
       \filldraw[white] (-1, 1) circle(2pt);
       \filldraw[white] (1, 1) circle(2pt);
       \filldraw[white] (-1, 2) circle(2pt);
       \filldraw[white] (1, 2) circle(2pt);
       \filldraw[white] (-1, 3) circle(2pt);
       \filldraw[white] (1, 3) circle(2pt);

       \filldraw[white] (0, 5) circle(2pt);
       \filldraw[white] (0, -1) circle(2pt);

\end{tikzpicture}
}

\centering (b2.2) $y_1\neq v_{n-2}$

\end{minipage}
\begin{minipage}{.22\textwidth}              
\resizebox{4cm}{!}{
\begin{tikzpicture}
\coordinate (center) at (0,0);
   \def\radius{2.2cm}
   \foreach \x in {0, 60,...,360} {
               }


       \coordinate (v1) at (-1, 0);\filldraw[black] (v1) circle(4pt);\node[left] at (v1)   {\LARGE${v_6}$};
       \coordinate (v2) at (1,0);\filldraw[black] (v2) circle(4pt);\node[right] at (v2)    {\LARGE${v_{7}}$};
       \coordinate (v3) at (-1, 1);\filldraw[black] (v3) circle(4pt);\node[left] at (v3)   {\LARGE${v_2}$};
       \coordinate (v4) at (1, 1);\filldraw[black] (v4) circle(4pt);\node[right] at (v4)   {\LARGE${v_4}$};
       \coordinate (v5) at (-1, 2);\filldraw[black] (v5) circle(4pt);\node[left] at (v5)   {\LARGE${v_5}$};
       \coordinate (v6) at (1, 2);\filldraw[black] (v6) circle(4pt);\node[right] at (v6)   {\LARGE${v_8}$};
       \coordinate (v7) at (-1, 3);\filldraw[black] (v7) circle(4pt);\node[left] at (v7)   {\LARGE${v_1}$};
       \coordinate (v8) at (1, 3);\filldraw[black] (v8) circle(4pt);\node[right] at (v8)   {\LARGE${v_{3}}$};

       \coordinate (v9) at (0, 3.33);\filldraw[black] (v9) circle(3pt);\node[above] at (v9)   {\LARGE${v_{x}}$};
       \coordinate (v10) at (1, 2.5);\filldraw[black] (v10) circle(3pt);\node[right] at (v10)   {\LARGE${v_{y}}$};

     \draw [line width=1.5,black](v8) -- (v6) -- (v4) -- (v2) -- (v3) -- (v5) -- (v7);
     \draw [line width=1.5,black](v1) -- (v3);
     \draw [line width=1.5,black](v1) -- (v4);
     \draw [line width=1.5,black](v1) -- (v6);
     \draw [line width=1.5,black](v2) -- (v5);

     \draw [line width=1.5,black](v7) to[bend right] (v8);
     \draw [line width=1.5,black](v8) to[bend right] (v7);

     \draw [line width=1.5,black](v9) to[bend right] (v10);

       \filldraw[white] (-1, 0) circle(2pt);
       \filldraw[white] (1, 0) circle(2pt);
       \filldraw[white] (-1, 1) circle(2pt);
       \filldraw[white] (1, 1) circle(2pt);
       \filldraw[white] (-1, 2) circle(2pt);
       \filldraw[white] (1, 2) circle(2pt);
       \filldraw[white] (-1, 3) circle(2pt);
       \filldraw[white] (1, 3) circle(2pt);

       \filldraw[white] (0, 5) circle(2pt);
       \filldraw[white] (0, -1) circle(2pt);

\end{tikzpicture}
}

\centering (b2.2) $y_1= v_{n-2}$

\end{minipage}
\caption{$G\/$ has no $2\/$-cycle and $G_e\/$ has a $2\/$-cycle}    \label{Case(1) of labeling (1.2)}

\end{center}
\end{figure}

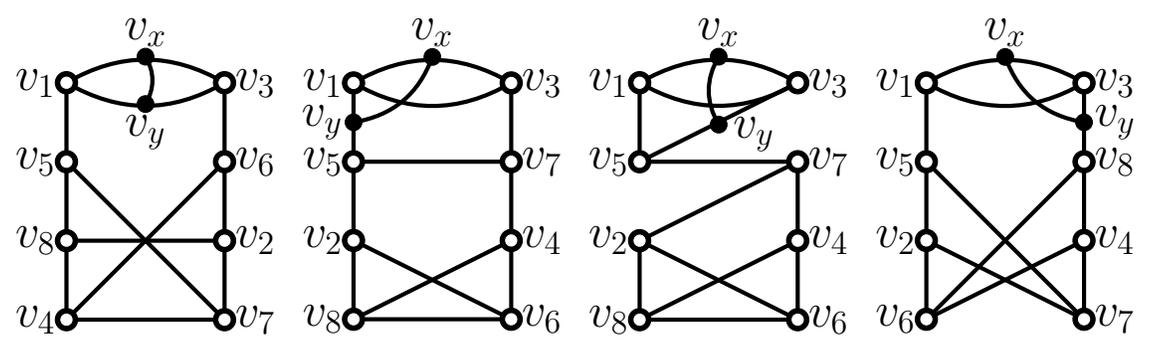
\end{document}